\documentclass[12pt, a4paper,reqno]{amsart}

\input{dibujos.tex}
\usepackage{centernot}
\usepackage{amsmath, amsthm, amssymb,amsfonts}
\usepackage{graphicx}
\usepackage{mathrsfs}
\usepackage{mathtools}
\usepackage{caption}
\usepackage{color}
\usepackage{lipsum,wasysym}
\usepackage[a4paper,margin=0.87in]{geometry}
\usepackage{geometry}
\usepackage{pst-all}
\usepackage{url}
\usepackage[colorlinks]{hyperref}
\usepackage[vcentermath]{youngtab}
\usepackage{young}
\usepackage{ytableau}
\usepackage{lineno}

\modulolinenumbers[1]

\definecolor{newred}{rgb}{1,0,0}

\numberwithin{equation}{section}

  \def\ag{{\mathfrak a}}  
  \def\bg{{\mathfrak b}}

\def\FM_pG{{\mathfrak F}}  \def\FM_pg{{\mathfrak f}}  \def\FM_pM{{\mathbb{F}}}
  
\def\HG{{\mathfrak H}}  \def\hg{{\mathfrak h}}

    \def\RM{{\mathbb{R}}}
\def\SG{{\mathfrak S}}  \def\sg{{\mathfrak s}}  
  \def\tg{{\mathfrak t}}  
  \def\ug{{\mathfrak u}}  
  \def\vg{{\mathfrak v}}

    \def\ZM{{\mathbb{Z}}}

    \def\AC{{\mathcal{A}}}

  \def\eb{{\mathbf e}}  
\def\FM_pB{{\mathbf F}}  \def\FM_pb{{\mathbf f}}  \def\FM_pC{{\mathcal{F}}}
    
\def\HC{{\mathbf H}}    \def\HC{{\mathcal{H}}}

    \def\RC{{\mathcal{R}}}

\def\UB{{\mathbf U}}

\def\FM_pS{{\EuScript F}}

\def\a{\alpha}

\newcommand{\mpar}{P^l(n)}
\newcommand{\ocmpar}{P^l_1(n)}
\newcommand{\ocmpardos}{P^2_1(n)}
\newcommand{\plam}{P(\lamn)}
\newcommand{\gblob}{B_{l,n}^p(\kappa)}
\newcommand{\gblobl}{B_{l,n}^{p,\lamn } (\kappa)}
\newcommand{\gblobdos}{B_{2,n}^p(\kappa)}
\newcommand{\gblobldos}{B_{2,n}^{p,\lamn } (\kappa)}

\newcommand{\bi}{\boldsymbol{i}}
\newcommand{\lau}{\mathcal{L}}

\newcommand{\lamn}{\boldsymbol{\lambda}}
\newcommand{\mun}{\boldsymbol{\mu}}
\newcommand{\nun}{\boldsymbol{\nu}}
\newcommand{\cellb}{\Delta^p(\lamn )}
\newcommand{\simb}{L^p({\lamn})}
\newcommand{\cellbl}{\Delta_{\lamn}^p(\mun)}
\newcommand{\simbl}{L_{\lamn}^p({\mun})}

\newcommand{\nc}{\newcommand} \newcommand{\renc}{\renewcommand}

\newcommand{\rdots}{\mathinner{ \mkern1mu\raise1pt\hbox{.}
    \mkern2mu\raise4pt\hbox{.}
    \mkern2mu\raise7pt\vbox{\kern7pt\hbox{.}}\mkern1mu}}

\def\res{{\mathrm{res}}}

\def\un{\underline}

\def\to{\rightarrow}
\def\FM_prom{\leftarrow}
\def\longto{\longrightarrow}

\nc{\triright}{\stackrel{[1]}{\to}}
\nc{\longtriright}{\stackrel{[1]}{\longto}}

\nc{\Hb}{H^\bullet}

\nc{\Br}{\mathcal{\color{blue}b}}
\nc{\HotRR}{{}_R\mathcal{K}_R}
\nc{\HotR}{\mathcal{K}_R}
\nc{\excise}[1]{}
\nc{\defect}{\text{df}}
\nc{\h}[1]{\underline{H}_{#1}}

\nc{\Ga}{\mathbb{G}_a} 
\nc{\Gm}{\mathbb{G}_m} 

\nc{\Perv}{{\mathbf{P}}}

\nc{\IH}{{\mathrm{IH}}}

\nc{\ic}{\mathbf{IC}}

\nc{\gl}{{\mathfrak{gl}}}
\renc{\sl}{{\mathfrak{sl}}}
\renc{\sp}{{\mathfrak{sp}}}

\renc{\Im}{\textrm{Im}}

\nc{\HCM}{H^{BM}}

 \DeclareMathOperator{\Hom}{Hom}

 \DeclareMathOperator{\ch}{ch}

\DeclareMathOperator{\End}{End} 

\DeclareMathOperator{\rank}{rk_v}

\numberwithin{equation}{section}

\newtheorem{teo}{Theorem}[section]
\newtheorem{lem}[teo]{Lemma}

\newtheorem{pro}[teo]{Proposition}
\newtheorem{defi}[teo]{Definition}
\newtheorem{cor}[teo]{Corollary}
\newtheorem{exa}[teo]{Example}
\newtheorem{remark}[teo]{Remark}

\newtheorem{algo}[teo]{Algorithm}
\newtheorem{notation}[teo]{Notation}

\newtheorem*{cona}{Blob vs Soergel Conjecture}
\newtheorem*{conb}{Blob vs Light Leaves Conjecture}
\newtheorem*{conc}{Blob vs Soergel Conjecture}
\newtheorem*{cond}{Categorical blob vs Soergel Conjecture}
\newtheorem*{consa}{Enhanced categorical blob vs Soergel Conjecture}

\newenvironment{dem}{\noindent \textit{Proof:}  \rm}{\quad \hfill $\square$}

\def\iff{\Leftrightarrow}

\nc{\simto}{\stackrel{\sim}{\to}}

\DeclareMathOperator{\Std}{Std}

\DeclareMathOperator{\gdim}{\dim_{v}}

\nc{\SD}{\mathcal{H}_{\mathrm{BS}}}

\newcommand{\ZZ}{\mathbb{Z}}

\begin{document}

\title[Blob algebra approach to modular representation theory]{Blob algebra approach to modular representation theory}

\author{Nicolas Libedinsky}
\address{Facultad de Ciencias, Universidad de Chile, Chile.}
\email{nlibedinsky@gmail.com}


\author{David Plaza}
\address{Instituto de Matem\'atica y F\'isica, Universidad de Talca, Chile.}
\email{dplaza@inst-mat.utalca.cl}






\begin{abstract}




Two decades ago P. Martin and D. Woodcock made a surprising and prophetic link between statistical mechanics and representation theory. They observed that the decomposition numbers of the blob algebra (that appeared in the context of transfer matrix algebras) are Kazhdan-Lusztig polynomials in type $\tilde{A}_1$. In this paper we take that observation far beyond its original scope. We conjecture that for $\tilde{A}_n$  there is an equivalence of categories between the characteristic $p$ diagrammatic Hecke category and a ``blob category'' that we introduce (using certain quotients of KLR algebras called \emph{generalized blob algebras}). Using alcove geometry we prove  the ``graded degree'' part of this equivalence for all $n$ and all prime numbers $p$. If our conjecture was verified, it would imply that the graded decomposition numbers of the generalized blob algebras in characteristic $p$ give the $p$-Kazhdan Lusztig polynomials in type $\tilde{A}_n$. We prove this for $\tilde{A}_1$, the only case where the $p$-Kazhdan Lusztig polynomials are known. 
\end{abstract}

\maketitle

\section{Introduction}

	\subsection{A new paradigm} 
	Kazhdan-Lusztig polynomials have been at the heart of representation theory since their discovery in $1979$. They have answered (often via geometric methods) an enormous number of deep questions concerning the characteristic zero and characteristic $p\gg 0$ representation theory of Weyl groups, Lie algebras, quantum groups, and reductive algebraic groups.
	 It was also widely believed (see, for example,  Lusztig and James's conjectures) that these polynomials control the modular (i.e., characteristic $p$) representation theory of these structures if $p$ is not too small. To be \emph{not too small} depends on the case and on the author. For example, in Lusztig conjecture for $SL_n(\overline{\mathbb{F}}_p),$  in Lusztig's original formulation the bound was  $p>2n-3$ and  the bound in Kato's version was $p>n$. In James's conjecture for $S_n$ the bound was $p>\sqrt{n}$.

A new paradigm has emerged in the last few years by the work of G. Williamson and his collaborators (see \cite{RW,JW,WExplosion,williamson2017algebraic,achar2017koszul}). Now we know that $p$-Kazhdan Lusztig polynomials are central objects of study in modular representation theory of Lie type objects. They provide solutions to old questions such as the modular representation theory of $SL_n(\mathbb{F}_q)$ (with $q$ a power of $p$) and give important information about modular representation theory of the symmetric group $S_n$. There is even a geometric interpretation of them!  Just as Kazhdan-Lusztig
polynomials have an interpretation in terms of the stalks of intersection cohomology complexes on flag varieties, $p$-Kazhdan-Lusztig
polynomials have one in terms of the stalks of parity sheaves with coefficients in a field of characteristic $p$.

 But there is still a big open problem in the theory. Kazhdan-Lusztig polynomials are obtained by calculating in an algebra (the Hecke algebra) but $p$-Kazhdan-Lusztig polynomials are obtained by calculating in a category (the Hecke category). This doesn't seem to be a satisfactory answer to the problem of how to calculate the $p$-Kazhdan-Lusztig polynomials. We need to understand  better these  polynomials in order to have a ``real'' answer\footnote{It is a matter of discussion what should one expect as a ``real" answer. Kazhdan-Lusztig polynomials are already quite difficult objects, and their $p$-version is harder. The authors hope for some kind of answer in the vein of Kazhdan-Lusztig polynomials, i.e., some recursive algorithm in some algebra or at least something more enlightening than calculating ranks of intersection forms.}. A step in this direction is given by the exciting conjecture (and theorem for very small $p$) by G. Lusztig and G. Williamson  \cite{lusztig2018billiards} where they predict that the characters of certain tilting modules for $SL_3$ (i.e., certain $p$-Kazhdan-Lusztig polynomials for $\tilde{A}_2$) are governed by a discrete dynamical system, that looks like billiards bouncing in equilateral triangles. But this conjecture also shows how incredibly complex  this quest is and how far we are from a full understanding of the $p$-canonical basis. 
 
  We see this paper as a link between physics and modular representation theory and we hope that it will help to raise intuitions from physics towards the (nowadays) obscure land of modular representation theory. The authors are working to find a statistical mechanical model for the generalized blob algebras in characteristic $p$ mimicking the situation of the usual blob algebras.

\subsection{Elias-Williamson's Hecke category}

The term $p$-canonical basis was coined by I. Grojnowski \cite{Gr} in 1999.  Its modern use (i.e., where one can actually calculate) comes from G. Williamson \cite{Wi12}. For this definition, we need first to explain a categorification of the Hecke algebra.

For any Coxeter system $(W,S)$, B. Elias and G. Williamson \cite{EW} defined the \emph{diagrammatic Hecke category} $\HC$ (or $\HC(W)$ when the group is not clear from the context). It is also called in the literature the Elias-Williamson's \emph{diagrammatic Soergel category}.

The definition goes roughly like this. Firstly we need to define the  \emph{Bott-Samelson diagrammatic category}  $\HC_{\mathrm{BS}}$ (for details see section \ref{real}). Assign a different color to each element of $S$. Objects in the category $\HC_{\mathrm{BS}}$ correspond to sequences of colored dots. For $\underline{w}= \red{s}{\blue{r}}\black{\cdots} \green{t}$ the corresponding object is 
$$ BS(\underline{w}):= {\color{red}{\bullet}}\, \blue{\bullet}\,  \black{\cdots} \,  \green{\bullet}\black{.}$$

Let $k$ be a field (we could even take $k$ to be the ring of integers in this definition and then extend scalars) and let $R:=k[\alpha_s]_{s\in S}$ be the polynomial ring in variables $\alpha_s$ parametrized by the set of simple reflections $S$.  The ring $R$ is endowed with an action of $W$ and with an operator $\partial_s:R\rightarrow R$ that we  make explicit in this section \ref{real}.

The morphisms in $\HC_{\mathrm{BS}}$ are (modulo some local relations) linear combinations over $R$ of isotopy classes of some decorated planar graphs embedded in the strip $\mathbb{R}\times [0,1]$. For example, a morphism between $\red{\bullet} \, \blue{\bullet}\, \red{\bullet}\, \green{\bullet}\, \blue{\bullet}\,\red{\bullet}$ (bottom) and $\green{\bullet}\, \blue{\bullet}\,\red{\bullet }$ (top) looks like
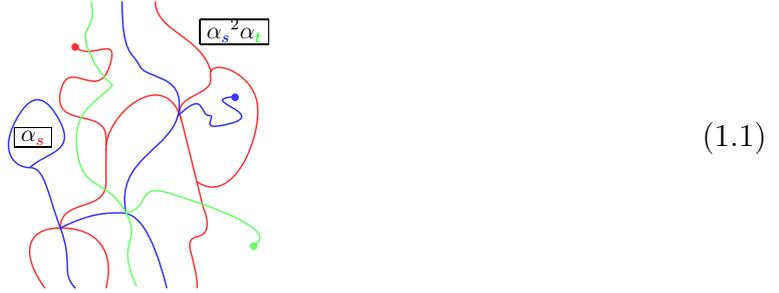
\begin{equation}
	\scalebox{.5}{\ExampleSoergel } 
\end{equation}

The edges of these graphs are colored as the elements of $S$ and they may end in a dot with the same color of the boundary of the strip. The connected components of the complement in $\mathbb{R}\times [0,1]$ of the embedded graph can be decorated by elements of $R$.

The generating morphisms, i.e., the kinds of vertices allowed in the graphs, are: 
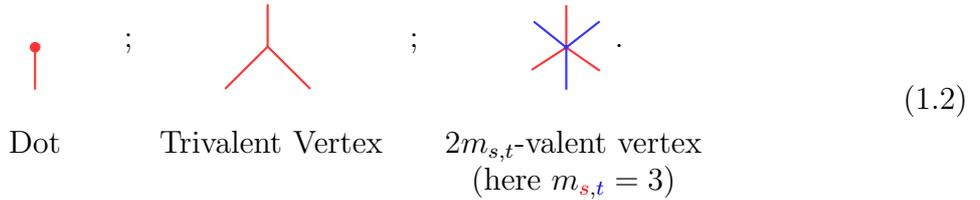
\begin{equation}
\begin{array}{ccccc}
\scalebox{0.7}{\NicoA} &  ; & \scalebox{0.7}{\NicoB} & ; &  \scalebox{0.7}{\NicoC}, \\
 &  &   &  &  \\
\mbox{Dot} &  &  \mbox{Trivalent Vertex} & &  2m_{st}\mbox{-valent vertex} \\
 &  &   &  & (\mbox{here } m_{{\red{s}}\blue{t}}=3) 
\end{array}
	\end{equation} 

 where $m_{st}$ denotes the order of $st$ in $W$. 
 
 Finally, to obtain the morphisms in $\HC_{\mathrm{BS}}$ one has to quotient the set of graphs obtained in this manner by the following local relations: 
\begin{itemize}
	\item One color relations: 
	\begin{equation}
	 \scalebox{0.6}{\NicoD}  = \scalebox{0.6}{\NicoE} \, \, ;  \, \, \scalebox{0.6}{\NicoF} =  \scalebox{0.6}{\NicoG} \, \, ; \, \,  \scalebox{0.6}{\NicoH} = 0\, \, ; \, \, \scalebox{0.7}{\NicoI} = \, \, \scalebox{0.7}{\NicoO} \, \, ;    \, \,  \scalebox{0.7}{\NicoE} \scalebox{0.7}{\NicoP} = \scalebox{0.7}{\NicoQ}  \scalebox{0.7}{\NicoE} + \, \, \scalebox{0.7}{\NicoR}\scalebox{0.7}{\NicoS}
\end{equation} 
	\item Two color relations (here we only illustrate  the case $m_{{\color{red}s}{\color{blue}t}}=3$): 
	\begin{equation}
 \scalebox{0.7}{\NicoJ}  = \scalebox{0.7}{\NicoK} \quad ; \qquad  \scalebox{0.7}{\NicoL} = \scalebox{0.7}{\NicoM}	+ \scalebox{0.7}{\NicoN}   
 \end{equation}
	\item Three color relations: We do not explain these relations, but instead refer the reader to \cite{EW}.  
\end{itemize}
To finish the definition of $\HC_{\mathrm{BS}}$, one needs to enrich the morphism spaces with a grading. They have a unique grading if we  prescribe that the dot, the trivalent and the $2m_{st}$-valent vertices have degrees $1, -1$ and $0$ respectively.

\begin{defi} \rm
 Let $\mbox{Kar}(\mathcal{H}_{BS})$ be the Karoubian envelope of $\mathcal{H}_{BS}$ (one formally adds direct summands of the objects, i.e. the objects in $\mbox{Kar}(\mathcal{H}_{BS})$ are pairs $(M, e)$, with $M $ an object of $\mathcal{H}_{BS}$ and $e$ an idempotent in $\End (M)$). The \emph{diagrammatic Hecke category} $\mathcal{H}$ is the additive graded envelope of $\mbox{Kar}(\mathcal{H}_{BS})$. That is, objects of $\mathcal{H}$ are formal direct sums of shifts of the objects in $\mbox{Kar}(\mathcal{H}_{BS})$, with obvious morphisms. Then $\mathcal{H}$ is a $k$-linear category.
\end{defi}

There are (at least) two other versions of this category. The first one is an algebraic version called the category of \emph{Soergel bimodules}. It is an additive monoidal full subcategory of the category of $\mathbb{Z}$-graded $R$-bimodules, where $R$ is some ring of polynomials acted upon by $W$.

The other  version is of geometric origin. Let $(W,S)$ be the Weyl group with simple reflections corresponding to a pair $(G,B)$, where $G$ is a complex reductive group  and  $B\subset G$ is a Borel subgroup. This version of the Hecke category is the additive, monoidal (under convolution) category
of $B$-biequivariant semi-simple complexes on $G$.

By explaining things like this, we are going backwards in  time. The geometric category was the first one that emerged. Soergel bimodules is a simplified version  with the advantages that 1) It is defined for all Coxeter groups and 2) It is way easier to compute with. The diagrammatic Hecke category is a simplified version of Soergel bimodules with the advantages that 1) It is well-behaved even over fields of positive characteristic and 2) It is defined over the integers, it is only after an extension of scalars that it is equivalent to Soergel bimodules (Versions of this category for particular groups were found before by \cite{LibRA}, \cite{EKh}, \cite{EDC}).


\subsection{The $p$-canonical basis}
It is a theorem by Elias and Williamson \cite{EW} (following a similar theorem by Soergel in the context of Soergel bimodules \cite{soergel2007kazhdan}) that the indecomposable objects $B_w$ of $\HC$ are parametrized (modulo grading shift) by the elements $w\in W$. These indecomposable objects are idempotents in $\mathrm{End}_{\HC_{BS}}BS(\underline{w})$ for any reduced expression $\underline{w}$ of $w$, so they are some (very difficult to calculate in practice) linear combination of colored graphs. Moreover, $B_w$ can be defined as the unique indecomposable object appearing in $BS(\underline{w})$ that has not appeared in any decomposition of $BS(\underline{x})$ for $x<w.$

We said before that the diagrammatic Hecke category behaves well even in positive characteristic. By that we meant that there is  a canonical ``character'' isomorphism\footnote{This character is also an isomorphism if one replaces $\HC$ by some categories of Soergel bimodules. It is the case for any Coxeter system  using the geometric representation in characteristic zero, as proved by Soergel \cite{SHC} and the first author \cite{Lieq}. It is also the case for Soergel bimodules associated to Weyl groups, using the Cartan matrix representation in positive characteristic $\neq 2, 3$, by \cite{LLC}. It is still unknown for a general Coxeter system if a ``reflection faithful representation'' exists in any characteristic. If one uses a representation that is not reflection faithful,   it is unknown if the character map is an isomorphim.} of $\mathbb{Z}[v,v^{-1}]-$algebras 
$\mathrm{ch}:[\HC]\rightarrow H$
between the split Grothendieck group of $\HC$ and the Hecke algebra $H$ of $(W,S)$ \cite{EW}. This is true for any characteristic of the field $k$.  Hence, $\HC$ provides a categorification of the Hecke algebra. In \cite{JW} it is proved that if $W$ is a crystallographic group, the characters of the indecomposable objects $B_x$   depend only on the characteristic of the field $k$. 

\begin{defi}\rm
	If $W$ is a crystallographic Coxeter group and $k$ has characteristic $p$, the $p$-\emph{canonical basis} is the subset  $\{\mathrm{ch}(B_x)\}_{x\in W}$ of the Hecke algebra. The coefficients $h^p_{x,y}$ of the base change matrix between
the $p$-canonical basis and the standard basis $$ \mathrm{ch}(B_y)=\sum_x h_{x,y}^pH_x$$ are the $p$-\emph{Kazhdan-Lusztig polynomials.}\end{defi}

\subsection{Some statistical mechanics}


The goal of classical statistical mechanics is to consider a model and by studying its behavior on some small scale (usually atomic) obtain some large-scale result (usually macroscopic). In this field of physics, there is a model, the Potts model, that in spite of not corresponding  to anything found in nature, is extremely interesting to study. This model consists roughly of a ``geometry" (where do the atoms sit) a ``configuration'' (an angular momentum for each atom) and a ``Hamiltonian'' (how much energy a given configuration has).

%
The Potts model is interesting mainly for two reasons. Firstly, because it gives accurately properties close to phase transitions (examples of phase transitions are boiling or melting water after a change of temperature, or, cooling enough some iron, it transforms itself into a permanent magnet, thus having \emph{spontaneous magnetization} below the Curie temperature). These properties  are (surprisingly) quite insensitive to the Hamiltonian, so even though the Hamiltonian of this model does not appear in nature, by calculating this model one can predict phase transitions properties that do appear in nature. Secondly, this model is interesting because it is extremely rich mathematically and in some particular cases (very few of them) one can find an analytical expression of the partition function, which is an extraordinary achievement in physics. 

Let us be (a bit) more precise. In the Potts model, you consider a graph, usually a subgraph of a Euclidean lattice. Let us consider for the rest of this section the  example of the graph in (\ref{Graph}) with vertices a set $I$ of atoms arranged in a square lattice of length $l$ and edges between \emph{neighbors}. Each vertex $i$ has some number of states  $s_i$ (a positive integer smaller than or equal to some constant $Q$) associated to  a ``spin". A \emph{configuration} $\mathbf{s}=\{s_i\}_{i\in I}$ is a prescription of a state for each atom, i.e., a function $\mathbf{s}: I \rightarrow \{1,2, \ldots, Q\}$.

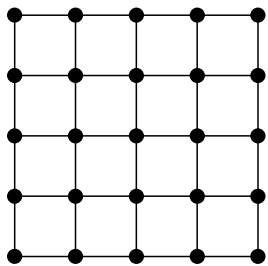
\begin{equation}\label{Graph}
	\Graph
\end{equation}

The partition function $Z: \mathbb{R} \rightarrow \mathbb{R}$ is a function of the parameter $\beta$ (the inverse of the thermal energy) and it defines the model. One magic feature of statistical mechanics is that if you can calculate the partition function, you are mostly done: you can calculate many important observables such as the free energy, the configuration with minimal energy where the system ``wants'' to go, etc. 

The Potts model tries to model ferromagnets (like iron). Due to physical considerations (exchange interactions tend to have short range) the partition function favor neighbors alignment. In formulas, we have 
$$Z(\beta) := \sum_{\substack{\mathbf{s} \\
\mathrm{configurations}}}  \prod_{\substack{(i,j) \\ \mathrm{neighbours}}}\mathrm{exp}(\beta \delta_{s_i,s_j}).$$
  The transfer matrix (if it exists) is a matrix $\tau$ satisfying the equality 
$$Z(\beta)=Tr(\tau^l).$$

The most interesting feature (see \cite{Martin4}) about writing the partition function in such a way is that all the physical information is condensed in the eigenvalues of $\mathcal{\tau}$.
In the particular case of  ``toroidal boundary conditions'', i.e. when $s_i=s_j$ if $i$ and $j$ are on the border of the graph and project to the same $x$ or $y$ component, then there is a transfer matrix, and it is  given explicitly by the $(Q\times Q)$-matrix with diagonal $\mathrm{exp}(\beta)$ and off-diagonal elements equal to $1$. This is a very relevant case for physics. 



It is in the aim of calculating  the transfer matrix (and ultimately, the partition function) that one introduces an algebra called the \emph{transfer matrix algebra} of the model  (that we will not define in this paper), defined by  ``easy" generators and some relations (we must say that there are other uses in physics for these algebras, such as identifying ``equivalent'' models). The transfer matrix turns out to be some specific element in the image of a representation of the transfer matrix algebra. 

\subsection{The blob algebra}

In the case of the Potts model on the graph in \eqref{Graph}, its transfer matrix algebra is the usual Temperley-Lieb algebra. If one considers the
Potts model but with a boundary that may have some extra  degrees of freedom (for details see \cite[Section 2.3]{gainutdinov2013physical}), its transfer matrix algebra is the ``one boundary Temperley Lieb algebra" or ``blob algebra".

The blob algebra was introduced by P. Martin and H. Saleur in  \cite{martin1994blob} as a two-parameter diagrammatic algebra. In the same paper they  constructed all the irreducible representations of that algebra when the ground field is of characteristic zero and the parameters are ``good enough'' (that is, when the algebra is semisimple). Some diagrammatic basis of the blob algebra turns out to be cellular in the sense of J. Graham and G. Lehrer \cite{graham1996cellular} in any characteristic and for any choice of the parameters. 

In \cite{martin2000structure} P. Martin and D. Woodcock completed the study of the representation theory of the blob algebra in characteristic zero  by treating the  ``bad parameters'' case. In particular they determined the decomposition numbers given by the cellular structure. In \cite{martin2003generalized} the same authors noticed that the  decomposition numbers of the blob algebra are given by KL polynomials (in Soergel's normalization) of type $\tilde{A}_1$ evaluated at $1$. This remark was upgraded by the second author and S. Ryom-Hansen  in \cite{plaza2013graded,plaza2014graded}. In these papers it is proven that the blob algebra is graded cellular (in the sense of J. Hu and A. Mathas \cite{hu2010graded}) and that its graded decomposition numbers (now Laurent polynomials) are the full KL polynomials, not just evaluated at $1$ (again, just for $\tilde{A}_1$). Over a field of positive characteristic the (ungraded) decomposition numbers  for the blob algebra were calculated by A. G. Cox, J. Graham and P. Martin in \cite{cox2003blob}.

\subsection{The generalized blob algebra}

Let $e,n,l>1$ be integers. Let $\kappa=(\kappa_1, \ldots, \kappa_l)$ be an $l$-tuple of elements of $\mathbb{Z}/e\mathbb{Z}$  ``separated from each other" in the sense that for every pair $i\neq j$, the element $\kappa_i$ does not belong to the set $\{\kappa_{j}-1,\kappa_j,\kappa_{j}+1\}$.

 Consider the algebra $A$ whose  elements are $\mathbb{F}_p$-linear combinations of isotopy classes of $n$-string diagrams (with dots or ``blobs") like the following one

 $$ \scalebox{.7}{\exampleKhL }  $$

  where in the bottom we place some $n$-tuple of elements of $\mathbb{Z}/e\mathbb{Z}$ that we call \emph{the bottom $n$-tuple of the diagram}. This determines in an obvious way \emph{the top $n$-tuple of the diagram}. Multiplication of two diagrams is defined by vertical concatenation if  the  top $n$-tuple of the lower diagram agrees with the bottom $n$-tuple of the top diagram. If not, the multiplication  is defined to be zero.

The \emph{generalized blob algebra} $B^p_{l,n}(\kappa)$
is the algebra $A$ modulo the following local relations
$$
\begin{array}{ll}
(1)\quad \CruceRelationB - \CruceRelationA = \CruceRelationBB - \CruceRelationAA = \delta_{i,j}\,\, \CruceRelationC \, \, ; & (2) \quad \scalebox{.8}{\ViernesA }\, =0.\\
(3)\quad \scalebox{.8}{ \ViernesC } \, =0 , \mbox{ if } i_1\neq \kappa_j \mbox{ for all }  j ; & (4) \quad\scalebox{.8}{\ViernesB } \, =0.
\end{array}
$$
and two other relations reminiscent from Reidemeister II and Reidemeister III (for details see section \ref{gba}).

The above definition presents $B^p_{l,n}(\kappa)$ as a quotient of a cyclotomic KLR algebra. This  does not correspond to the original definition of the generalized blob algebra as introduced by P. Martin and D. Woodcock in \cite{martin2003generalized}. Indeed,  $B^p_{l,n}(\kappa)$ was introduced as a   quotient of a cyclotomic Hecke algebra.  The fact that both definitions coincide is a consequence of Brundan and Kleshchev's isomorphism theorem \cite[Theorem 1.1]{BrKl}. A proof of the coincidence of both definitions is given in \cite[Section 12]{bowman2017family} and \cite[Section 7]{lobos2018graded}.  

A graded cellular basis for the generalized blob algebra was constructed independently by C. Bowman \cite{bowman2017many} and D. Lobos and S. Ryom-Hansen \cite{lobos2018graded}. Furthermore, C. Bowman showed in \cite{bowman2017many} that if $p=0$ then the graded decomposition numbers of $B^p_{l,n}(\kappa)$ 
  are given by KL polynomials in type $\tilde{A}_{l-1}$. 

Recently, C. Bowman and A. G. Cox \cite{bowman2017modular} proposed a conjecture that roughly speaking says  that if $p>0$ is \emph{not too small}\footnote{In this case the authors do not provide an explicit definition of what is meant by \emph{not to small}.} then graded decomposition numbers for $B^p_{l,n}(\kappa)$  coincide with KL polynomials in type $\tilde{A}_{l-1}$. 

\subsection{Alcove geometry}

 Let $W_l$ be the affine Weyl group of type $\tilde{A}_{l-1}$. We define a simply transitive action of $W_l$ on the set of alcoves (that are some subsets of  $E_l:=\mathbb{R}^l/\langle (1,1,\ldots, 1)\rangle$). Our action, with respect to the usual one, is dilated by a factor of $e$ and $\kappa$ prescribes the distance from the origin to the walls: there is a wall in the fundamental alcove ``distant"\footnote{We mean the geodesic distance in the graph with vertices $\mathbb{Z}^l$ and where two vertices are connected by an edge if they differ by some $\epsilon_i$. } $\kappa_j-\kappa_{j-1}$ from the origin, for each $j$ (for precise definitions see section \ref{alcove}).  We call the alcove containing the class of $0$ the \emph{fundamental alcove} $A_0$.

 We define $P_1^l(n)$ as the set of $l$-tuples of non-negative integers adding up to $n$. We can see an element  $\lamn\in P_1^l(n)$ as an element of $\mathbb{Z}^l\subset\mathbb{R}^l$ (whose coordinates add to $n$) and thus as an element in $E_l$. If the latter element belongs to the interior of an alcove, we call it \emph{regular}. We denote by $w_{\lamn}$ the unique element in $W_l$ such that $\lamn\in w_{\lamn}\cdot A_0$. For each $w\in W_l$ there is at least one (usually many) regular $\lamn$ such that $w=w_{\lamn}$.

If $\lamn \in  P_1^l(n)$ is regular, we will give  in Section \ref{alcove} a combinatorial definition (that we will omit in this introduction) of  $\underline{w}_{\lamn}$,
a particular reduced expression of $w_{\lamn}$.

 Another combinatorial definition given in Section \ref{alcove} that we omit in this introduction is, for a given generalized blob algebra and a regular $\lamn \in  P_1^l(n)$, an $n$-tuple $i^{\lamn}$ of elements of $\mathbb{Z}/e\mathbb{Z}.$

\subsection{The main conjecture and the main result}

 Let us fix the data for a generalized blob algebra, i.e., $(e,l,n, \kappa)$, and let us also fix $\lamn \in  P_1^l(n)$ regular.

  \begin{itemize}
  	\item Define $\HC_{BS}^l(\leq \lamn)$ as the full subcategory of $\mathbb{F}_p\otimes_R \HC_{BS}(W_l)$ having as  objects all  $BS(\underline{w}_{\mun})$ where $\mun$ is in the orbit of $\lamn$ and $w_{\mun}\leq w_{\lamn}$ in the Bruhat order. It is clear that this category has a finite number of objects. 
  	\item  Define $\mathrm{Blob}^l({\leq \lamn})$ as the category with objects all elements $\mun \in \ocmpar$ in the orbit of $\lamn$ such that $w_{\mun} \leq w_{\lamn }$. 
  Given two objects $\mun$ and $\nun$ in $\mathrm{Blob}^l({\leq \lamn})$,  the morphism space from $\mun$ to $\nun$ is the space $B^{p,\mun, \nun}_{l,n}$, defined as the subspace of $B^{p}_{l,n}$ with bottom $n$-tuple given by $i^{\mun}$ and top $n$-tuple given by $i^{\nun}$. The composition of morphisms is just multiplication in $B^{p}_{l,n}$. 
  \end{itemize}

\begin{cond}
For  each integer $l\geq 2$, there is an equivalence of categories $$\mathrm{Blob}^l({\leq \lamn})    \cong \HC_{BS}^l(\leq \lamn)$$ sending $\mun$ to $BS(\underline{w}_{\mun})$.
\end{cond}

We remark that, as both categories involved are skeletal (i.e. do not have different isomorphic objects), this equivalence of categories gives in fact an isomorphism of categories. 

The main result of this paper is the ``graded degree'' part of this conjecture.

\begin{teo}\label{maint} 

Let $n$ be a positive integer and $\lamn \in \ocmpar$  regular.  Let  $\mun$ and $\nun$ be objects of $\mathrm{Blob}^l({\leq \lamn})$. Then there is an isomorphism of 
 graded $\mathbb{F}_p$-vector spaces 
	\begin{equation}
\mathrm{Hom}_{\mathrm{Blob}^l({\leq \lamn})}(\mun, \nun)
\cong		 \mathrm{Hom}_{\HC_{BS}^l(\leq \lamn)}(BS(\underline{w}_{\mun}), BS(\underline{w}_{\nun })). 
	\end{equation}
\end{teo}  

\subsection{Lifting the conjecture}\label{subsection an equivalent conjecture}

Consider the Coxeter group $W_l$ as a poset with the Bruhat order. Consider a set  $\{\lamn_w\}_{w\in {W_l}} $ of elements $\lamn_w\in P_1^l(n_w)$ in the same orbit under the action of $W_l$, satisfying that  $w_{\lamn_w}:=w$ and that $n_{u}<n_v$ if $u<v$. 

For each pair $u,v\in W_l$ such that $u<v$, there is a fully-faithful functor 
$B_{u}^v: \mathrm{Blob}^l({\leq \lamn_{u}})\rightarrow \mathrm{Blob}^l({\leq \lamn_{v}})$ which maps an object $\mun $ of  $\mathrm{Blob}^l({\leq \lamn_{u}})$ (and therefore in $P_1^l(n_u)$) to the object $\mun'$ of  $\mathrm{Blob}^l({\leq \lamn_{v}})$ (and therefore in $P_1^l(n_{v})$) which is obtained from $\mun$ by adding $(n_{v}-n_u)/l$ to each component. The description of the action of  $B_u^v$ on Hom-spaces is given in terms of the graded cellular basis of   $B_{l,n}^p$ introduced in section \ref{sec gcb for blob}. For this reason we postpone such a description for  section \ref{Section Idempotent truncations}. This gives a direct system of categories (so, in particular there is an equality of functors $B_{u}^v\circ B_{v}^w= B_{u}^w$ for all $u\leq v \leq w$), which clearly has a limit. To see this, just remark that after forgetting the Hom's one has a direct system of sets (the objects) which is known to have a limit, and for the Hom spaces, after fixing two objects, it reduces to the existence of a direct limit of vector spaces, which is also known to exist. Let us denote this direct limit by 
\begin{equation}\label{superconj}
	\mathrm{Blob}^l({\infty}):= \varinjlim \mathrm{Blob}^l({\leq \lamn_w}).
\end{equation}

On the other hand, one has a similar, but easier direct limit for the Hecke category, given by fully faithful embeddings $H_u^v: \HC_{BS}^l(\leq \lamn_u) \rightarrow \HC_{BS}^l(\leq \lamn_v)$. In this case we can give a direct description of the limit: the full subcategory of $\mathbb{F}_p\otimes_R \HC_{BS}(W_l)$ with objects the $\underline{w}_{\lamn_w}.$ We denote this direct limit by$$\HC_{BS}^{l,\infty}=\varinjlim \HC_{BS}^l(\leq \lamn_w).$$ 
We conjecture that one can lift  Categorical blob vs Soergel conjecture (Cbvs for short) to an equivalence of categories  
 \begin{equation}\label{bigconj}\mathrm{Blob}^l({\infty})        \cong  \HC_{BS}^{l,\infty}.
\end{equation}
 We believe that to prove this conjecture one should prove that  the Cbvs conjecture is a natural equivalence. Rephrasing, if we call $F_{\lamn}:\HC_{BS}^l(\leq \lamn)\rightarrow \mathrm{Blob}^l({\leq \lamn})$ the Cbvs equivalence, the naturality means that for each $u<v$ one has an equality of functors $B^v_u\circ F_{\lamn_u} = F_{\lamn_v}\circ H_u^v$. This would imply the equivalence \eqref{bigconj}, just by applying the universal property of direct limits and seeing that the functor produced is essentially surjective and fully faithful. In general contexts, it is a rare thing to find an equality of functors, but in this case there is some evidence in type $\tilde{A}_1$ pointing  out in that direction by unpublished work of D. Lobos, S. Ryom-Hansen and the second author (given the diagrammatic nature of the Hom spaces).

One could also think of other ways of proving this conjecture by relaxing that above equality of functors. For example, adding self-equivalences $h_{\lamn_w}:\HC_{BS}^l(\leq \lamn_w)\rightarrow \HC_{BS}^l(\leq \lamn_w)$ and $b_{\lamn_w}:\mathrm{Blob}^l({\leq \lamn_w})\rightarrow \mathrm{Blob}^l({\leq \lamn_w})$ such that the equality above is relaxed to: 
\begin{equation}
	B^v_u\circ F_{\lamn_u}\circ h_{\lamn_u} = b_{\lamn_v}\circ F_{\lamn_v}\circ H_u^v .
\end{equation}
If $F_{\lamn_u}\circ h_{\lamn_u}=b_{\lamn_u}\circ F_{\lamn_u}$ then it is again easy to see that the equivalence Cbvs lifts to the equivalence \eqref{bigconj}. There are still other alternatives, such as a natural transformation from the functor $B^v_u\circ F_{\lamn_u}$ to the functor $F_{\lamn_v}\circ H_u^v$ for each $u\leq v$ satisfying some higher congruence relations, but we do not want to enter into details about this approach. 
In any case, if one proves \eqref{bigconj}, the equivalence is preserved if one takes the Karoubi envelope (being this a categorical construction) 
\begin{equation}\label{4}  \mathrm{Blob}^l({\infty})^e\cong (\HC_{BS}^{l,\infty})^e.\end{equation}
 This latter equivalence is preserved under taking the additive closure (objects in the additive closure are formal finite direct sums of objects in the original category)
\begin{equation}\label{5}  \mathrm{Blob}^l({\infty})_{\oplus}^e\cong (\HC_{BS}^{l,\infty})_{\oplus}^e.\end{equation}
By the description of $(\HC_{BS}^{l, \infty})$ given before, one has an equivalence of categories: 
\begin{equation}\label{poto}
	(\HC_{BS}^{l,\infty})^e_{\oplus}\cong(\mathbb{F}_p\otimes_R \HC_{BS}(W_l))^e_{\oplus}.
\end{equation}

Let us define the \emph{Blob category} $\mathrm{Blob}^l$ by the equality $$\mathrm{Blob}^l:= \mathrm{Blob}^l({\infty})^e_{\oplus}.$$

By the preceding paragraph, if one proves \eqref{bigconj}, it implies

\begin{consa}
For  each  integer $l>1$, there is an equivalence of categories	$$ \mathrm{Blob}^l\cong \mathbb{F}_p\otimes_R \HC(W_l). $$ 
\end{consa}

\subsection{About this conjecture}

When a conjecture is proposed there are two basic questions that one should try to answer. Namely, why should one care and why is it reasonable to believe in it. 

With respect to the first question we have the following reasons in increasing order of importance. 
\begin{enumerate}
	\item This would give a new Hecke category and it is widely known that each Hecke category (graded BGG category $\mathcal{O}$, Soergel bimodules, Elias-Williamson diagrammatic category, 2-braid groups, perverse (or parity) sheaves on affine Grassmannians, sheaves on moment graphs, etc.) has been very important and has brought deep new insights into the theory.
	\item In the blob category it is self-evident what the anti-spherical category is (see \cite{libedinsky2017anti}), and how to calculate on it, in contrast with other Hecke categories. The anti-spherical category controls the representation theory of algebraic groups as conjectured (and proved in type A) by Riche and Williamson \cite{RW} and proved in any type by  Achar, Makisumi, Riche and Williamson in the preprint \cite{achar2017koszul}.
	\item If one confirms this conjecture, one could calculate the $p$-canonical basis using the language  of generalized blob algebras (we explain this in more detail later in Section \ref{section algorithm to compute gdn blob}). As each of these algebras is  some KLR algebra modulo two relations, this would translate the fundamental problem of understanding the $p$-canonical basis into  a language that the KLR community can work with.
	\item In the same vein as the last point, this equivalence would be related to, but it would be much more concrete computationally than  the ``Categorical Schur-Weyl duality'' (which is implicit in \cite{rouquier20082} and a version of it is  explained in \cite[Theorem 8.1]{RW}). 
	\item Calculations that are quite short using Soergel bimodules are divided into several steps in the blob category. This is reminiscent to the category of singular Soergel bimodules studied by G. Williamson in his thesis. While  calculations are longer in the blob algebra, they are more transparent.  
	\item The last is our most important reason. We believe that the generalized blob algebras are transfer matrix algebras of some statistical mechanics system, as blob algebras are. If this was true, then one could bring physical intuitions into the problem of calculating the $p$-canonical basis. This problem doesn't seem to be near to its end, so this kind of intuitions might be extremely fruitful.  
\end{enumerate}

For the second question,  our main answer  is of course Theorem \ref{maint}, but we will explain in the next section a decategorified version of the conjecture. The reasons to believe in that simplified version also give us evidence for the categorical version.

\subsection{Decategorified conjecture}

As we said before, generalized  blob algebras are   graded cellular. Recall that $P_1^l(n)$ is the set of $l$-tuples of non-negative integers adding up to $n$.  Graded cellularity of $B^p_{l,n}(\kappa)$ gives, by the general theory of graded cellular algebras  and some extra effort, that $B^p_{l,n}(\kappa)$ is equipped with  ``graded cell modules'' $\Delta^p(\lamn)$ and ``graded simple modules'' $L^p(\lamn)$, for each $\lamn \in P_1^l(n)$\footnote{In fact, the general theory tells us this only for elements $\lamn$ in a subset of $P_1^l(n)$, but in this case one can prove that this subset is $P_1^l(n)$ itself.}. The \emph{graded decomposition numbers } $d_{\lamn, \mun}^p(v)$ (which despite their name, are not numbers, but polynomials in $\mathbb{Z}[v,v^{-1}]$) give graded multiplicities  of simple modules in composition series of cell modules. More precisely, one can define them  by the formula in the Grothendieck group

$$ [\Delta^p(\lamn)]=\sum_{\mun \in P_1^l(n)}  d_{\lamn, \mun}^p(v)[L^p(\mathbf{\mun})].$$

Now we can present the decategorified version of our conjecture.

  \begin{conc}\label{central}
 	Let 
 	$\lamn, \mun\in P_1^l(n)$ be regular, with  $\mun$ in the orbit of $\lamn$, then
  	 $$h^p_{w_{\lamn},w_{\mun}}= d^p_{\lamn, \mun}\in \mathbb{Z}[v,v^{-1}].$$
  	\end{conc}

Martin and Woodcock's main motivation to introduce generalized blob algebras was to try to find  algebras in which  the polynomials occurring in Soergel's tilting algorithm \cite{SoeKL} had a representation theoretic meaning. In particular, they were looking for algebras whose decomposition numbers were given by Kazhdan-Lusztig polynomials of type $\tilde{A}_{l-1}$ evaluated at $1$. Using our notation, their hope can be expressed as 
\begin{equation}\label{MWs conjecture}
h^0_{w_{\lamn},w_{\mun}}(1)= d^0_{\lamn, \mun}(1).	
\end{equation}

 Although this was never formally documented, it is sometimes referred as \emph{Martin and Woodcock's Conjecture}. 
In this setting, the Blob vs Soergel Conjecture can be seen as a ``graded'' generalization (for any prime number $p$) of \eqref{MWs conjecture}.

\begin{remark} \rm
	As we already pointed out, generalized blob algebras are quotients of cyclotomic KLR algebras. On the other hand, cyclotomic KLR algebras are isomorphic to cyclotomic Hecke algebras \cite{BrKl}. In a recent paper Elias and Losev \cite{elias2017modular} showed that the decomposition numbers of cyclotomic Hecke algebras can be obtained as evaluation at $1$ of certain parabolic $p$-Kazhdan-Lusztig polynomials. It would be tempting to believe that the previous conjecture could be deduced from Elias and Losev's work. However, this is not the case since in general cell (standard) modules for cyclotomic Hecke algebras are not the pullback of cell modules for generalized blob algebras (see \cite{ryom2010ariki}), since the cellular structure that they use (the one defined in \cite{dipper1998cyclotomic}) is not compatible with the one used in the generalized blob algebra.
\end{remark}

As we have said before, the $p$-canonical basis is very hard to calculate. For the symmetric group we can calculate it using computers until $S_{10}$ but for  type $\tilde{A}_n$ there is only a formula for $n=1$ (and the aforementioned conjecture by Lusztig and Williamson for some finite family of elements in $\tilde{A}_2$). The full complexity of the case $\tilde{A}_1$ (where the $p$-canonical basis relates to the canonical basis in a fractal-like way, see \cite{JW}) can be seen in the paper \cite{ElLi}, where enormous amount of effort is needed (including the development of the theory of $n$-colored Temperley-Lieb $2$-categories) to prove that the projectors defined there are indeed indecomposable when some quantum numbers are not zero. All this to say that, even though in usual representation theory the case $\tilde{A}_1$
is relatively easy, in this context that is far from being true. 

\begin{teo}\label{uno}
Blob vs Soergel Conjecture is verified in the case $l=2$ if $p\neq 2$. 	
\end{teo}

\subsection{Organization of the paper}
In section \ref{Section dos} we introduce the Hecke algebra and the Hecke category as well as intersection forms and the $p$-canonical basis (using diagrams). In section \ref{gba} we introduce generalized blob algebras and explain their  graded cellular algebra structure. We also explain the combinatorial way to associate $\lamn \rightsquigarrow \tg^{\lamn} \rightsquigarrow \boldsymbol{i}^{\lamn}$ and introduce the algebra $\gblobl $. In section \ref{alcove} we  study some  alcove geometry in order to find the map $\lamn \rightarrow \underline{w}_{\lamn}$ and prove Theorem \ref{maint}.
Finally, in section \ref{proof} we prove Theorem \ref{uno}.

\subsection{Acknowledgments: } We would like to warmly thank Felipe Torres, Paul Martin and Herbert Saleur for their explanations about the physics behind the blob algebra. We would also like to thank Ben Elias, Ivan Losev, Stephen Griffeth, Geordie Williamson, Francesco Brenti, Paolo Sentinelli and Gaston Burrull for their useful corrections and/or comments on earlier versions of the paper.

The  first author was partially supported by FONDECYT project 1160152 and Anillo project ACT1415. The second author was partially supported by  FONDECYT project 11160154 and Inserci\'on en la Academia-CONICYT project 79150016.

\section{The $p$-canonical basis via Soergel calculus}\label{Section dos}

\subsection{The Hecke algebra}  We follow Soergel's conventions for the Hecke algebra and Kazhdan-Lusztig basis. Let $(W,S)$ be a Coxeter system and $(m_{sr})_{s,r \in S}$ its Coxeter matrix. Let $l:W\rightarrow \mathbb{N}$ be the corresponding length function, $e\in W$ the identity element  and $\leq$ the Bruhat order on $W$. We use the convention that if $\underline{w}$ is a reduced expression, then $w$ is the corresponding element in $W$.

Consider the ring $\mathcal{L}=\mathbb{Z}[v^{\pm1}]$  of Laurent polynomials with integer coefficients in one variable $v$. 
The \emph{Hecke algebra} $\mathcal{H}=\mathcal{H}(W,S)$ of a Coxeter system $(W,S)$ is the associative algebra over $\mathcal{L}$ with generators $\{H_s\}_{s\in S}$, quadratic relations $(H_s+v)(H_s-v^{-1})=0$ for all $s\in S$, and for every pair $s, r \in S$ with $m_{sr}<\infty$,  braid relations $$H_sH_rH_s\cdots=H_rH_sH_r\cdots$$ with $m_{sr}$ terms on each side of the equation.

Consider $x\in W.$ To any $sr\cdots t$ reduced expression  of $x$ one can
associate the element $H_sH_r\cdots H_t\in \mathcal{H}.$   
H. Matsumoto proved that this element, that we call $H_x,$ is independent of the choice of reduced
expression of $x$. N. Iwahori proved that $$\mathcal{H}=\bigoplus_{x\in W}\mathcal{L} H_x.$$

\subsection{Kazhdan-Lusztig polynomials}

Let us define for $s\in S$, the element $\underline{H}_s:=H_s+v.$ There is a unique ring homomorphism (moreover, an involution) $H \mapsto \overline{H}$ on $H$ such that $\overline{v}=v^{-1}$ and $\overline{H_x}=(H_{x^{-1}})^{-1}$.  We will call an element
\emph{self-dual} if it is invariant under $\overline{(-)}$.

The fact that 
 for every element $x\in W$ there is a unique self-dual element $\underline{H}_x\in {H}, $  such that  $\underline{H}_x\in H_x+\sum_{y\in W}v\mathbb{Z}[v]H_y,$ was proved by D. Kazhdan and G. Lusztig in \cite{KLPolynomials}.   We call the set $\{ \underline{H}_x\}_{x\in W}$  the \emph{ Kazhdan-Lusztig basis} or the \emph{canonical basis} of  $\mathcal{H}$. It is a basis of the Hecke algebra as an $\mathcal{L}$-module. In formulas, $$\mathcal{H}=\bigoplus_{x\in W}\mathcal{L} \underline{H}_x.$$  For each couple of elements $x,y\in W$ we define $h_{y,x}\in \mathcal{L}$ by the formula $$\underline{H}_x=\sum_yh_{y,x}H_y.$$

\begin{remark}
The \emph{Kazhdan-Lusztig polynomials} (as defined in \cite{KLPolynomials}) are given by the formula $$P_{y,x}=(v^{l(y)}-v^{l(x)})h_{y,x},$$ and they are polynomials in $q:=v^{-2}$.
\end{remark}

\begin{defi}
Let  $\underline{w}=sr\cdots t$ be an expression of $w\in W$. Then we define      $$\underline{H}_{\underline{w}}:=\underline{H}_s\underline{H}_t\cdots \underline{H}_t.$$
\end{defi}

For Weyl groups or affine Weyl groups and for a given prime number $p$, there is another family of bases of the Hecke algebra, the so-called $p$-canonical basis $\underline{H}_x^p$  $$\mathcal{H}=\bigoplus_{x\in W}\mathcal{L} \underline{H}^p_x.$$

To this date there is no known way to calculate the $p$-canonical basis entirely within the Hecke algebra; all the known ways rely on categorical calculations. We will define the diagrammatic Hecke Category of Elias and Williamson, that, for many purposes (in particular to prove our main theorem \ref{uno} in Section \ref{proof})  is the best way to calculate the $p$-canonical basis. From now on we will concentrate on the case of affine Weyl groups of type $A$; we consider for the rest of the paper $W=\tilde{A_n}$.

\subsection{The Hecke category $\HC$} \label{real}

Recall that  a realization, as defined in \cite[\S 3.1]{EW})  consists of a commutative ring $\Bbbk$ and a
free and finitely generated $\Bbbk$-module $\hg$ together with subsets
\[
\{ \a_s \}_{s \in S} \subset \hg^* \quad \text{and} \quad \{ \a^\vee_s \}_{s \in S} \subset \hg
\]
of ``roots'' and ``coroots'' such that $\langle \alpha_s,
\alpha_s^\vee \rangle = 2$ for all $s \in S$ and such that
  the formulas 
\[
s(v) := v - \langle \alpha_s , v \rangle \alpha^\vee_s \hspace{.5cm} \mathrm{for }\  s \in S\ \mathrm{and }\ 
v \in \hg,
\]
define an action of $W$ on $\hg$. A realization also has to satisfy some other technical conditions on $2$-colored quantum numbers. To define Elias-Williamson Hecke category one needs two   conditions, namely ``Demazure surjectivity'' and ``balancedness''. 

Unless otherwise stated we will assume in this paper that $\Bbbk = \mathbb{F}_p$ and that $\hg$ is the Cartan matrix representation of $W$, 
i.e.,  $\hg = \bigoplus_{s \in S} \Bbbk \alpha_s^{\vee}$ and the elements $\{ \alpha_s \} \subset \hg^*$ are defined by the equations
\begin{equation} \label{eq:justcos}
\langle \alpha_t^\vee, \alpha_s \rangle = -2\cos(\pi/m_{st})
\end{equation}
(by convention $m_{ss} =1$ and $\pi/\infty = 0$).  This data satisfies all the technical conditions that we have mentioned to define Elias-Williamson Hecke category, so this is the realization that we will use in this paper. 

Let  $R =S(\hg^*)$ be the ring of regular functions on $\hg$ or, equivalently, the symmetric algebra of $\hg^*$ over
$\Bbbk$. We view $R$ as a graded $\Bbbk$-algebra by declaring  $\deg \hg^*  =
2$. The action of  $W$ on $\hg^*$, extends to an action on $R$, by functoriality.
For $s \in S$, let  $\partial_s :
R \to R[-2]$ be the \emph{Demazure operator} (where the shift $[-2]$ indicates that $\partial_s$ is a degree $-2$ map) defined by 
\[
\partial_s(f) = \frac{f - sf}{\alpha_s}.
\]
In \cite[\S
3.3]{EW} the authors prove that this operator is well defined under our assumptions.

In the introduction we gave a rough definition of the \emph{Bott-Samelson Hecke category} $\HC_{\mathrm{BS}}$. We refer the reader to \cite{EW} for details on that construction.

If $M=\bigoplus_iM^i$ is a $\mathbb{Z}$-graded object, we denote by $M(1)$
its grading shift defined by $M(1)^i=M^{i+1}.$ For an additive category
$\AC$ we denote by $[\AC]$ its split Grothendieck group.
If in addition $\AC$ has Hom spaces enriched in graded vector
spaces we denote by $\AC^\oplus$ its additive graded closure, i.e.  objects are formal finite direct sums $\bigoplus a_i(m_i)$ for certain objects
$a_i \in \AC$ and ``grading shifts'' $m_i \in \ZM$ and \[
\Hom_{\AC^\oplus}( \bigoplus a_i(m_i), \bigoplus b_j(n_j)) := \bigoplus
\Hom(a_i, b_j)(n_j - m_i).
\]
The category $\AC^\oplus$ 
is equipped with a grading shift
functor defined on objects by $\bigoplus a_i(m_i) \mapsto
\bigoplus a_i(m_i+1)$.  We define
$\AC^e$ to be the Karoubian envelope of 
$\AC^{\oplus}$. 
We can finally define the  \emph{Hecke category}  $\HC:=\HC_{\mathrm{BS}}^e.$ 



\subsection{Light leaves}\label{basic}
 Let $\un{x} = s_1s_2\dots s_m$
 be an expression (a sequence of elements in $S$). A
 \emph{subsequence} of  $\un{x} = s_1s_2\dots s_m$ is a sequence $\pi_1 \pi_2 \dots \pi_m$ such that $\pi_i \in \{ e, s_i \}$ for all $1 \le i \le m$. Instead of working with subsequences,
we work with the equivalent datum of a sequence $\eb = \eb_1 \eb_2 \dots \eb_m$ of 1's and 0's giving the indicator function of a subsequence, which we refer to as a $01$\emph{-sequence}.
Associated to this, one has the  \emph{Bruhat stroll}. It is the sequence $x_0, x_1,
\dots, x_m$ defined inductively by $x_0=e$ and  \[ x_i := s_1^{\eb_1} s_2^{\eb_2} \dots s_i^{\eb_i} \] for $1 \le i \le m$. We call $x_i$ the $i^{\mathrm{th}}$-\emph{point} and $x_m$ the \emph{end-point} of the Bruhat stroll. We denote the end-point by $\un{x}^{\eb}$. 



\emph{ Light leaves}  and \emph{Double leaves} for Soergel bimodules
were introduced  in \cite{LLL} and \cite{LLC}. They give bases, as $R$-modules of the Hom spaces between
Bott-Samelson bimodules. We suggest reading \cite[\S~6.4--6.5]{Li17}  in order to get used  to these combinatorial
objects (in that paper the dot is called $m_s$, the trivalent vertex is called $j_s$ and the four and six valent vertices are called $f_{sr}$). We also suggest reading   \cite[\S~6.1--6.3]{EW}, where these bases are explained in the diagrammatic language used in this paper.
 
We denote by $\mathbb{L}_{\underline{w}}(x)$ (resp. $\mathbb{L}^d_{\underline{w}}(x)$) the set of light leaves (resp. light leaves of degree $d$) with source $\underline{w}$ and target  some expression of $x\in W$. In particular, the cardinality of $\mathbb{L}_{\underline{w}}(x)$ is equal to the number of $01$-sequences with end-point equal to $x$. 

Let $X$ be a set of homogeneous elements of graded vector spaces. We define the
\emph{degree} of $X$ by the formula
$$d(X)=\sum_{x\in X}v^{\mathrm{deg}(x)}.$$

In \cite[\S~5.2]{LLC} the following  is proved
\begin{equation} \label{eq uno}
	\underline{H}_{\underline{w}} =\sum_{x\leq w} d(\mathbb{L}_{\underline{w}}(x)) H_x.
\end{equation}
  
Let $(-)^a:\mathcal{H}\rightarrow \mathcal{H}$ be the endofunctor that flips upside down a diagram. Let $\underline{w}$ and $\underline{u}$ be some expressions. Then, the set

\begin{equation} \label{equation double leaves basis}
	\{ \mathbb{L}_{\underline{u}}(x)^a \circ \mathbb{L}_{\underline{w}}(x) \, |\, x\in W  \}
\end{equation}

is an $R$-basis of $\mbox{Hom}_{\mathcal{H}} (BS(\underline{w}) , BS (\underline{u}))$. It is called the \emph{double leaves basis}.

\subsection{Intersection forms}
In this section  $\Bbbk=\mathbb{Z}$ and  $d$ is any integer. Let  $n$ and $m$ be integers such that  $\mathbb{L}^d_{\underline{w}}(x)=\{l^{+}_1, \ldots, l^{+}_{n}\}$ and $\mathbb{L}^{-d}_{\underline{w}}(x)=\{l^{-}_{1}, \ldots, l^-_m\}$. Expand the element $l^-_j\circ (l^+_i)^a\in \mathrm{End}(\underline{x})$ in the double leaves basis:
 $$l^{-}_{j}\circ (l^+_i)^a=n_{ij}\mathrm{id}+\sum_k ll_k^a\cdot ll'_k,$$
	where $ll_k, ll'_k$ belong to some $\mathbb{L}_{\underline{x}}(y_k)$ with $y_k<x.$ Define the \emph{$d$-th grading piece of the intersection form} as the matrix $I_{\underline{w},x,d}=(n_{ij})_{i,j}$. Finally, we define the intersection form simply as the matrix direct sum of the degree pieces $$I_{\underline{w},x}:=\bigoplus_dI_{\underline{w},x,d}.$$
	
	
	If $p$ is a prime, the matrix $I_{\underline{w},x}^p$  is the reduction of the matrix $I_{\underline{w},x}$ modulo $p$. In \cite{JW}  the following formula is proven
	\begin{equation}\label{eq dos}
	\underline{H}_{\underline{w}} =  \sum_{y\leq w} \mbox{rk}_v \left( I_{\underline{w},y}^p\right)\underline{H}_y^p .  
\end{equation}

In \cite[Definition 6.24]{EW} the authors define a character map $\ch:
[\HC] \to H$ and in  \cite[Corollary 6.27]{EW} they prove that it is an isomorphism sending $[B_s]$ to $\underline{H}_s$ and $[v]$ to the empty word shifted by $1$. This is the reason to call  $\HC$ the
Hecke category.


\subsection{$p$-canonical basis}

Following Soergel's
classification of indecomposable Soergel bimodules \cite{soergel2007kazhdan}, Elias and Williamson  proved in   \cite{EW} that the indecomposable objects in $\HC$
are indexed by $W$ modulo shift, and they call $B_w$ the
indecomposable object corresponding to $w\in W$.  It happens that the object $B_s$ is the sequence with one element $(s)\in \HC.$ Because of this, if $\underline{w}=(s,r,\cdots, t)$ we will sometimes denote by $BS({\underline{w}}):=B_sB_r\cdots B_t$  the element $\underline{w}\in \HC$. The right hand side is the monoidal product in the (strict) monoidal category $\HC$.

 In  \cite{EW2} the authors prove  that, if $\Bbbk$ is the field $\mathbb{R}$, then $\ch([B_w]) = \underline{H}_w$. (More
 precisely, to obtain this result one must combine the equivalence
 between $\HC$ and the category of Soergel bimodules proved in \cite{EW} with the main results of \cite{EW2} and
 \cite{Lieq}.) Thus the indecomposable objects in $\HC$ categorify the
 Kazhdan-Lusztig basis when we work over the real numbers. But when  $\Bbbk$ is the field $\mathbb{F}_p$, then the character of $B_w$ is no longer forced to be the canonical basis.

\begin{notation} Let $\Bbbk=\mathbb{F}_p$. To emphasize the field in which we are working, we will call $B_w^p$ the indecomposable objects of $\mathcal{H}$ (instead of the usual $B_w$). We use the notation $\ch([B_w^p]) = \underline{H}_w^p\in \mathcal{H}.$  The \emph{$p$-canonical basis} is the set $\{\underline{H}_w^p\}_{w\in W}$, which can be proved to be a basis of the Hecke algebra. If we write the $p$-canonical basis in terms of the standard basis
	\begin{equation}\label{eq tres} \underline{H}_y^p =\sum_{x\leq y} h_{x,y}^pH_x, 
	\end{equation}
the polynomials $ h_{x,y}^p$  are called \emph{$p$-Kazhdan-Lusztig polynomials.}
 \end{notation}
 
\subsection{An important formula} \label{section an important formula}

By (\ref{eq uno}),  (\ref{eq dos}) and  (\ref{eq tres})  we obtain
\begin{eqnarray}\label{eq cuatro}
\sum_{x\leq w} d(\mathbb{L}_{\underline{w}}(x)) H_x &=&  \sum_{y\leq w} \mbox{rk}_v\left( I_{\underline{w},y}^p\right)\sum_{x\leq y} h_{x,y}^pH_x  \\
&=&\sum_{x\leq w} \left( \sum_{x\leq y\leq w} \mbox{rk}_v \left( I_{\underline{w},y}^p\right) h_{x,y}^p \right) H_x.\nonumber
\end{eqnarray} 
	
By equating the coefficients we obtain, for each $x\leq w$ the equation 

\begin{equation}\label{eq cinco}
	d(\mathbb{L}_{\underline{w}}(x))=   \sum_{x\leq y\leq w} \mbox{rk}_v\left( I_{\underline{w},y}^p\right) h_{x,y}^p  .
\end{equation}

As
 $h_{y,y}^p=\mbox{rk}_v \left( I_{\underline{w},w}^p\right)=1$, for all $y\in W$ and any reduced expression $\underline{w}$ of $w$, by expanding and rearranging the terms of Equation \eqref{eq cinco}, we obtain 

\begin{equation}\label{eq siete}
		d(\mathbb{L}_{\underline{w}}(x))-\sum_{x< y< w} \mbox{rk}_v\left( I_{\underline{w},y}^p\right) h_{x,y}^p= \mbox{rk}_v \left( I_{\underline{w},x}^p\right)  +  h_{x,w}^p.
	\end{equation}

The term $d(\mathbb{L}_{\underline{w}})$ is obtained by a simple calculation in the Hecke algebra, so we consider it as always known.

Let us suppose for a moment that $p=0$. As $\mbox{rk}_v \left( I_{\underline{w},y}^0\right)$ is self-dual under the involution $\overline{(-)}:v\mapsto v^{-1}$ and $h^0_{x,w}=h_{x,w}\in v\mathbb{Z}[v]$, if we know the left-hand side of Equation \ref{eq siete}, we can easily calculate  both $\mbox{rk}_v \left( I_{\underline{w},y}^0\right)$ and $h^0_{x,w}$. So by induction on  $l(w)-l(x)$ one can obtain all Kazhdan-Lusztig polynomials and all the graded ranks of the intersection forms with this algorithm. 

In positive characteristic things don't work as smoothly as in characteristic zero. But the $\tilde{A}_1$ case is still nice, because in that case $h^p_{x,y}\in \mathbb{Z}[v]$\footnote{This inclusion is not even true  in finite type $A$, as proved recently in \cite{LW17}.}. This inclusion can be seen by the  fact that in that case there exist no light leaves of negative degree \cite[Corollary 3.8]{ElLi}. This inclusion implies that if one knows $h^p_{x,y}(0)$ for all pair $x,y$, one can calculate as before all the $p$-Kazhdan-Lusztig polynomials recursively. A closed formula for $h^p_{x,y}(0)$ is given in \cite[Lemma 5.1]{JW}. In order to introduce such a formula we need a bit more of notation. In type $\tilde{A}_1$ the corresponding affine Weyl group is the infinite dihedral group $W=\langle s,t\mbox{ }|\mbox{ } s^2=t^2=e \rangle$. In this group, each element different from the identity has a unique reduced expression of the form
\begin{equation}
	k_s:=sts\ldots \, \, \, (k\mbox{-terms}) \qquad \mbox{ and } \qquad k_t:= tst\ldots \, \,  (k\mbox{-terms}),
\end{equation}
for some integer $k\geq 1$. We use  the convention $0_s=0_t=e$. \\
Given two non-negative integers $a$ and $b$ we write their $p$-adic expansion as
\begin{equation}\label{p-adic}
  a=a_0+a_1p+a_2p^2+\cdots +a_{r}p^{r} \qquad \mbox{ and } \qquad  b= b_0+b_1p+b_2p^2+\cdots +b_{s}p^{s},
\end{equation}
where $0\leq a_i < p$, $0\leq  b_{i} <p$, $a_r\neq 0$ and $b_s\neq 0$. Then we say that $a$ contains $b$ to base $p$ if $s<r$ and $b_i=0$ or $b_{i}=a_{i}$, for all $0\leq i \leq s$. We then define the function $f_{p}: \mathbb{Z}_{\geq 0} \times \mathbb{Z}_{\geq 0} \rightarrow \{0,1\} $ given by
\begin{equation}\label{function p contains}
  f_p (a,b)= \left\{
               \begin{array}{ll}
                 1, & \mbox{if } a+1 \mbox{ contains } b \mbox{ to base } p;   \\
                 0, & \hbox{otherwise.}
               \end{array}
             \right.
\end{equation}

Let $y=k_s\in W$, for some $k\geq 1$. Then, we have
\begin{equation} \label{equation p-poly evaluated at zero}
h^p_{x,y}(0)=\left\{ \begin{array}{rl}
f_p(k-1,j),	& \mbox{ if } x=(k-2j)_s \mbox{ for some } 0\leq j \leq \lceil \frac{k-2}{2}\rceil ; \\
0, &   \mbox{ otherwise.}  
\end{array}  \right.
 \end{equation}

Of course, equation (\ref{equation p-poly evaluated at zero}) remains true if we replace $s$ by $t$. 

\section{Generalized blob algebras are graded cellular}\label{gba}

We will give two equivalent versions of generalized blob algebras. Both of them will be needed in the sequel. We will work throughout this paper with Bowman's homogeneous presentation \cite[Section 12]{bowman2017many} of these algebras rather than their original definition by Martin and Woodcock \cite[Section 1.2]{martin2003generalized}.

\subsection{Generalized blob algebra, algebraic definition}
From now on, fix  integers \linebreak $e,l, n>1$ and set $I_e=\ZZ / e\ZZ$.  We refer to the elements of  $ I_e^l$ as \emph{multicharges}. An \emph{adjacency-free multicharge} $\kappa =(\kappa_1,\ldots , \kappa_l)$ is a multicharge such that  $\kappa_i \not \in \{\kappa_j,\kappa_j+1,\kappa_j-1\}$ for all $i\neq j$. For such a multicharge to exist, the inequality $e\geq 2l$ has to hold.  Given $i\in I_e $ we define
\begin{equation}
	\langle i\,  |\, \kappa \rangle =|\{\, j \, |\,  1\leq j\leq l\mbox{, }\kappa_j=i\, \}|.
\end{equation} 
This number if always $0$ or $1$ if $\kappa$ is an adjacency-free multicharge. 
\begin{defi}
	Given integers $e,l, n>1$ and an adjacency-free multicharge $\kappa \in I_e^l$ the \emph{generalized blob algebra $\gblob $ of level $l$ on $n$ strings} is defined to be the unital, associative  $\mathbb{F}_p $-algebra with generators
	$$  \{ \psi_1,\ldots , \psi_{n-1} \}\cup \{y_1,\ldots , y_n\}\cup \{ e(\bi) \mbox{ }|\mbox{ } \bi \in I_e^n\}  $$

and relations

\begin{align}
\label{kl2} e(\bi )e(\boldsymbol{j})& =\delta_{\textbf{i,j}} \, e(\bi), & \\
\label{kl3}\sum_{\bi \in I_e^n} e(\bi)& =1, \\
\label{kl4}y_{r}e(\bi)& =e(\bi)y_r ,& \\
\label{kl5}\psi_r e(\bi )& =e(s_r\, \bi) \psi_r, \\
\label{kl6}y_ry_s& = y_sy_r,&
&   \\
\label{kl7}\psi_ry_s& = y_s\psi_r,&
 \mbox{ if } s\neq r,r+1 \\
\label{kl8}\psi_r\psi_s & = \psi_s\psi_r,&
 \mbox{ if } |s-r|>1 \\
\label{kl9}\psi_ry_{r+1}e(\bi )& =
(y_r\psi_r -\delta_{i_r,i_{r+1}})e(\bi )  &  \\
\label{kl10} y_{r+1}\psi_re(\bi )& =(\psi_ry_r -\delta_{i_r,i_{r+1}})e(\bi )  &  \\
\label{kl11}\psi_r^{2}e(\bi)& =\left\{ \begin{array}{l}
0   \\
e(\bi)  \\
(y_{r+1}-y_{r})e(\bi ) \\
(y_{r}-y_{r+1})e(\bi )  \\
\end{array}
\right. &  \begin{array}{l}
\mbox{if   }  i_r=i_{r+1}  \\
\mbox{if   }   i_r \neq , i_{r+1}, i_{r+1} \pm 1    \\
\mbox{if   }  i_{r+1}= i_r+1   \\
\mbox{if   }  i_{r+1}= i_r-1  \\
\end{array} 
\end{align}

\begin{align}
\label{kl12}\psi_r\psi_{r+1}\psi_re(\bi )& =\left\{ \begin{array}{l}
(\psi_{r+1}\psi_r\psi_{r +1} -1)e(\bi )  \\
(\psi_{r+1}\psi_r\psi_{r +1} +1)e(\bi )    \\
(\psi_{r+1}\psi_r\psi_{r +1} )e(\bi )    \\
\end{array}
\right. &  \begin{array}{l}
\mbox{if   }  i_{r+2}=i_r=i_{r+1}-1   \\
\mbox{if   } i_{r+2}=i_r=i_{r+1}+1      \\
\mbox{otherwise   }        \\
\end{array}  \\
 \label{Cyclotomic relation} y_1^{\langle i_1| \kappa \rangle}e(\bi )&=0  & \\
 \label{Blob relation}  e(\bi )&=0 & \mbox{ if } i_2=i_1+1.
 \end{align}

where for a sequence $\bi \in I_e^n$, $i_j\in I_e $ denotes the $j$'th coordinate of $\bi$. 

\vspace{0.3cm}
We consider $\gblob$ as a graded algebra by decreeing 

$$   \begin{array}{ccccc}
       deg\mbox{  } e(\bi )=0,  &  &  deg\mbox{  } y_r=2,  &  &
 deg\mbox{  }\psi_se(\bi )=\left\{
      \begin{array}{rl}
        -2, & \hbox{if } i_s=i_{s+1} \\
        1, & \hbox{if } i_s = i_{s+1}\pm 1  \\
        0, & \hbox{otherwise }. 
      \end{array}
    \right. \end{array}
$$
 \end{defi}
 

\begin{remark} \label{remark generalized generalized} \rm
Bowman defined in  \cite[Section 12.2]{bowman2017many}  a family of algebras $\mathcal{Q}_{l,h,n}(\kappa)$ which are a one parameter generalization of the Generalized blob  algebras. More precisely, $\gblob$ is the specialization of $\mathcal{Q}_{l,h,n}(\kappa)$ at $h=1$.  
\end{remark}

\subsection{Generalized blob algebras, diagrammatic definition}

We will explain a way to ``draw" the elements of $\gblob$. This can be achieved by using a variant of the diagrammatic calculus of Khovanov and Lauda \cite{KhLa}. 

A \emph{Khovanov-Lauda diagram on $n$-strings} (or simply a diagram when no confusion is possible) consists of $n$ points on each of two parallel edges (the top edge and bottom edge) and $n$ arcs connecting the $n$ points on one edge with the $n$ points on the other edge. Arcs can intersect, but no triple intersections are allowed. Each arc can be decorated by a finite number of dots, but dots cannot be located on an intersection of two arcs. Finally, each diagram is labelled by a sequence $\boldsymbol{i}=(i_1,\ldots ,i_n)\in I_e^n$, which is written below the bottom edge. An example of such diagrams is depicted below. 

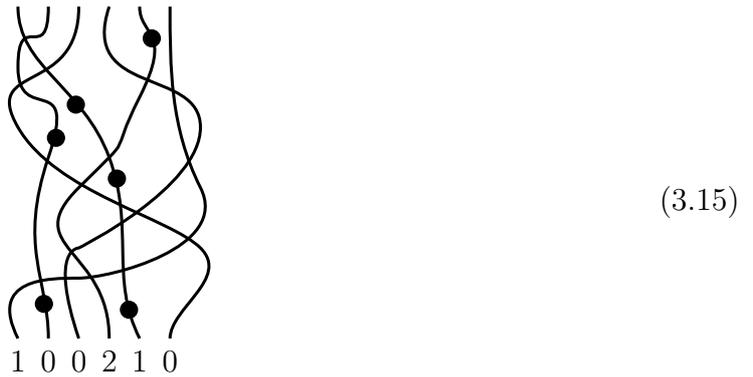
\begin{equation}\label{exa KhL}
	\exampleKhL
\end{equation}

Given a diagram $D$ we denote by $b(D)$ the bottom sequence labelling $D$ (read from left to right). This sequence determines, in the obvious fashion, a top sequence, which we denote by $t(D)$. For example, if $D$ is the diagram in (\ref{exa KhL}) then $t(D)=(1,0,0,0,2,1)$.

Let us define the diagrammatic algebra $\gblob'$. As a vector space it consists of the $\mathbb{F}_p-$linear combinations of Khovanov-Lauda diagrams on $n$-strings modulo planar isotopy and modulo the following  relations:

\begin{equation} \label{Diagrammatic Crossing dot}
	\CruceRelationA = \CruceRelationB - \delta_{i,j}  \quad \CruceRelationC\qquad \qquad 	\CruceRelationAA = \CruceRelationBB - \delta_{i,j}  \quad \CruceRelationCC    
\end{equation}

\vspace{.5cm}

\begin{equation}  \label{Diagrammatic Quadratic Relation}
	\QuadraticRelationA =\left\{\begin{array}{rl}
	0\quad , &  \mbox{if } i=j; \\
	\QuadraticRelationB \quad ,  & \mbox{if } |i-j|>1;   \\
	  &   \\
	\QuadraticRelationD \quad - \quad \QuadraticRelationC \quad,   & \mbox{if } j=i+1;\\
		  &   \\
	\QuadraticRelationC \quad - \quad \QuadraticRelationD \quad ,   & \mbox{if } j=i-1.
	\end{array} \right.
\end{equation}

\vspace{.5cm}

\begin{equation} \label{Trenza diagramtical}
	\TrenzaRelationA \quad = \quad \TrenzaRelationB \quad + \alpha \quad \TrenzaRelationC \quad ,
\end{equation}

\vspace{.3cm}
where $\alpha =-1$ if $i=k=j-1$, $\alpha = 1$ if $i=k=j+1$ and $\alpha =0$ otherwise.

\begin{align*}
	\UltimasRelaciones  & = 0,  & & \mbox{ if } i_1=\kappa_j \mbox{ for some } 1\leq j \leq l; \\
			&   &  & \\ 
	\UltimasRelacionesA & = 0, &  &\mbox{ if } i_1\neq \kappa_j \mbox{ for all } 1\leq j \leq l;    \\
		&   &  & \\ 
	\UltimasRelacionesA & =0, &  & \mbox{ if } i_2=i_1+1.
\end{align*}

\vspace{0.4cm}

The multiplication $DD'$ between two diagrams $D$ and $D'$ in $\gblob'$ is defined by vertical concatenation ($D$ on top of $D'$) when $b(D)=t(D')$. If this last equality  is not satisfied, the product is defined to be zero. We extend this product to all of $\gblob'$ by linearity. With this we end the definition of $\gblob'$. The following result can be obtained by a combination of Theorem 6.13 and Corollary 12.3 of  \cite{bowman2017many}.

\begin{teo}
	The map given by 
	 \begin{equation}
  	e( \boldsymbol{i}) \mapsto \DiagrammaticGeneratorA  \qquad y_re( \boldsymbol{i}) \mapsto \DiagrammaticGeneratorB \qquad \psi_r e(\boldsymbol{i}) \mapsto \DiagrammaticGeneratorC
  	  \end{equation}
  	  \smallskip
  	  induces an $\mathbb{F}_p$-algebra isomorphism between $\gblob$ and $\gblob'$. 
\end{teo}

In view of the above result we do not distinguish between the algebraic and the diagrammatic presentation of the generalized blob algebras.

\subsection{Graded cellular algebras} 
We briefly recall the theory of graded cellular algebras, which were defined by Hu and Mathas \cite{hu2010graded} as a graded  version of Graham and Lehrer's cellular algebras \cite{graham1996cellular}.

\begin{defi} \label{Definition cellular algebras}
   Let $A$ be an associative, $\ZZ$-graded  and finite-dimensional $\mathbb{F}_p $-algebra. Then, $A$  is called a graded cellular algebra with graded cell datum $(\Lambda , T, C, \deg)$ if the following conditions are satisfied:
  \begin{itemize}
   \item[(a)] $(\Lambda , \geq ) $ is a finite poset. For each $\lambda \in \Lambda$ there is a finite set $T(\lambda)$. Furthermore, the algebra $A$ has an $\mathbb{F}_p$-basis
                      $ \mathcal{C}=\{ c_{\mathfrak{st}}^{\lambda} \mbox{ } | \mbox{ } \lambda \in \Lambda\mbox{, } \mathfrak{s,t} \in T(\lambda)  \}$.
     \item[(b)] The $\mathbb{F}_p$-linear map $\ast $ determined by $(c_{\mathfrak{st}}^{\lambda})^{\ast}=c_{\mathfrak{ts}}^{\lambda}$ is an $\mathbb{F}_p$-algebra anti-isomorphism of $A$.
     \item[(c)] Let $\lambda \in \Lambda$, $\mathfrak{s}, \mathfrak{t} \in T(\lambda)$ and $a\in A$. Then, there are scalars $r_{\mathfrak{sv}}(a) \in \mathbb{F}_p$ (not depending on $\mathfrak{t}$) such that
         \begin{equation}\label{multiplication in a cellular algebra}
           ac_{\mathfrak{st}}^{\lambda} \equiv \sum_{\mathfrak{v} \in T(\lambda)} r_{\mathfrak{sv}}(a) c_{\mathfrak{vt}}^{\lambda} \mod A^{>\lambda},
         \end{equation}
         where $A^{>\lambda}$ is the subspace of $A$ spanned by $\{ c_{\mathfrak{ab}}^{\mu} \mbox{ }|\mbox{ } \mu > \lambda\mbox{, } \mathfrak{a}, \mathfrak{b} \in T(\mu)  \}$.
         \item[(d)] $\deg: \amalg_{\lambda \in \Lambda} T(\lambda) \rightarrow \ZZ$ is a function and the  elements of $\mathcal{C}$ are homogeneous with $\deg (c_{\mathfrak{st}}^{\lambda} ) = \deg  (\mathfrak{s}) +\deg (\mathfrak{t})$.
           \end{itemize}
The set $\mathcal{C}$ is called a \emph{graded cellular basis} of $A$.
\end{defi}

 Let $A$ be a graded cellular algebra. For $\lambda\in \Lambda$ let  $\Delta_k^\lambda$ be  the $\mathbb{F}_p$-linear span of  some formal symbols 
 $\{c_{\mathfrak{s}}^{\lambda} \mbox{ }| \mbox{ } \mathfrak{s} \in T(\lambda) \mbox{ and } \deg ({\mathfrak{s}} )=k  \}.$
 Consider the direct sum \begin{equation}
 	\Delta^{\lambda} :=\bigoplus_{k\in \ZZ} \Delta_k^\lambda.
 \end{equation}
 One defines an action of  $A$ on $\Delta^{\lambda}$  by
\begin{equation}\label{action cell modules}
  ac_{\mathfrak{s}}^{\lambda } =\sum_{\mathfrak{v} \in T(\lambda)} r_{\mathfrak{sv}}(a) c_{\mathfrak{v}}^{\lambda },
\end{equation}
where $r_{\mathfrak{sv}}(a)$ are the scalars that appear in Definition \ref{Definition cellular algebras}(c). We call the graded left $A$-module $\Delta^{\lambda}$ 
  a \emph{graded cell module}, where of course $\Delta_k^\lambda$ lives in degree $k$.

  Given a $\ZZ$-graded $A$-module $M=\oplus_{k\in \ZZ} M_k$ we define its \emph{graded dimension} as
 \begin{equation}\label{graded dimension}
   \gdim M = \sum_{k\in \ZZ} \left( \dim M_{k} \right) v^{k} \in \ZZ[v,v^{-1}].
 \end{equation}
It is clear from the definitions that
\begin{equation} \label{dim of a graded cell module}
	\gdim \Delta^{\lambda} = \sum_{\mathfrak{s}\in T(\lambda)}v^{\deg (\mathfrak{s})}.
\end{equation}

Each  graded cell module $\Delta^{\lambda}$ comes equipped with a symmetric and associative bilinear form $\langle \mbox{ , } \rangle_{\lambda}$ determined by 
\begin{equation}
	c_{\ag \sg}^\lambda c_{\tg \bg}^\lambda \equiv \langle c_{\sg}^\lambda , c_{\tg}^\lambda \rangle_\lambda c_{\ag \bg}^\lambda \mod A^{>\lambda}. 
\end{equation}

  The radical of this form,
  $\mbox{Rad}(\Delta^{\lambda}) := \{ x\in \Delta^{\lambda} \mbox{ }| \mbox{ } \langle x,y \rangle_{\lambda }=0\mbox{, for all } y\in \Delta^{\lambda} \}$, is a graded $A$-submodule of $\Delta^{\lambda}$.  Given $\lambda \in \Lambda$, we define $L^{\lambda}:=\Delta^{\lambda} /\mbox{Rad}(\Delta^{\lambda})  $ and $$\Lambda_0=\{\lambda \in \Lambda \mbox{ }| \mbox{ }L^{\lambda}\neq 0 \}.$$\\

\begin{teo} \label{Theorem clasification simples in a graded cellular algebra}\cite[Theorem 3.4]{graham1996cellular} \cite[Theorem 2.10]{hu2010graded}
  The set $\{ L^{\lambda}(j)\mbox{ }| \mbox{ } \lambda \in \Lambda_0 \mbox{ and } j\in \ZZ \}$ is a complete  set of pairwise non-isomorphic graded simple $A$-modules.
\end{teo}

 Let $\mathcal{G}(A)$ be the graded Grothendieck group of $A$. This is the  $\lau$-module generated by the symbols $[M]$, where $M$ runs over the set of all finite-dimensional graded right $A$-modules, subject to the following relations:
\begin{itemize}
  \item $[M(j) ] =v^{j}[M]$, for all $j\in \ZZ$.
  \item $[M]=[N]+[P]$, if  $0\rightarrow N \rightarrow M \rightarrow P \rightarrow 0 $ is a short exact sequence of graded right $A$-modules.
\end{itemize}

\begin{cor}\label{theorem classification in terms og Grothendieck}
Let $A$ be a graded cellular algebra. Then,  $\{ [L^{\lambda}] \mbox{ } | \mbox{ } \lambda \in \Lambda_0  \}$ is an $\lau $-basis of $\mathcal{G}(A)$.
\end{cor}  

For any  $\lambda \in \Lambda$ and $\mu \in \Lambda_{0}$, we define $d_{\lambda, \mu}(v) \in \lau$, the \emph{graded decomposition numbers of} $A$,  by the formula
\begin{equation}\label{expansion cell in terms of simples}
  [\Delta^{\lambda}] = \sum_{\mu \in \Lambda_0} d_{\lambda, \mu}(v) [L^{\mu}].
\end{equation}
By general theory of graded cellular algebras \cite[Lemma 2.13]{hu2010graded} we know that
\begin{equation}\label{properties decomposition numbers}
  d_{\mu , \mu}(v) =1 \mbox{ and that } d_{\lambda,\mu}(v) =0 \mbox{ if }  \lambda \ngeq \mu ,
  \end{equation}
for all $\lambda \in \Lambda$ and $\mu \in \Lambda_0$. 

\begin{defi}
	Let $A$ be a graded algebra and  $M$ be a graded $A$-module. We say that $M$ is \emph{positively graded} if $\gdim M \in \mathbb{Z}_{\geq 0}[v]$. In addition, we say  that $M$ is \emph{pure of degree} $d$ if $\gdim M = nv^d$, for some $n\in \mathbb{N}$. A \emph{positively graded cellular algebra} is one in which the image of the function $\mathrm{deg}$ is in $\mathbb{Z}_{\geq 0}$.
\end{defi}
The following is a recollection of some useful facts about these kind of algebras.  
\begin{lem}\cite[Lemma 2.20]{hu2010graded}  \label{lemma properties positively graded algebras}
	Let $A$ be a positively graded cellular algebra with graded cell datum as in Definition \ref{Definition cellular algebras}. Then, 
	\begin{enumerate}
	\item  $\Delta^\lambda$ is positively graded, for all $\lambda \in \Lambda$.
	\item If $\mu \in \Lambda_0$ then there exists $\mathfrak{t}\in T(\mu)$ such that $\deg (\mathfrak{t})=0$. 
	\item $L^{\mu} $ is pure of degree zero, for all $\mu \in \Lambda_0$.
	\item $d_{\lambda , \mu }(v)\in \mathbb{Z}_{\geq 0}[v]$, for all $\lambda \in \Lambda$ and $\mu \in \Lambda_0$.
	\end{enumerate}
\end{lem}

\subsection{A graded cellular basis for $B_n(\kappa)$.} \label{sec gcb for blob}
In this section we recall the construction by Bowman \cite{bowman2017many} of a graded cellular basis of $\gblob$. His construction is for a larger family of algebras (see Remark \ref{remark generalized generalized}). We omit the subtleties of this more general construction that do not arise in our context. In order to define the graded cellular basis for $\gblob$ we need to introduce  some combinatorial objects. 

\subsubsection{One column multipartitions}

For $n\geq 0,$ a \emph{partition} $\lambda $ of $n $ is a weakly decreasing sequence of non-negative integers $(\lambda_1,\lambda_2, \ldots)$ such that $|\lambda|:=\lambda_1 + \lambda_2+ \cdots =n$. An $l$-\emph{multipartition} $\lamn=(\lambda^{(1)}, \lambda^{(2)}, \ldots, \lambda^{(l)}) $ of $n$, is an $l$-tuple of partitions such that $|\lambda^{(1)}|+|\lambda^{(2)}|+\cdots |\lambda^{(l)}|=n$. 
The set of all $l$-multipartitions of $n$ is denoted by $\mpar$.

In this paper, we are exclusively interested in \emph{one-column $l$-multipartitions of $n$}, that is, $l$-multipartitions $\lamn =(\lambda^{(1)}, \lambda^{(2)}, \ldots, \lambda^{(l)})$ such that $0\leq \lambda^{(m)}_1\leq 1$, for all $1\leq m \leq l$. The set of all one-column $l$-multipartitions of $n$ is denoted  by $\ocmpar $. Another way of looking at this is to say that an element of $\ocmpar$ is a sequence of $l$ terms $(1^{a_1}, \cdots, 1^{a_l})$ where $a_1+\cdots +a_l=n$ and $a_i\in \mathbb{Z}_{\geq 0}$. For example, $\lamn=(1^3,1^4, 1^0,1^3, 1^1)\in P_1^5(11).$

For an element $\lamn=(1^{a_1}, \cdots, 1^{a_l}) \in \ocmpar$  the \emph{Young diagram} $[\lamn]$ of $\lamn$  is the set of ``boxes" or ``nodes"
\begin{equation}	\bigcup_{1\leq m \leq l}\{ (r,m) \mbox{ }|\,   1\leq r\leq a_m  \}.
\end{equation}
We say that a node $(r,m)\in [\lamn ]$ is in the $r^{\mathrm{th}}$ row  of the $m^{\mathrm{th}}$ component of $\lamn$.

	Let $(r,m)$ and $(r',m')$ be two boxes. We say that $(r,m)$ \emph{dominates} $(r',m')$ and write $(r,m)\rhd (r',m')$ if $(r,m)$ is smaller than $(r',m')$ in the lexicographical order, i.e.,  if $r<r'$ or  $r=r'$ and $m<m'$. In the Young diagram, a box $B$ is bigger than a box $B'$ if $B'$ is in the same row than $B$ but strictly to the right or in a strictly lower row. Let $\lamn, \mun \in \ocmpar $. We say that $\lamn $ \emph{dominates} $\mun$ and write $\lamn \rhd \mun $ if for each box $(r,m)\in [\mun]  $ the number of boxes in $[\lamn] $ that dominate $(r,m)$ is greater  or equal than the number of boxes in $[\mun]$ that dominate $(r,m)$.

Let $\lamn \in \ocmpar$. A \emph{tableau of shape} $\lamn$ is a bijection $\tg :[\lamn]\rightarrow \{1,\ldots, n\}$. We think of $\tg$ as filling in the boxes of $[\lamn]$ with the numbers $1,2, \ldots , n$. For $1\leq i\leq n$,  $\tg^{-1}(i)$ is the box of the Young diagram $[\lamn]$ filled in by the number $i$.

  A tableau is called \emph{standard} if its entries increase from top to bottom on each component. We denote by $\Std(\lamn )$ the set of all standard tableaux of shape $\lamn$.  For $\tg \in \Std (\lamn)$ we write $\mbox{Shape}(\tg )=\lamn $. Given $\tg \in \Std(\lamn )$ and $1\leq k \leq n$ we define $\tg \downarrow_k$ as the standard tableau obtained from $\tg$ by erasing the boxes with entries strictly greater than $k$.

Let $\lamn \in \ocmpar $. We say that a box $(r,m)\in [\lamn] $ is \emph{removable} from  $\lamn $ if
 $[\lamn] \backslash \{(r,m)\}  $ is the Young diagram of some element in $P_1^l(n-1)$. Similarly, we say that a box  $(r,m)\not \in [\lamn] $ is \emph{addable} to $\lamn$ if $[\lamn] \cup \{(r,m)\}  $ is the Young diagram of some element in  $P_1^l(n+1)$. We denote by $\RC (\lamn)$ (resp. $\AC (\lamn)$) the set of all removable (resp. addable) boxes of $\lamn$. 
 
 Let us define $\mbox{res}(r,m)$  the \emph{residue of the box} $(r,m)$ by 
\begin{equation}\label{res}
	\mbox{res}(r,m) = \kappa_m +1 -r \in I_e.
\end{equation}

 \subsubsection{Definition of $\mathrm{deg}$ in $B_n(\kappa)$} This definition is not needed for the rest of the paper. We give it for completeness.

 Given $1\leq k \leq n$ and a standard tableau $\tg$ we let $\RC_\tg(k)$ (resp. $\AC_\tg(k)$)  denote the set of all boxes $(r,m)$ satisfying the following three conditions:
\begin{enumerate}
    \item $ (r,m)\in \RC (\mbox{Shape}(\tg\downarrow_k  ))$ (resp. $\AC (\mbox{Shape}(\tg\downarrow_k  )) )$
	\item $\mbox{res}(r,m) = \mbox{res}(\tg^{-1}(k))$
	\item $\tg^{-1}(k) \rhd (r,m)$
\end{enumerate}

	Let $\lamn\in \ocmpar$ and $\tg \in \Std (\lamn)$. The \emph{degree} of $\tg$ is defined as 
	\begin{equation}
		\deg (\tg) = \sum_{k=1}^n \left( |\AC_\tg(k)|-|\RC_\tg(k)| \right).
 	\end{equation}

\subsubsection{Residue sequences and $\tg \rightsquigarrow  d_{\tg}$} For each $\lamn \in \ocmpar$ we define a tableau $\tg^{\lamn} \in \Std(\lamn )$ called the \emph{dominant tableau} ({{for  example, see the first two tableaux in }} Picture \ref{Exa tableaux}). It  is determined by the following  rule: 
\begin{equation}
	(\tg^{\lamn})^{-1}(k)\rhd (\tg^{\lamn})^{-1}(j) \mbox{ if and only if } k<j.  
\end{equation}

Given $\lamn \in \ocmpar$ and $\tg \in \Std(\lamn )$ we denote by $d_{\tg}\in \SG_n$ the permutation determined by $d_{\tg} \tg^{\lamn} = \tg$. 

	For each $\tg \in \Std(\lambda )$, we define its \emph{residue sequence} (see definition \ref{res}) as follows
	\begin{equation}
		\boldsymbol{i}^{\tg} = (\mbox{res} (\tg^{-1}(1)),\mbox{res} (\tg^{-1}(2)), \ldots ,\mbox{res} (\tg^{-1}(n)) )\in I^n_e.
	\end{equation}
	For notational convenience, we write 
	$\boldsymbol{i}^{\lamn}:= \boldsymbol{i}^{\tg^{\lamn}}$.

\begin{defi} \label{defi elements psi}
Let $\lamn \in \ocmpar $ and $\sg , \tg \in \Std(\lamn )$. We fix reduced expressions $d_\sg = s_{i_1}\cdots s_{i_a}$ and $d_{\tg}=s_{j_1}\cdots s_{j_{b}}$. Then, we define 
\begin{equation}
	\psi_{\sg\tg}^{\lamn}:= \psi_{i_1}\cdots \psi_{i_a} e(\boldsymbol{i}^{\lamn } ) \psi_{j_b} \cdots \psi_{j_1} \in \gblob .
\end{equation}	
\end{defi}

\subsubsection{The graded cellular basis}
\begin{pro} \label{Theorem Graded cellular basis blob}\cite[Theorem 7.1]{bowman2017many}
	The set $\{  \psi_{\sg\tg}^{\lamn} \mbox{ }| \mbox{ } \lamn\in \ocmpar\mbox{, } \sg , \tg \in \Std(\lamn) \}$ is a graded cellular basis of the generalized blob algebra $\gblob$ with respect to the dominance order on $\ocmpar $, the degree function defined above and the involution $\ast$ determined by flipping upside down the diagrams.  
\end{pro}

Once we have specified a graded cellular structure on $\gblob$ we have automatically defined graded cell modules and graded simple modules. Given $\lamn \in \ocmpar $, we will denote them by $\cellb $ and $\simb$, respectively. Since $\psi_{\tg^{\lamn}\tg^{\lamn}}^{\lamn}=e(\boldsymbol{i}^{\lamn} )$ one can see that the bilinear form defined by the cellular structure on $\cellb$ is distinct from zero for each $\lamn \in \ocmpar$. In other words, $\simb \neq 0$, for each  $\lamn \in \ocmpar$. Then, we have well-defined graded decomposition numbers $d_{\lamn, \mun}^p\in \ZM [v,v^{-1}]$, for each pair  $\lamn , \mun \in \ocmpar $. 

\begin{exa} \label{example weno weno weno} \rm
Let $n=23$, $l=4$, $e=8$ and $\kappa=(0,2,4,6)\in I_{8}^4$.  Let $\lamn=(1^1,1^{13},1^1,1^8)$, $\mun = (1^5,1^5,1^6,1^7) \in P_1^4(23)$. In (\ref{Exa tableaux}) we have drawn the dominant tableaux $\tg^{\lamn}$ and $\tg^{\mun}$ corresponding to $\lamn$ and $\mun$. We have also drawn another tableau $\tg\in \Std (\mun)$.    

\begin{equation} \label{Exa tableaux} 
\tg^{\lamn }= 
\begin{ytableau}
1 
\end{ytableau}
\mbox{ }
\begin{ytableau}
2 \\
5 \\
7\\
9\\
11\\
13\\
15\\
17 \\
19\\
20\\
21\\
22\\
23
\end{ytableau}
\mbox{ }
\begin{ytableau}
3
\end{ytableau}
\mbox{ }
\begin{ytableau}
4 \\
6 \\
8\\
10\\
12\\
14\\
16\\
18
\end{ytableau}
\qquad 
\tg^{\mun} = 
\begin{ytableau}
1 \\
5 \\
9 \\
13 \\
17
\end{ytableau}
\mbox{ }
\begin{ytableau}
2 \\
6 \\
10\\
14\\
18
\end{ytableau}
\mbox{ }
\begin{ytableau}
3 \\
7 \\
11 \\
15 \\
19 \\
21
\end{ytableau}
\mbox{ }
\begin{ytableau}
4 \\
8 \\
12\\
16\\
20\\
22\\
23
\end{ytableau}
\qquad
\tg = 
\begin{ytableau}
1 \\
9 \\
11 \\
13 \\
15  
\end{ytableau}
\mbox{ }
\begin{ytableau}
2 \\
5 \\
7\\
18 \\
23 
\end{ytableau}
\mbox{ }
\begin{ytableau}
3 \\
10 \\
12 \\
14 \\
16 \\
22
\end{ytableau}
\mbox{ }
\begin{ytableau}
4 \\
6 \\
8\\
17\\
19\\
20\\
21
\end{ytableau}
\end{equation}	

The permutation $d_\tg \in \mathfrak{S}_{23}$ is given by  
$$  \left( \begin{array}{ccccccccccccccccccccccc}
  1&2&3&4&5&6&7&8&9&10&11&12&13&14&15&16&17&18&19&20&21&22&23  \\
 1&2&3&4&9&5&10&6&11&7&12&8&13&18&14&17&15&23&16&19&22&20&21
\end{array}
    \right) .  $$

Of course, there are many reduced expressions for $d_\tg$. However, any fixed but arbitrary choice of such a expression can be used to construct a graded cellular basis for  $\gblob$.

It is not hard to see that for any one-column $l$-multipartition $\lamn$  the sets $\mathcal{A}_{\tg^{\lamn }} (k) $ and 
$\mathcal{R}_{\tg^{\lamn } }(k)$ are empty, for all $1\leq k\leq n$. Therefore, we have  $\deg (\tg^{\lamn})= \deg (\tg^{\mun})=0$. On the other hand, for the tableau  $\tg $ all the sets $\mathcal{R}_\tg(k)$ are empty, but there are six non-empty sets $\mathcal{A}_\tg(k)$, namely, 
$$
\begin{array}{rrr}
	\mathcal{A}_\tg(9)=\{ (4,2) \}, &	\mathcal{A}_\tg(10)=\{ (4,4) \}, & 	\mathcal{A}_\tg(17)=\{ (6,1) \}, \\ 	\mathcal{A}_\tg(18)=\{ (6,3) \}, & \mathcal{A}_\tg(22)=\{ (8,4) \},  & 	\mathcal{A}_\tg(23)=\{ (7,3) \}.
\end{array}
$$

Therefore,  $\deg (\tg )=6$. \\
The residue sequences for the tableaux $\tg^{\lamn}$ and $\tg^{\mun}$ are given by
$$   \bi^{\lamn} =(0,2,4,6,1,5,0,4,7,3,6,2,5,1,4,0,3,7,2,1,0,7,6), $$
$$   \bi^{\mun}  =(0,2,4,6,7,1,3,5,6,0,2,4,5,7,1,3,4,6,0,2,7,1,0). $$

We  remark that the residue sequence associated to $\tg$ coincides with the one associated to $\tg^{\lamn}$, that is, $\bi^{\lamn} = \bi^{\tg}$. A main part of the forthcoming  Section 4 is devoted to understand the set of all standard tableaux with the same residue sequence as some dominant tableau. For the moment, we want to roughly explain why are there $2^6$ standard tableaux with residue sequence $\bi^{\lamn}$.

Let $\sg $ be a standard tableau with $\bi^\sg = \bi^{\lamn}$. The first eight entries in $\sg$ are determined by  $\bi^{\lamn}$, so that  $\sg \downarrow_8$ corresponds to the tableau on top of the Figure (\ref{Example Tableaux same residues}). In  $\mbox{Shape} (\sg \downarrow_8 )$ we have two addable boxes with residue $7$ (which is the ninth entry in  $\bi^{\lamn}$), namely, the boxes $(2,1)$ and $(4,2)$. It follows that we can locate the number $9$ in any of these two boxes keeping the same residue sequence as $\bi^{\lamn}$. Then we have two  options for $\sg \downarrow_{9}$  which are drawn in the middle of (\ref{Example Tableaux same residues}). A similar situation occurs for the number $10$, thus we have four options for $\sg\downarrow_{10} $, which are depicted in the bottom of (\ref{Example Tableaux same residues}). For any of the four options for  $\sg \downarrow_{10}$ we can see that the number $11$  is forced to be located below the number $9$. This is because for any of these tableaux the box below the number $9$  is the unique  addable box with residue $6$ (which is the eleventh entry in $\bi^{\lamn}$). Then, we still have only four options for $\sg\downarrow_{11}$. By continuing in this way we can see that in each step we will have one or two options to locate the next numbers. The numbers in which we will have two possibilities are $17$, $18$, $22$ and $23$.  Therefore, we have  $2^6$ standard tableaux $\sg$ with $\bi^\sg=\bi^{\lamn}$.

\begin{equation} \label{Example Tableaux same residues}
\begin{array}{ccc}
  & \scalebox{.8}{ \begin{ytableau}
1   
\end{ytableau}
\,
\begin{ytableau}
2 \\
5 \\
7 
\end{ytableau}
\,
\begin{ytableau}
3
\end{ytableau}
\,
\begin{ytableau}
4 \\
6 \\
8
\end{ytableau} } &  \\
 \scalebox{.8}{\begin{ytableau}
1\\
9
\end{ytableau}
\mbox{ }
\begin{ytableau}
2 \\
5 \\
7 
\end{ytableau}
\mbox{ }
\begin{ytableau}
3
\end{ytableau}
\mbox{ }
\begin{ytableau}
4 \\
6 \\
8
\end{ytableau} } & &  \scalebox{.8}{ \begin{ytableau}
1
\end{ytableau}
\mbox{ }
\begin{ytableau}
2 \\
5 \\
7 \\
9
\end{ytableau}
\mbox{ }
\begin{ytableau}
3
\end{ytableau}
\mbox{ }
\begin{ytableau}
4 \\
6 \\
8
\end{ytableau} } 
\end{array}
\end{equation}

$$
\begin{array}{c|c|c|c}
\scalebox{.8}{\begin{ytableau}
1\\
9
\end{ytableau}
\mbox{ }
\begin{ytableau}
2 \\
5 \\
7 
\end{ytableau}
\mbox{ }
\begin{ytableau}
3 \\
10
\end{ytableau}
\mbox{ }
\begin{ytableau}
4 \\
6 \\
8
\end{ytableau} } &  \scalebox{.8}{\begin{ytableau}
1\\
9
\end{ytableau}
\mbox{ }
\begin{ytableau}
2 \\
5 \\
7 
\end{ytableau}
\mbox{ }
\begin{ytableau}
3
\end{ytableau}
\mbox{ }
\begin{ytableau}
4 \\
6 \\
8 \\
10
\end{ytableau}}  & \scalebox{.8}{\begin{ytableau}
1
\end{ytableau}
\mbox{ }
\begin{ytableau}
2 \\
5 \\
7  \\
9
\end{ytableau}
\mbox{ }
\begin{ytableau}
3 \\
10
\end{ytableau}
\mbox{ }
\begin{ytableau}
4 \\
6 \\
8 
\end{ytableau}}  & \scalebox{.8}{\begin{ytableau}
1
\end{ytableau}
\mbox{ }
\begin{ytableau}
2 \\
5 \\
7  \\
9
\end{ytableau}
\mbox{ }
\begin{ytableau}
3 
\end{ytableau}
\mbox{ }
\begin{ytableau}
4 \\
6 \\
8 \\
10
\end{ytableau}}
\end{array}
$$
 
\end{exa}


\subsection{Idempotent truncations of $\gblob$.}  \label{Section Idempotent truncations}

In this section we consider certain subalgebras of $\gblob$  which are obtained by idempotent truncation.  As before,  $e, l, n>0$ and $\kappa$ is an increasing adjacency-free multicharge.

 We define the \emph{idempotent truncation of} $\gblob$ \emph{at} $\lamn$ as 
$\gblobl := e(\boldsymbol{i}^{\lamn}) \gblob e(\boldsymbol{i}^{\lamn } )$, for each $\lamn\in \ocmpar$.  More generally, given $\lamn$ and $\mun$ in $\ocmpar$ we define $B_{l,n}^{p,\lamn,\mun} := e(\boldsymbol{i}^{\mun})  B_{l,n}^{p} e(\boldsymbol{i}^{\lamn})$. At the diagrammatic level  $B_{l,n}^{p,\lamn,\mun}$ corresponds to the subspace of $B_{l,n}^{p}$ generated by diagrams with bottom $n$-tuple given by $\boldsymbol{i}^{\lamn}$
 and top $n$-tuple given by $\boldsymbol{i}^{\mun}$. It follows from the relations (\ref{kl2}), (\ref{kl5}) and from Definition \ref{defi elements psi} that 
\begin{equation} \label{orthogonality blob basis}
	e(\boldsymbol{i}^{\mun}) \psi_{\sg\tg}^{\nun }e(\boldsymbol{i}^{\lamn})= 
	\left\{ 
	\begin{array}{rl}
		\psi_{\sg \tg}^{\nun}, & \mbox{ if } \boldsymbol{i}^{\sg} =\boldsymbol{i}^{\mun} \mbox{ and } \boldsymbol{i}^{\tg} =\boldsymbol{i}^{\lamn}; \\
		0, & \mbox{ otherwise.}
	\end{array}
	\right.
\end{equation}

Therefore,  a basis for $\gblobl $ (resp. $B_{l,n}^{p,\lamn,\mun}$) can be obtained by considering a subset of the graded cellular basis of $\gblob $. Let us be more precise. For any $\mun \in \ocmpar$ we define $\Std_{\lamn}(\mun ) = \{ \tg \in \Std (\mun )\mbox{ }|\mbox{ } \boldsymbol{i}^{\tg} =
	\boldsymbol{i}^{\lamn } \}$.  We also define $ \plam =\{ \mun\in \ocmpar \mbox{ } |\mbox{ } \Std_{\lamn}(\mun ) \neq \emptyset \}$.

\begin{lem} \label{lemma basis idempotent truncation}
For each $\lamn \in \ocmpar $ the set $\{ \psi_{\sg\tg}^{\mun} \mbox{ } |\mbox{ } \mun \in \plam\mbox{,  } \sg , \tg \in \Std_{\lamn}(\mun ) \}$ is a graded cellular basis of $\gblobl  $ with respect to the order, degree and involution inherited of  $\gblob $.  More generally, a basis for $B_{l,n}^{p,\lamn,\mun}$ is given by 
\begin{equation} \label{basis hom spaces truncated}
\{ \psi_{\sg\tg}^{\nun} \,|\, \nun \in P(\lamn) \cap P(\mun); \, \sg \in \Std_{\mun} (\nun) , \, \tg \in \Std_{\lamn} (\nun)    \}.
\end{equation}
	\end{lem}

\begin{dem}
	It follows directly from (\ref{orthogonality blob basis}). 
\end{dem}

Given $\mun \in \plam $, we denote by $\cellbl  $ and $\simbl$ the corresponding graded cell and simple modules of $\gblobl$, respectively. It is not hard to see  that $\cellbl \simeq e(\boldsymbol{i}^{\lamn}) \Delta^p(\mun)  $ and $\simbl \simeq  e(\boldsymbol{i}^{\lamn}) L^p(\mun)$. As $\simbl$ is a quotient of $\cellbl$,
\begin{equation}\label{cell zero implies simple zero}
	\cellbl =0 \Rightarrow \simbl =0. 
\end{equation}
The following formula is Equation \eqref{dim of a graded cell module} applied to this set up.
\begin{equation}  \label{dim cell module idempotent}
	\gdim \cellbl =\sum_{\tg\in \Std_{\lamn} ( \mun )} v^{\deg(\tg)}. 
\end{equation} 
Let $\mun , \nun \in P(\lamn)$, with $L_{\lamn}(\nun) \neq 0$. We  denote by $d_{\mun , \nun}^{p,\lamn}$ the graded decomposition number associated to $\mun$ and $\nun$ in $\gblobl$.  Graded decomposition numbers of $\gblob$ and $\gblobl$ are compatible, this is, $d_{\mun , \nun}^p =d_{\mun , \nun}^{p,\lamn}  $ if $d_{\mun , \nun}^{p,\lamn}$  is defined. This  can be proved  using essentially the same arguments as in the ungraded case; see \cite[Appendix A1(4)(vii)]{donkin1998q}. \\
On the other hand, suppose we want to calculate $d_{\mun , \nun}^p$, for some pair $(\mun , \nun )\in \ocmpar^2$. In this case we work with $B_{l,n}^{p,\nun}(\kappa)$. We have $\Std_{\nun}(\nun )=\{ \tg^{\nun}  \}$. Therefore, 
\begin{equation}
	\Delta_{\nun}(\nun) \cong L_{\nun} (\nun ) \cong \mbox{Span}_{\mathbb{F}_p}\{ \psi_{\tg^{\nun}} \}. 
\end{equation}   
In particular, $L_{\nun} (\nun )\neq 0$. If $\mun \in P(\nun)$  then $d_{\mun, \nun}^{p,\nun}$ is defined, so that $d_{\mun , \nun}^{p,\nun} =d_{\mun , \nun}^p$. 
If $ \mun \not \in P(\nun)$ (and therefore $e(\bi^{\nun }) \Delta^p(\mun)=0$) then we can  again use \cite[Appendix A1(4)(vii)]{donkin1998q} to obtain  $d_{\mun , \nun}^p=0$. Summing up, we conclude that  each  non-zero graded decomposition number of $\gblob $ arises as a graded decomposition number of some idempotent truncation. For this reason, in the forthcoming sections we study idempotent truncations of $\gblob$ rather than  $\gblob$ itself.

\medskip
We conclude this section by making explicit the description of the functor  $F_u^v$ mentioned in Section \ref{subsection an equivalent conjecture}. Let $\lamn_u \in P_1^l(n_u)$ and $\lamn_{v}\in P_1^l(n_{v})$ be elements in the same orbit with $n_u<n_{v}$ and $w_{\lamn_{u}} <w_{\lamn_{v}}$ (according to the notation considered in Section \ref{subsection an equivalent conjecture} this last condition is the same as $u<v$). Let $\mun $ and $\nun$ be elements in $P(\lamn_u)$. Let $\overline{\mun}$ and $\overline{\nun}$ be the images of $\mun $ and $\nun$ under the action of $F_u^v$, respectively. We recall that $\overline{\mun}$ (resp.  $\overline{\nun}$) is obtained from $\mun$ (resp. $\nun$) by adding $(n_{v}-n_u)/l$ boxes in each component. Let $\boldsymbol{\rho} \in P(\mun) \cap P( \nun) $ and $\mathfrak{s} \in \Std (\boldsymbol{\rho}) $.  We can define an element $\overline{\mathfrak{s}} \in \Std (\overline{\boldsymbol{\rho} } )$ as follows. First, we add $(n_{v}-n_u)$ to each one of the entries of $\mathfrak{s}$. Then, we add $(n_{v}-n_u)/l$ empty boxes at the top of each component. We finally fill the empty boxes with the numbers $1, 2, \ldots , n_{v}-n_u$ from left to right and top to bottom, i.e., in the same way as the numbers are filled into a dominant tableau. Let $\psi_{\mathfrak{st}}^{\boldsymbol{\rho}}$ be an element of the basis of $B_{l,n_i}^{p,\mun , \nun}$ given in (\ref{basis hom spaces truncated}). We define $F_u^v(\psi_{\mathfrak{st}}^{\boldsymbol{\rho}}) = \psi_{\overline{\mathfrak{s}}\overline{\mathfrak{t}} }^{\overline{\boldsymbol{\rho }}}$ and extend it by linearity. 
	
	 The fact that $F_u^v: \mbox{Blob}^l(\leq \lamn_u) \rightarrow \mbox{Blob}^l(\leq \lamn_{v})$ preserves composition of morphisms follows once we notice that $\psi_{\overline{\mathfrak{s}}\overline{\mathfrak{t}} }^{\overline{\boldsymbol{\rho }}}$ is obtained from $\psi_{\mathfrak{st}}^{\boldsymbol{\rho}} $ by adding $(n_{v}-n_u)$ straight lines on its leftmost region. Finally, the fully faithfulness of $F_{u}^v$ follows from the fact that it sends a basis of the Hom space in the source category to a base in the corresponding Hom space of the target category.

\section{Alcove geometry}\label{alcove}
The main  goal  of this section is to prove that  $\gdim \cellbl$  coincides with $\deg (\mathbb{L}_{\underline{w_{\lamn}}} (w_{\mun}) )$, where $w_{\lamn}$ and $w_{\mun}$ are some elements in $W_l$ that we will introduce and $\underline{w}_{\lamn}$ is a particular reduced expression of $w_{\lamn}$. To do this it is convenient to interpret standard one-column tableaux {
as paths in an alcove geometry. 

\subsection{Affine reflections and alcoves} Let $E_l$ denote the quotient space $\RM^l/\langle \epsilon_1+\epsilon_{2}+\cdots +\epsilon_l \rangle $, where $\epsilon_i$ denotes the $i$-th coordinate vector of $\RM^l$.  Let $\odot $ denote the origin in $E_l$. We identify elements  $\lamn \in \ocmpar$ with points (resp. classes) in $\mathbb{R}^l$ (resp. $E_l$) via the map $\lamn =(1^{a_1}, \ldots 1^{a_l}) \mapsto \sum_{i=1}^l a_i\epsilon_i$.  For each $1\leq i<j\leq l$ and each $m\in \ZM$ define the affine hyperplane
\begin{equation}
	\hg_{i,j}^m= \{ x\in E_l \mbox{ } |\mbox{ } x_i-x_j=\kappa_i-\kappa_j +me \}
\end{equation}
and $s_{i,j}^m $ the corresponding affine reflection. In formulas 
\begin{equation}
	s_{i,j}^m (x)= x-( (x_i-x_j)-(\kappa_i-\kappa_j +me))(\epsilon_i-\epsilon_j). 
\end{equation}

Let $W_{l}$ be the affine Weyl group of type $\tilde{A}_{l-1}$. We identify $W_l$ with the group generated by the affine reflections via the assignment $s_{0}\mapsto s_{1,l}^{1}$ and  $s_{i}\mapsto s_{i,i+1}^0$, for $1\leq i<l$. 


Let $\HG' $ be the set of affine hyperplanes and $\HG $ be their union. The connected components of $E_l\backslash \HG $ are called \emph{alcoves}. Let $\mathcal{A}$ denote the set of alcoves.  A point $x\in E_l$ is called \emph{regular} if it belongs to  some alcove and $\lamn \in P_1^l(n)$ is called \emph{regular} if its image in $E_l$ is regular. By the restriction imposed on $\kappa$ it is clear that $\odot $ is regular.

 The alcove containing $\odot$ is called the fundamental alcove and is denoted by $A_0$. It is well-known  that $A_0$ has  $l$ walls which are supported by the hyperplanes $\{\hg_{i,i+1}^0\mbox{ } |\mbox{ }  1\leq i<l\}$ and $\hg_{1,l}^1$. The action of $W_{l}$ on $E_l$ extends to an action on $\AC$. Given $A\in \AC $ there exists a unique $w\in W_{l}$ such that $w\cdot A_0=A$. If this is the case, we denote $A:=A_w$.

\subsection{Paths and sequences of hyperplanes associated to standard tableaux}
 \begin{defi}
	Given  $\lamn \in \ocmpar$ and $\tg\in \Std (\lamn )$ we define the path associated to $\tg$ as the piecewise linear path $p_{\tg}:[0,n]\rightarrow E_l$ with vertices $p_\tg (0)=\odot$ and $p_\tg (k)= \mbox{Shape}(\tg\downarrow_k)\in E_l$, for $1\leq k \leq n$. We denote $p_{\lamn}:= p_{\tg^{\lamn}}$.
\end{defi}

We remark that if $\lamn \in \ocmpar$ and  $\tg\in \Std (\lamn) $ then $p_{\tg}(n)=\lamn$. In particular, $p_{\lamn} (n)=\lamn$. The following lemma is \cite[Lemma 4.7]{bowman2017family}.    
 
\begin{lem}  \label{lemma residues and paths}
Let $\lamn \in \ocmpar  $ and $\tg \in \Std(\lamn )$. Then, the vertex  $p_{\tg}(k)$ of the path $p_{\tg}$ lies in a hyperplane $\hg_{i,j}^m$ for some $m\in \mathbb{Z}$ if and only if the addable boxes in the $i$-th and $j$-th components of  $\mbox{Shape}(\tg\downarrow_{k})$ have the same residue. 
 \end{lem}

From now and on we fix $\lamn \in \ocmpar$ regular. In accordance with Lemma \ref{lemma residues and paths},  the regularity of $\lamn$ is equivalent to the fact that the  boxes addable to $\lamn$ have different residues.  We denote by  $w_{\lamn} \in W_l$  the unique element such that $\lamn \in A_{w_{\lamn }}$. In order to establish a connection between this alcove geometry and the light leaves basis we need a particular reduced expression for $w_{\lamn}$. To do this we begin by associating to $\lamn$ a sequence of hyperplanes $\hg_{\lamn}$ given, in order, by the new hyperplanes touched by the path $p_{\lamn}$. Let us be more precise. 
For $1 \leq k \leq n$ we define 
	\begin{equation}
		\hg(k,\lamn)= \{  \hg \in \HG' \mbox{ } | \mbox{ } p_{\lamn} (k) \in \hg   \} \backslash  \{  \hg \in \HG' \mbox{ } | \mbox{ } p_{\lamn} (k-1) \in \hg    \}. 
	\end{equation} 

We want to prove that the set 	$\hg(k,\lamn)$  consists of at most one element, for all $1\leq k \leq n$. We will first need a preliminary result.

\begin{lem} \label{Lemma At most two addable}
	Let $\lamn \in \ocmpar $ regular. Then, $\mbox{Shape}(\tg^{\lamn}\downarrow_k)$ has at most two addable boxes with the same residue, for each $0\leq k\leq n$.
\end{lem}

\begin{dem}
Suppose there exists $k$ such that $\mbox{Shape}(\tg^{\lamn}\downarrow_k)$ has three addable boxes $A_1$, $A_2$ and $A_3$ with the same residue. Let $r_i$ be the row to which $A_i$ belongs and $r$ the row to which the box occupied by $k$ in $\tg^{\lamn}$ belongs. Since $\kappa$ is assumed to be adjacency-free we know that $r_1$, $r_2$ and $r_3$  are different and not adjacent, that is, $r_i\neq r_j\pm 1$. In particular, two of the elements of $\{r_1,r_2 ,r_3\}$ are less than $r-1$. By definition of  $\tg^{\lamn}$,  the numbers greater than $k$ are placed in rows lower than or equal to the $r^{\mathrm{th}}$ row. It follows that at least two of the boxes $A_1$, $A_2$ and $A_3$ are addable to $\mbox{Shape}(\tg^{\lamn}\downarrow_{k'})$, for all $k\leq k'\leq n$. In particular, we have that $\lamn = \mbox{Shape}(\tg^{\lamn}\downarrow_{n})$ has two addable boxes with the same residue, which contradicts the hypothesis of $\lamn$ being regular. \end{dem}

\begin{cor}  \label{coro unique hyperplane}
	Let $\lamn \in \ocmpar  $ be regular. Then,  the set $\hg(k,\lamn)$  consists of at most one element, for all $1\leq k \leq n$. 
\end{cor}	
	
\begin{dem} 
Let $1\leq k\leq n$. Let $A$ be the box occupied by $k$ in $\tg^{\lamn}$ and let $i$ be the component to which $A$ belongs.  The paths $p_{\lamn }(k-1)$ and $p_{\lamn }(k)$ only differ in the  $i$-th component, so the unique addable box to 
$\mbox{Shape}(\tg^{\lamn}\downarrow_k )$ which is not addable to $\mbox{Shape}(\tg^{\lamn}\downarrow_{k-1})$ is the one located just below $A$.

If  $\hg(k,\lamn)$ has two or more elements, then  Lemma \ref{lemma residues and paths} implies $ \mbox{Shape}(\tg^{\lamn}\downarrow_k)$ has three or more addable boxes with the same residue, contradicting  Lemma  \ref{Lemma At most two addable}.
\end{dem}

Let $k_1<k_2< \ldots < k_r$ be the integers such that $\hg({k_i},\lamn )\neq \emptyset$. In accordance with Corollary \ref{coro unique hyperplane}, we define $\hg_i$ to be the unique hyperplane belonging to $\hg({k_i},\lamn )$. Then, we define the \emph{sequence of hyperplanes associated to} $\lamn$ by $\hg_{\lamn} =(\hg_1,\ldots , \hg_r)$.  By convention, if $\lamn \in A_0$ then $\hg_{\lamn}$  is the empty sequence.

 It is possible that some vertex $p_{\lamn} (k)$ of $p_{\lamn}$ belongs to the intersection of two or more hyperplanes. For instance, in Example \ref{example weno weno weno}, $p_{\lamn} (8)$ belongs to $\hg_{1,2}^0\cap \hg_{3,4}^0$. However, if some vertex of $p_{\lamn} $ is on  the intersection of some hyperplanes, say $\hg_{i_1,j_1}^{m_1}$ and $\hg_{i_2,j_2}^{m_2}$, then, because of Lemma \ref{lemma residues and paths} and  Lemma \ref{Lemma At most two addable}, we have that $\{ i_{1},j_1  \}\cap \{i_2,j_2 \}=\emptyset $. In this case, following \cite[Definition 2.7]{bowman2017family}, we say the hyperplanes $\hg_{i_1,j_1}^{m_1}$ and $\hg_{i_2,j_2}^{m_2}$ are \emph{orthogonal}.  



 \begin{exa} \rm
 	We keep the same parameters as in Example \ref{example weno weno weno}, so that  $\kappa= (0,2,4,6)$, $e=8$ and $\lamn = (1^1,1^{13},1^1,1^8)$. We want to determine the integers $k_i$ such that $\hg({k_i},\lamn )\neq \emptyset$. According to Lemma \ref{lemma residues and paths}, these integers will be the ones satisfying that $\mbox{Shape}(\tg^{\lamn}\downarrow_{k_i} ) $ has a pair of addable boxes of the same residue but one of the boxes in this pair is not addable to 
 	$\mbox{Shape}(\tg^{\lamn}\downarrow_{k_i-1} )$. So, we have
$$ 	k_1= 7,\,  k_2=8, \,  k_3=15 , \,  k_4=16, \,  k_5=21, \,  k_6 =22.   $$
 The sequence $\hg_{\lamn}$  is given by $ \hg_{\lamn} =(\hg_{1,2}^0, \hg_{3,4}^0 ,\hg_{2,3}^1,\hg_{1,4}^0,\hg_{1,2}^{-1},\hg_{2,4}^1)   $.	
 \end{exa}

\subsection{Some properties of $\hg_{\lamn}$}
Each hyperplane $\hg_{i,j}^m\in \HG'$ splits $E_l$ into two half-spaces
\begin{equation} \label{half-spaces}
\small
	E_{i,j}^m(+):=\{ x\in E_l \mbox{  }|\mbox{ } x_{i}-x_j >\kappa_{i}-\kappa_j+me \}  \mbox{ and }  E_{i,j}^m(-):=\{ x\in E_l \mbox{  }|\mbox{ } x_{i}-x_j <\kappa_{i}-\kappa_j+me \}. 
\end{equation} 
We say that a hyperplane \emph{separates} two alcoves if they belong to distinct half-spaces. It is a known fact that $l(w)$ coincides with the number of hyperplanes separating $A_w$ from $A_0$ (see, for example, \cite{humphreys1992reflection}).

\begin{lem}  \label{lemma hyperplanes}
	Let $\lamn=(1^{a_1},\ldots, 1^{a_l}) \in \ocmpar $ regular and let $\hg_{\lamn}=(\hg_1,\ldots \hg_r)$ be the corresponding sequence of hyperplanes. Then,
	\begin{enumerate}
		\item The hyperplanes $\hg_i$ are all distinct.
		\item Each $\hg_{i}$ separates $A_{w_{\lamn}}$ and $A_0$.
		\item Any hyperplane separating $A_{w_{\lamn}}$ and $A_0$ appears in $\hg_{\lamn}$.
		\item $l(w_{\lamn})=r$.
	\end{enumerate}  
\end{lem}

\begin{dem}
Let us prove (1). Let us suppose that the hyperplane $\hg_{i,j}^m $ appears in $\hg_{\lamn}$. 
We must have  $a_i\neq a_j$ because if $a_i= a_j$, by the way that the numbers are filled in $\tg^{\lamn}$, for any $k$ the addable boxes in the $i$-th and $j$-th components of $\mbox{Shape}(\tg^{\lamn}\downarrow_{k})$  would never have the same residue   since $\kappa$ is adjacency-free. This, by Lemma \ref{lemma residues and paths}  contradicts the hypothesis  that the hyperplane $\hg_{i,j}^m $ appears in  $\hg_{\lamn}$.

 Without loss of generality, assume $a_i >a_j$. For $1\leq k\leq n$, we set 
 \begin{equation}
 	(1^{\mu^k_1},\ldots ,1^{ \mu^k_l}):=\mbox{Shape}(\tg^{\lamn}\downarrow_{k})=p_{\lamn}(k).
 \end{equation}
  Let $a$ be the minimal integer satisfying $\hg(a,\lamn )=\{  \hg_{i,j}^m   \}$ 
  and let $b$ be the integer which occurs  in the box located just below the box occupied by $a$ in $\tg^{\lamn}$ (in the $i$-th component). This number must exist, otherwise $\lamn $ would lie on $ \hg_{i,j}^m  $ and it would not be regular. Let us remark that since $a_i >a_j$ we have $a_j=\mu^k_{j}$ for all $a\leq k\leq n$.
  
  \begin{itemize}
  	\item If $a<k<b$ then $p_{\lamn}(k)$ and $p_{\lamn}(k-1)$ lie on $\hg_{i,j}^m $, and therefore $\hg_{i,j}^m \notin  \hg (k,\lamn )$.
  	\item  If $b\leq k\leq n$ then $ \mu^k_i-\mu^k_j >\mu^a_i-\mu^a_j =\kappa_i-\kappa_j+me$, thus $ p_{\lamn}(k)\notin \hg_{i,j}^m $.
  \end{itemize}
  We conclude that $\hg_{i,j}^m$ cannot occur twice in $\hg_{\lamn}$, thus proving (1).

\medskip
Let us prove (2). Let $\hg_{i,j}^m$ be a hyperplane occurring in $\hg_{\lamn }$. Since $\kappa$ is increasing and adjacency-free we know that
$-e< \kappa_i-\kappa_j <0$. Therefore, $\odot\in E_{i,j}^m (+)$ if $m<1$ and  $\odot \in E_{i,j}^m (-)$ if $m\geq 1$. Let $a$ be the (unique) integer given by $\hg(a,\lamn)=\{ \hg_{i,j}^m \}$. As before, we set $(1^{\mu^a_1},\ldots , 1^{\mu^a_l})=\mbox{Shape}(\tg^{\lamn}\downarrow_{a})$. We have
\begin{equation}\label{eq mu separates}
	\mu^a_i-\mu^a_j = \kappa_i-\kappa_j +me. 
\end{equation}
Suppose that $m\geq 1$ (the case $m<1$ is similar).  By (\ref{eq mu separates}) we obtain $\mu^a_i>\mu^a_j$. As in point (1) $\mu^a_j=a_j$ and $\mu^a_i <a_i $. It follows that $a_i-a_j >\mu_i^a -\mu_j^a = \kappa_i-\kappa_j +me$. 
Thus, $\lamn \in E_{i,j}^{m}(+)$ and $\hg_{i,j}^m$ separates $A_{w_{\lamn}}$ and $A_0$. This proves (2).

\medskip

If a hyperplane separates $A_0$ and $A_{w_{\lamn}}$, it must be hit by $p_{\lamn}$ (a continuous path) thus it appears in  $\hg_{\lamn}$ and (3) follows. (4) is a direct consequence of the previous statements.
\end{dem}

\vspace{0.4cm}

\subsection{The alcove sequence $\ag_{\lamn}$}

Let $\hg_{\lamn}=(\hg_1, \cdots, \hg_r)$ be the sequence of hyperplanes associated to $\lamn$. Let us define a  sequence of alcoves $\ag_{\lamn} = (\ag_0, \ag_1, \ldots,\ag_r  ) $ as follows.  Firstly define $\ag_0=A_0$. The following alcoves in the sequence are inductively defined by
 $\ag_{i} = \rho_{i}\ag_{i-1}$, where $\rho_i$ denotes the orthogonal reflection through the hyperplane $\hg_i$.  
 
 Let $k_1<k_2<\ldots < k_r$ be the integers such that $\hg (k_i,\lamn) =\{\hg_i\}$. Set $p_0=p_{\lamn}$. We are going to define paths $p_1, p_2, \ldots , p_r$ as follows. Suppose the path $p_i$ has been already defined. Then, we define $p_{i+1}$ as the path obtained from $p_{i}$ by applying the reflection $\rho_1\rho_2 \cdots \rho_i \cdots \rho_2 \rho_1$ to all the vertices of $p_{i}$ after the vertex $p_{i}(k_{i+1})$.  For instance, keeping the same notation and parameters as in Example \ref{example weno weno weno} the path $p_r$ is the one which corresponds to the tableu $\tg$. 
  Since $\lamn$ is regular we know that $k_r<n$, so that in each step of the above construction the relevant reflection was applied to last vertex. In other words, 
 \begin{align*}
 	p_r(n) & = (\rho_1\rho_2 \cdots \rho_r \cdots \rho_2 \rho_1)\cdots (\rho_1\rho_2\rho_1) (\rho_1)(p_0(n))\\
 	   & = (\rho_1\rho_2 \cdots \rho_r \cdots \rho_2 \rho_1)\cdots (\rho_1\rho_2\rho_1) (\rho_1)(\lamn )\\
 	   & = \rho_1\rho_2 \cdots \rho_r (\lamn). 
 \end{align*} 
Since $\hg_{\lamn}=(\hg_1, \cdots, \hg_r)$ corresponds to the sequence of all hyperplanes touched by $p_{\lamn}=p_0$ in order, the path $p_r$ is completely contained in the closure of the fundamental alcove. In particular, $p_r(n)$ belongs to the closure of the fundamental alcove.  As $\lamn$ is regular and regularity is not affected by the action of $W_l$ we conclude that $p_r(n)=\rho_1\rho_2 \cdots \rho_r (\lamn)$ belongs to the fundamental alcove. Therefore, $w_{\lamn} = \rho_r\cdots \rho_2\rho_1$.
 

Given two alcoves $A$ and $B$ we say that they are \emph{adjacent} if $B$ is obtained from $A$ by applying an affine reflection with respect to one of the walls of $A$. In other words, two alcoves are adjacent if they are different and share a common wall. We remark that two adjacent alcoves share $l-2$ common walls. A sequence of alcoves $(\bg_0,\bg_{1}, \ldots, \bg_r)$ is called an \emph{alcove path} if $\bg_0=A_0$ and $\bg_{i-1}$ and $\bg_i$ are adjacent, for $1\leq i\leq r$.  An alcove path   $(\bg_0,\bg_{1}, \ldots, \bg_r)$ is called \emph{reduced} if $r$ is minimal among the set of all alcove paths ending in $\bg_r$.

\begin{lem}
	Let $\lamn\in \ocmpar$ regular. The alcove sequence $\ag_{\lamn} =(\ag_0,\ag_{1}, \ldots ,\ag_r   )$ of $\lamn$ is a reduced alcove path.
\end{lem}  

\begin{dem}
	Let $\hg_{\lamn} =(\hg_1, \ldots , \hg_r)$ be the hyperplane sequence associated to $\lamn$. Recall that $\rho_i$ is the reflection with respect the hyperplane $\hg_i$. Since by definition $\ag_0=A_0$, in order to prove that $\ag_{\lamn}$ is an alcove path it is enough to show that $\hg_{i}$ is a wall of $\ag_{i-1}$.  We proceed by induction. The hyperplane $\hg_1$ is the first hyperplane intersected by $p_{\lamn}$.  As $p_{\lamn}$ is continuous and it begins at $\odot \in A_0$ we conclude that $\hg_1$ must be a wall of the fundamental alcove $A_0=\ag_0$. This yields the base of our induction. 
	
We now suppose that $i>1$ and that $\hg_{i-1}$ is a wall of $\ag_{i-2}$. As $\ag_{i-1} =\rho_{i-1}\ag_{i-2}$ we know that $\hg_{i-1}$ is a wall of $\ag_{i-1}$ as well. Let $k_{i-1}$ and $k_{i}$ be the integers such that $\hg (k_{i-1},\lamn)=\{ \hg_{i-1}  \}$ and $\hg(k_{i},\lamn )=\{\hg_i  \}$, respectively. Then, $p_{\lamn} (k_{i-1})\in \hg_{i-1}$ and $p_{\lamn} (k_i)\in \hg_{i}$. The proof splits into two cases.

\medskip


\medskip
\textbf{Case A: $p_{\lamn} (k_i) \in \hg_{i-1}$.}  In this case we have $p_{\lamn} (k_i) \in \hg_{i-1} \cap \hg_i$ and therefore  $\hg_{i-1} $ and $ \hg_i $ are orthogonal. Let $\hg$ be a wall of $\ag_{i-2}$ distinct of $\hg_{i-1}$ and orthogonal to $\hg_{i-1}$. Since $\ag_{i-1} =\rho_{i-1}\ag_{i-2}$ we obtain that $\hg =\rho_{i-1}\hg $ is a wall of  $\ag_{i-1}$.  So, it suffices to prove that $\hg_i$ is a wall of $\ag_{i-2}$. This follows by noticing that from the vertex $k_{i-1}$ to the vertex  $k_i$ the path $p_{\lamn}$ is always on $\hg_{i-1}$ and it does not intersect a new hyperplane till its intersection with $\hg_i$. Then, as $\hg_{i-1} $ is a wall of $\ag_{i-2}$, $\hg_{i}$ must also be a wall of $\ag_{i-2}$. 

\medskip
\textbf{Case B: $p_{\lamn} (k_i) \notin \hg_{i-1}$.} In this case the path $p_{\lamn}$ leaves the hyperplane $\hg_{i-1}$ before the vertex $k_i$. Actually, by the way the numbers are filled in  $\tg^{\lamn}$, $p_{\lamn}$ must leave any hyperplane before the vertex $k_i$. In other words, $p_{\lamn} (k_i-1)$ belongs to an alcove. By construction this alcove must be $\ag_{i-1}$. Therefore, the next hyperplane that $p_{\lamn}$ intersects must be a wall of $\ag_{i-1}$. But this hyperplane is by definition $\hg_i$. 

\medskip
By  \cite[Lemma 5.3]{lenart2007affine} we know that the number of alcoves involved in a reduced alcove path ending  in $A_{w_{\lamn}}$ is $l(w_{\lamn})+1 $. Then, $\ag_{\lamn}$ is reduced by Lemma \ref{lemma hyperplanes}.
\end{dem}

\medskip
We are now in position to obtain a reduced expression for $w_{\lamn}$.  Recall that $\hg_{\lamn}=(\hg_1,\ldots ,\hg_r)$ and $\ag_{\lamn} =(\ag_0,\ag_1, \ldots , \ag_r)$ are the hyperplane and alcove sequences associated to $\lamn$. Recall that $\rho_j$ is the reflection through $\hg_j$. By definition of $\ag_{\lamn}$ we know that $\ag_j=\rho_j\rho_{j-1} \cdots \rho_1\ag_0$, for each $1\leq j\leq r$. We remark that if $A$ and $B$ are adjacent alcoves and $w\in W_l $  then $w A$ and $w B$ are also adjacent. Then, applying $\rho_1 \cdots \rho_{j-1}$ to the pair of adjacent alcoves $\ag_{j-1} $ and $\ag_j$, we obtain the pair of adjacent alcoves $\rho_1 \cdots \rho_{j-1}\ag_{j-1} =\ag_0$ and $\rho_1 \cdots \rho_{j-1}\ag_{j}= \rho_1 \cdots \rho_{j-1} \rho_j  \rho_{j-1}\cdots \rho_1\ag_0$. It follows that  $\rho_1 \cdots \rho_{j-1} \rho_j  \rho_{j-1}\cdots \rho_1$ is a reflection with respect  some wall of the fundamental alcove. So, for each $1\leq j\leq r$, there is an $i_j \in \{0,1,\ldots, l-1\}$ satisfying that $\rho_1 \cdots \rho_{j-1} \rho_j  \rho_{j-1}\cdots \rho_1$ is the simple reflection $s_{i_j}$. Clearly, we have 
$$s_{i_1} \cdots s_{j_r}= (\rho_1) (\rho_1\rho_2\rho_1)\cdots (\rho_1\rho_{2} \cdots \rho_{r-1}\rho_r\rho_{r-1}\cdots\rho_{2}\rho_1)= \rho_r\rho_{r-1}\cdots\rho_{2}\rho_1 =w_{\lamn}.$$

 Therefore, $\underline{w_{\lamn}}= s_{i_1} \cdots s_{i_r}$ is a reduced expression of $w_{\lamn} $ which we call the \emph{principal  reduced expression of} $w_{\lamn}$.  
 
 \begin{exa}\rm
 	With the same notation and parameters as in Example \ref{example weno weno weno} we have 
 	$$\underline{w_{\lamn}}= s_1s_3s_0s_2s_3s_2.$$
 \end{exa}
 
 \begin{teo}  \label{teo dim cell module hecke}
 	Let $\lamn ,\mun\in \ocmpar $. Suppose that $\lamn$ is regular.
 	\begin{enumerate}
 		\item  If $\mun$ does not belong to the orbit of $\lamn$ then $\Std_{\lamn}(\mun)=\emptyset$.
 		\item  If $\mun $ belongs to the orbit of $\lamn$ then $	\gdim \cellbl $ coincides with the coefficient of $H_{w_{\mun}}$ in the expansion of $\underline{H}_{\underline{w_{\lamn}}}$ in terms of the standard basis of the Hecke algebra of $W_{l}$. Consequently, 
 	\begin{equation}  \label{Equation blob vs light leaves}
 		\gdim \cellbl = \deg \mathbb{L}_{\underline{w_{\lamn}}} (w_{\mun} ).  
 	\end{equation}

 	\end{enumerate}
  \end{teo}

 \begin{dem}
 Let $k_1<k_2<\ldots <k_r$ be the integers such that $\hg(k_i,\lamn) \neq \emptyset$. We construct a perfect binary tree with nodes decorated  by paths starting at $\odot$. We construct it by induction on the level of depth. In depth zero the unique node is decorated by $p_{\lamn}$. Suppose now a node of depth $s $ has been decorated by a path $p$.  One of the two child nodes is decorated by $p$ itself. The other one is decorated by the path $p'$ which is obtained from $p$ by applying a reflection through the hyperplane $\hg$ to all the vertices of $p$ after the vertex $p(k_s)$. The hyperplane $\hg$ corresponds to the  $s$-th hyperplane touched by $p$. It is shown in \cite[Proposition 4.11]{bowman2017family} that the paths decorating  the leaves of the tree are exactly those associated to standard tableaux with residue sequence $\bi^{\lamn}$. The endpoint of all of these paths belongs to $W_l\cdot \lamn $, thus proving (1).\\ 
 By (\ref{dim cell module idempotent}) we know that
 \begin{equation} \label{cell blob dimension sum}
  \gdim \cellbl =\sum_{\tg\in \Std_{\lamn} ( \mun )} v^{\deg(\tg)}. 	
 \end{equation}

That the right-hand side coincides with the coefficient of $H_{w_{\mun}}$ in the expansion of $\underline{H}_{\underline{w_{\lamn}}}$ in terms of the standard basis of the Hecke algebra of $W_{l}$ follows by combining  Proposition 2.14 and Proposition 4.11 in \cite{bowman2017family}. Equation (\ref{Equation blob vs light leaves}) is now a consequence of (\ref{eq uno}).   
 \end{dem}


\begin{cor} \label{coro bruhat order y plambda uno}
	Let $\lamn, \mun \in \ocmpar$ with $\lamn$ regular. Then, $\Std_{\lamn}(\mun)\neq \emptyset $ if and only if $\mun$ is in the orbit of $\lamn$ and $w_{\mun} \leq w_{\lamn}$. 
\end{cor}

\begin{dem}
We first notice that 
\begin{equation}\label{eq light leaves are not empty}
	\mathbb{L}_{\underline{w_{\lamn}}} (w_{\mun} )\neq \emptyset \iff w_{\mun} \leq w_{\lamn}.
\end{equation}
Suppose that $\Std_{\lamn}(\mun)\neq \emptyset $. By Theorem \ref{teo dim cell module hecke}(1) we know that $\mun$ belongs to the orbit of $\lamn$. Furthermore, by combining \eqref{Equation blob vs light leaves} and  \eqref{cell blob dimension sum} we obtain
\begin{equation}
	\deg \mathbb{L}_{\underline{w_{\lamn}}} (w_{\mun} )= \gdim \cellbl = \sum_{\tg\in \Std_{\lamn} ( \mun )} v^{\deg(\tg)}   \neq 0 .
\end{equation}
Therefore, $\mathbb{L}_{\underline{w_{\lamn}}} (w_{\mun} )\neq \emptyset$. We conclude via \eqref{eq light leaves are not empty}  that $w_{\mun} \leq w_{\lamn}$.\\ 
Conversely, suppose that $\mun$ is in the orbit of $\lamn$ and $w_{\mun} \leq w_{\lamn}$.  A combination of \eqref{Equation blob vs light leaves}, \eqref{cell blob dimension sum} and \eqref{eq light leaves are not empty} yields
\begin{equation}
\sum_{\tg\in \Std_{\lamn} ( \mun )} v^{\deg(\tg)} = \gdim \cellbl =	\deg \mathbb{L}_{\underline{w_{\lamn}}} (w_{\mun} )   \neq 0.
\end{equation}
Therefore, $ \Std_{\lamn} ( \mun ) \neq \emptyset$. 
\end{dem}

\medskip
The following theorem is the main result of this paper. It is a slight generalization of Theorem \ref{maint} in the introduction.

\begin{teo} \label{coro graded vector space isomorphism}

Let $n\leq n'$ be positive integers. Let $\lamn\in \ocmpar $ and $\mun \in P_1^l(n') $.  Suppose that $\lamn$ and $\mun$ are regular and belong to the same orbit. Let $\lamn'$ be the unique element in $P_1^l(n')$   such that $\lamn' =\lamn $ in $E_l$. Then,

	\begin{equation}
		\mathbb{F}_p\otimes_R \mathrm{Hom}_{\HC}(BS(\underline{w_{\lamn}}), BS(\underline{w_{\mun }})) \cong B^{p,\lamn',\mun}_{l,n'}(\kappa), 
	\end{equation}
	as graded vector spaces. In particular, if $n=n'$  then we obtain that the Hom-spaces involved in the \textbf{Categorical blob vs Soergel conjecture} are isomorphic as graded vector spaces. 
\end{teo}

\begin{dem}
We start by recalling that $B^{p,\lamn',\mun}_{l,n'}(\kappa) = e(\boldsymbol{i}^{\mun} ) B_{l,n'}^p(\kappa) e(\boldsymbol{i}^{\lamn'} ) $.  By using the basis of $B^{p,\lamn',\mun}_{l,n'}(\kappa)$ given in Lemma \ref{lemma basis idempotent truncation} we have
\begin{equation}
	\gdim B^{p,\lamn',\mun}_{l,n'}(\kappa) = \sum_{\nun \in P(\lamn') \cap P(\mun)} \gdim \Delta_{\lamn'}^p(\nun) \gdim \Delta_{\mun}^p(\nun).
\end{equation}
On the other hand, we notice that $\lamn'$ is obtained by adding  $(n'-n)/l$ to each component of $\lamn$.  So that, $p_{\lamn'}$ is a concatenation of a path contained in the fundamental alcove  with $p_{\lamn}$.  It follows that the hyperplane sequence of $\lamn$ and $\lamn'$ coincide. Therefore, $\underline{w_{\lamn'}}=\underline{w_{\lamn}}$ (as reduced expressions). Then, by using (\ref{equation double leaves basis}), we obtain
\begin{equation}
	\gdim \mathbb{F}_p\otimes_R \mathrm{Hom}_{\HC}(BS(\underline{w_{\lamn}}), BS(\underline{w_{\mun }})) = \sum_{x\leq w_{\lamn},w_{\mun}} \deg \mathbb{L}_{\underline{w_{\lamn'}} }(x) \deg \mathbb{L}_{\underline{w_{\mun}} }(x).
\end{equation}
The result is now a consequence of Theorem \ref{teo dim cell module hecke}.	
\end{dem}

\begin{lem} \label{lema residuo implica mas chico}
	Let $\lamn\in \ocmpar$ regular and $\mun\in \plam $. Then, $\lamn \unlhd \mun $. 
\end{lem}

\begin{dem}
Suppose by contradiction that there exists $\tg \in \Std_{\lamn}(\mun )$ with  $\lamn \centernot{\unlhd} \mun$. Let $A$ and $B$ be the boxes occupied by $n$ in   in $\tg^{\lamn}$ and  $\tg$, respectively. Note that $A$ is the least dominant box in $[\lamn]$ and that $\res (A) = \res (B)$ since $\boldsymbol{i}^\tg= \boldsymbol{i}^{\lamn} $.   Let $\bar{\lamn}$ be the partition obtained from $\lamn$ by removing the box $A$. Similarly, let $\bar{\tg}$ be the tableau obtained from $\tg $  by removing the box $B$ and $\bar{\mun} = \mbox{Shape}(\bar{\tg})$. Let $C$ and  $D$ be the least dominant boxes in $[\bar{\lamn}]$ and  $[\bar{\mun }]$, respectively. If  $\bar{\lamn} \unlhd \bar{\mun}$ then $C\unlhd D$. Therefore, the addable boxes to $[\bar{\mun}]$ are less dominant than $A$ or are located one row below than $A$. In the later case these boxes cannot have the same residue as $A$ since our multicharge $\kappa$ is adjacency-free. So, relocating the boxes $A$ and $B$ we would obtain that $\lamn \unlhd \mun $, contradicting our assumption. Then, we must have $\bar{\lamn}\centernot{\unlhd } \bar{\mun}$. Note that $\boldsymbol{i}^{\bar{\tg}}= \boldsymbol{i}^{\bar{\lamn}}  $. So, we are in the same situation as in the beginning of the proof but with one box less. The result follows now by induction.
\end{dem}

\begin{cor} 
	Let $\lamn \in \ocmpar$ regular. Then, 
	\begin{equation} \label{description plambda}
		\plam = \{ \mun\in \ocmpar \mbox{ }| \mbox{ } \mun \in W_l \cdot \lamn \mbox{ and } w_{\mun}\leq w_{\lamn} \} 
	\end{equation} 
	and the  order inherited by $\plam$ from the dominance order on $\ocmpar $ is a refinement of reversed Bruhat order. Concretely, if  $\mun , \nun \in \plam $ and $w_{\nun} \leq w_{\mun} $ then $\mun \unlhd \nun$.
\end{cor}

\begin{dem}
	Equation (\ref{description plambda}) is just a restatement of Corollary \ref{coro bruhat order y plambda uno}. Let $\mun , \nun \in \plam $. Suppose that $w_{\nun} \leq w_{\mun} $. Then, $\nun \in P(\mun)$ and Lemma \ref{lema residuo implica mas chico} implies $\mun \unlhd \nun$. 
\end{dem}

\subsection{A process to compute graded decomposition numbers.}  \label{section algorithm to compute gdn blob}

In this section we explain a process which calculates graded decomposition numbers for $\gblob$. At this point the reader should compare this process with the one described in Section \ref{section an important formula} to compute $p$-Kazhdan-Lusztig polynomials. 

\medskip
Let $\lamn\in \ocmpar$ regular and $\mun\in \plam$.  By combining (\ref{properties decomposition numbers}) and Lemma \ref{lema residuo implica mas chico} we obtain
\begin{equation} \label{equation algorithm one}
\gdim \cellbl = \sum_{\substack{\nun \in \plam \\ \lamn \unlhd \nun \unlhd \mun}}	 d_{\mun, \nun}^p \gdim L_{\lamn}^p(\nun). 
\end{equation}
Let $\nun \in \plam$.  By (\ref{description plambda}) we know that $w_{\nun}\leq w_{\lamn}$. On the other hand, if $\nun  \unlhd \mun $ but $w_{\mun } \not \leq w_{\nun }$ we would have $ \mun \notin P(\nun )$. Then, $e(\nun) \Delta^p(\mun) =0$ and therefore $d_{\mun , \nun}^p=0$. So that, the terms on the right-hand side of (\ref{equation algorithm one}) associated to $\nun $ such that $w_{\mun } \not \leq w_{\nun }$ vanish and we can rewrite such an equation as 
\begin{equation} \label{equation algorithm two}
\gdim \cellbl = \sum_{\substack{\nun \in W\cdot \lamn \\ w_{\mun} \leq w_{\nun} \leq w_{\lamn}}}
d_{\mun, \nun}^p \gdim L_{\lamn}^p(\nun). 
\end{equation}

As $d_{\mun,\mun}^p= \gdim L_{\lamn}^p(\lamn) =1$, by expanding and rearranging the terms in (\ref{equation algorithm two}), we obtain

\begin{equation} \label{equation algorithm three}
\gdim \cellbl - \sum_{\substack{\nun \in W\cdot \lamn \\ w_{\mun} < w_{\nun} < w_{\lamn}}}d_{\mun, \nun}^p \gdim L_{\lamn}^p(\nun) = \gdim L_{\lamn}^p(\mun)+d_{\mun , \lamn}^p.
\end{equation}

Let us suppose for a moment $p=0$. In this case it follows by \cite[Theorem 3.15]{bowman2017family} and \cite[Corollary 12.3]{bowman2017many} that $d_{\mun, \lamn}^p\in v\mathbb{Z}[v]$, if $\mun \neq \lamn$. On the other hand, by general theory of graded cellular algebras \cite[Proposition 2.18]{hu2010graded}, we know that $\gdim L_{\lamn}^p(\nun)$ is invariant under the involution $v\mapsto v^{-1}$. We conclude, via  the algorithm described in Section \ref{section an important formula} and Theorem \ref{teo dim cell module hecke}, that
\begin{equation} \label{equality in char zero}
d_{\mun ,\lamn}^p=h_{w_{\mun},w_{\lamn}}^p \qquad \mbox{ and } \qquad  \gdim \simbl = \rank I_{\underline{w_{\lamn}}}^p(w_{\mun}) .
\end{equation}
 
 The equality on the left-hand side of (\ref{equality in char zero}) was proven in \cite[Theorem 4.16]{bowman2017family}. As far as we know, the equality on the right-hand side of (\ref{equality in char zero}) had not been observed before. 
 
 \medskip
 Let us return to the case $p>0$. In this setting, the inclusion $d_{\mun, \lamn}^p\in v\mathbb{Z}[v]$ is no longer true. So, that in order to calculate graded decomposition numbers we must compute graded dimensions of simple modules. We stress that we are in the same situation as in the $p$-Kazhdan-Lusztig theory. Namely, in order to calculate $p$-Kazhdan-Lusztig polynomials we have to compute graded ranks of intersection forms.
 
 \medskip
  In light of the above, we introduce the following conjectures. 

\begin{cona}\label{cona}
	Let $\lamn\in\ocmpar$ regular. Suppose that $\mun \in \plam$. Then, $$d_{\mun ,\lamn}^p=h_{w_{\mun},w_{\lamn}}^p.$$ 
\end{cona}

\begin{conb}\label{conb}
	Let $\lamn \in\ocmpar$ regular. Suppose that $\mun \in \plam$. Then, $$\gdim \simbl = \rank I_{\underline{w_{\lamn}}}^p(w_{\mun}) .$$ 
\end{conb}

By the previous discussion and Theorem \ref{teo dim cell module hecke} it is clear that	\textbf{Blob vs Soergel Conjecture} and \textbf{Blob vs Light leaves Conjecture} are equivalent.

\section{Proof of blob vs Soergel conjecture for $\tilde{A}_1$.}\label{proof}

In this section we prove that graded decomposition numbers of $\gblobdos$ coincide with $p$-Kazhdan-Lusztig polynomials of $\tilde{A}_1$, that is, we prove conjecture Blob vs Soergel (and therefore Conjecture Blob vs Light Leaves) for $\tilde{A}_1$. So, for the rest of this section  $l=2,$  $e,n\in\mathbb{N}$   and $\kappa =(\kappa_1, \kappa_2) $ is an increasing adjacency-free multicharge.

\subsection{Paths in the Pascal triangle and hooks} In type $\tilde{A}_1$ the corresponding affine Weyl group is the infinite dihedral group $W=\langle s,t\mbox{ }|\mbox{ } s^2=t^2=e \rangle$. In this group, each element different from the identity has a unique reduced expression of the form
\begin{equation}
	k_s:=sts\ldots \, \, (k\, \mbox{ terms)} \qquad \mbox{ and } \qquad k_t:= tst\ldots \, \, (k\, \mbox{  terms)},
\end{equation}
for some integer $k\geq 1$. We use  the convention $0_s=0_t=e$. 

On the other hand, the ambient space for our alcove geometry is $\mathbb{R}^2/\langle \epsilon_1 +\epsilon_2 \rangle 	\simeq \mathbb{R}$.  Therefore, a path $p_{\tg}$ associated to a standard tableau $\tg$ is drawn as a one-dimensional continuous path starting at $\odot$ and then going to the left or to the right according to the component occupied  by the  numbers in $\tg$. To draw such paths in $\mathbb{R}$ we will represent them as paths in the Pascal triangle. Given a standard tableau $\tg $ we use the convention that the $k$-th step in the path $p_\tg$ is drawn to the right (resp. left) if $k$ is located in the first (resp. second) component of $\tg$.

Each point in the Pascal triangle is determined by a \emph{level} and a \emph{weight}. The highest point in the  Pascal triangle is at level $0$. Levels increase by one from top to bottom. At level $n$ the leftmost point has weight $-n$.  Weights increase by two from left to right.

Clearly, two standard tableaux have the same shape if and only if their endpoints coincide. For this reason, we often identify this common endpoint with $\lamn$.

\begin{exa}\label{Exa tableau azul verde} \rm 
Let $\lamn=(1^{4},1^{8})\in P_1^2(12)$. On the left-hand side of (\ref{tableaux to paths}), we have two standard tableaux $\red\sg$ and $\blue \tg$ of shape $\lamn$. On the right-hand side of (\ref{tableaux to paths}), we have a  Pascal triangle.  For instance, the endpoint of both paths in the picture is at level $12$ and weight $-4$.   
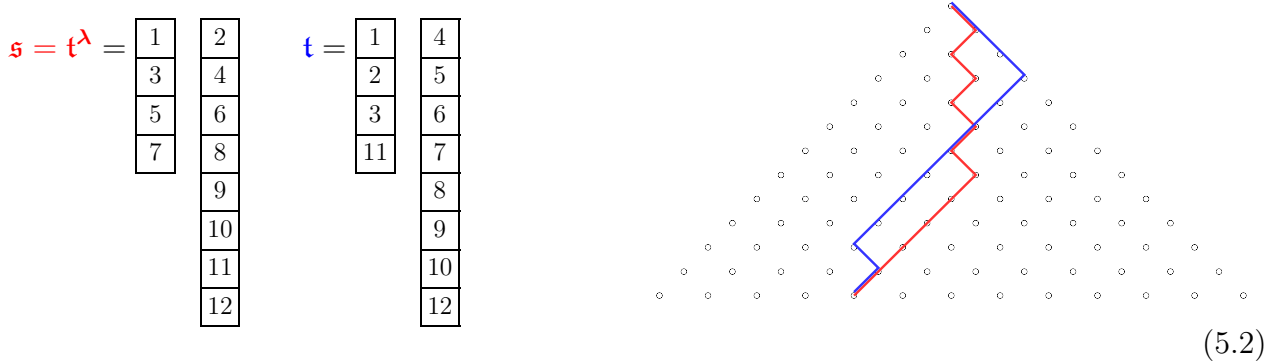
\begin{equation}\label{tableaux to paths}
	{\red \sg =\tg^{\lamn}} =\scalebox{.8}{\begin{ytableau}
1 \\
3\\
5\\
7
\end{ytableau}
\mbox{ }
\begin{ytableau}
2\\
4\\
6\\
8\\
9\\
10\\
11\\
12
\end{ytableau}}
\qquad
	{\blue \tg }=\scalebox{.8}{\begin{ytableau}
1\\
2\\
3\\
11
\end{ytableau}
\mbox{ }
\begin{ytableau}
4\\
5\\
6\\
7\\
8\\
9\\
10\\
12
\end{ytableau}}
\qquad
\qquad
\qquad
\scalebox{.8}{\UlSecA }
\end{equation} 
\end{exa}

It is clear that for an arbitrary $\lamn \in \ocmpardos$  the path $p_{\lamn}$ associated to the dominant tableau $\tg^{\lamn}$ is the one that first zig-zags (right-left) and then finishes  in a straight line ending up in $\lamn$. 

\medskip
Let  $\lamn \in \ocmpardos$ and $\tg\in \Std (\lamn)$. Suppose that $k$ and $k+1$ are located in different components of $\tg$. Let $\sg$ be the standard tableau obtained from $\tg$ by interchanging $k$ and $k+1$. In this case, we say that $\sg$ is obtained from $\tg$ by \emph{making a hook at level} $k$. This name is justified by the way in which paths $p_\tg$ and $p_\sg$ are related. Indeed,  if the $k$-th and $(k+1)$-th steps in $p_{\tg}$ form a subpath of the form  $\langle$  then $p_{\sg}$ is obtained from $p_{\tg}$ by making a hook at level $k$ if $p_{\sg}$ coincides with $p_{\tg}$ except in the subpath formed by the $k$-th and $(k+1)$-th steps where $\langle$  is  replaced by $\rangle$. One could say a similar thing for a subpath of the form $\rangle$.

\subsection{An algorithm to obtain $d_{\tg}$}
We recall that to each standard tableau we have associated an element  $d_\tg \in \mathfrak{S}_n$  defined by the rule $d_\tg\,  \tg^{\lamn}=\tg$. An interesting feature of the representation of $\tg$  as a path  in the Pascal triangle is that we can obtain in a beautiful way a reduced expression for $d_\tg$ as follows:
\begin{algo}  \label{Algoritmo} \rm 
In order to obtain $d_{\tg}$  we must follow the next steps:
	\begin{itemize}
	\item Draw the paths $p_\tg$ and $p_{\lamn}$ associated to $\tg$ and $\tg^{\lamn}$, respectively. 
	\item We will define a tuple of elements of the symmetric group $(d_\tg(0), d_\tg(1),\cdots, d_\tg(j))$ and a tuple of paths in the Pascal triangle $(p_0, p_1, \cdots, p_j)$ for some $j$ that will be defined in the next step of the algorithm.  Set $d_\tg(0)=\mathrm{id}\in S_n$ and $p_0:=p_{\lamn}$.  Assume that $p_{i-1}$ has been defined. Define $p_{i}$ as the path obtained from $p_{i-1}$ by making a hook at any level $k$ satisfying  that the area bounded by $p_\tg$ and $p_{i}$ is smaller than the one bounded by $p_{\tg} $ and $p_{i-1}$. Set $d_\tg(i)=s_k
	d_\tg(i-1)$. 
	\item Repeat the previous step until $p_j=p_\tg$, for some $j$. Then, $d_\tg =d_\tg (j)$.
\end{itemize}
\end{algo}

The fact that this algorithm indeed provides a reduced expression for $d_\tg$ is proved in \cite[Section 4]{plaza2014graded}. The second step of this algorithm can be performed in several  ways, so the reduced expression of $d_\tg$ obtained with  this process is not unique. However, in the case considered in this section, that is  $l=2$, it is easy to see that permutations $d_\tg$ are  $321$-avoiding  and therefore  the elements of our graded cellular basis of $\gblobdos$ do not depend on the particular choice of a reduced expression.

\begin{exa}\rm
	Keeping the  notation of Example \ref{Exa tableau azul verde} we can  obtain many results, for example 
	\begin{equation}
		d_{\blue \tg} = s_{10}s_9s_8s_7s_3s_4s_2  \qquad \mbox{ or } \qquad d_{\blue \tg} = s_3s_2s_4s_{10}s_9s_8s_7.		
 	\end{equation} 
 	The reduced expression of $d_{\blue \tg}$ on the left (resp. right) is obtained by applying the second step in the algorithm to the highest (resp. lowest) level possible.  
\end{exa}

\subsection{Description of standard tableaux with the same residue as $\lamn$}

For $l=2$ our alcove geometry lives in the real line and the set of hyperplanes reduces to the set of points 
$\{\hg_{1,2}^m:=\kappa_1-\kappa_2+me$, $m\in \mathbb{Z}\}$. In the Pascal triangle these hyperplanes are drawn as vertical lines located on weights corresponding to the aforementioned points. 

In this setting, the fundamental alcove corresponds to the one which contains the symmetry axis  of the Pascal triangle. We identify affine Weyl group generators $s$ and $t$ with the reflections through the hyperplanes that delimit  the fundamental alcove. By convention, the left hyperplane is identified with $s$ and the right hyperplane is identified with $t$. Therefore, alcoves to the left (resp. right) of the fundamental alcove are labelled by elements of the form $sts\ldots$ (resp. $tst\ldots$). 

\begin{exa}
	Let $\kappa = (1,4)$ and $e=5$. In the figure there is a Pascal triangle truncated at level $15$. We have drawn the six relevant hyperplanes and the labels of the alcoves appearing in the picture.  
	 \begin{equation}
	 	\scalebox{1.3}{\UlSecB }
	 \end{equation}
\end{exa}

Given a path $p$ we say that a subset of consecutive steps of $p$ is a \emph{wall to wall step} if these steps form a straight line between two consecutive hyperplanes.   In the following lemma (proved in \cite[Lemma 4.7]{plaza2013graded}), for a path $p(t)$  we  think of $t$ as a variable representing time.

\begin{lem} \label{lemma description tableaux with residue}
Let $\lamn \in \ocmpardos$ be regular and $\tg$ a standard tableau with $n$ boxes. We have that $\bi^\tg=\bi^{\lamn }$ if and only if  $p_\tg$  satisfies all of the following conditions:
\begin{enumerate}
	\item The paths $p_\tg$ and $p_{{\lamn}}$ coincide until the moment of the first contact of $p_{\lamn}$ with a hyperplane. 
	\item Then, until the moment when $p_{{\lamn}}$ touches for the last time a hyperplane, 
	  the path $p_\tg $ only makes  wall to wall steps (as many  as those made by $p_{\lamn}$).
	\item After $p_{{\lamn}}$ touches for the last time a hyperplane, $p_{\tg}$ is completed by a straight line until it stops (at level $n$). 
\end{enumerate} 
\end{lem}

\begin{exa}\rm
	Let $\kappa=(1,4)$ and $e=5$. In the figure  is shown a Pascal triangle truncated at level $30$. Let $\lamn=(1^2,1^{28})\in P_1^2(30)$. The path $p_{\lamn}$ is drawn in black. The $2^5$ descending paths following gray lines are those which correspond to standard tableaux  with  residue sequence equal to  $\bi^{\lamn}$. 
	
	\begin{equation}
	 	\scalebox{1.3}{\UlSecC }
	\end{equation}
\end{exa}

\vspace{0.2cm}

\subsection{Degree-zero subalgebra of $\gblobldos$}
Let $\lamn\in \ocmpardos$. The Pascal triangle also provides us an easy way to compute the degree of the standard tableaux with residue sequence $\bi^{\lamn}$. Let $\tg$ be a standard tableau with $n$ boxes such that $\bi^\tg=\bi^{\lamn}$. We define $\delta(\tg )=1$ (resp. $\delta(\tg )=0$) if the straight line completing $p_\tg$ after the wall to wall steps points towards (resp. away from) the symmetry axis of the Pascal triangle. We also define $w(\tg)$ as the number of wall to wall steps in $p_\tg$ crossing through the fundamental alcove. It was proved in  
\cite[Lemma 4.9]{plaza2013graded} that
\begin{equation} \label{equation degree}
	\deg(\tg) = w(\tg)+\delta (\tg).
\end{equation}
Consequently, $\deg (\tg) \geq 0$, for each $\tg$ with $\bi^\tg=\bi^{\lamn}$. An immediate consequence of this fact is the following 
\begin{lem}\label{pg}
Let $\lamn\in \ocmpardos$ be regular. Then the algebra $\gblobldos $ is positively graded. 
\end{lem}

By Lemma \ref{lemma properties positively graded algebras} we know that graded decomposition numbers of $\gblobldos $ (and therefore those of $\gblobdos $) belong to $\mathbb{N}[v]$. Therefore, in order to apply the process outlined in Section \ref{section algorithm to compute gdn blob} to calculate graded decomposition numbers, we only need to know the value of $d_{\lamn ,\mun}^p(0)$. To do this it is enough to compute  decomposition numbers of the  degree-zero component of $\gblobldos$, which we denote by  
$\gblobldos_0$.

 By Lemma \ref{pg}, a cellular basis for $\gblobldos_0$ is given by the elements $\psi_{\sg\tg}^{\mun}$ such that $\bi^\sg =\bi^\tg=\bi^{\lamn}$ and $ \deg(\sg) =\deg (\tg)=0$. In order to provide a better description of such a basis we need  more notation.

	Let $\lamn\in \ocmpardos$ regular. Given $\mun \in \ocmpardos $ we define 
	
	\begin{align}
		 \Std_{\lamn}^0 (\mun) &=\{ \tg \in \Std_{\lamn}(\mun) \mbox{ } | \mbox{ } \deg (\tg )=0 \}   \\
		P^0(\lamn)&=\{\mun \in \ocmpardos \mbox{ } |\mbox{ }  \Std_{\lamn}^0 (\mun) \neq \emptyset \}  \\
		\Std_{\lamn}^0  & = \bigcup_{\mun \in P^0(\lamn) } \Std_{\lamn}^0 (\mun).
	\end{align} 
We can now reformulate the preceding paragraph.
\begin{lem} \label{Lemma basis degree algebra}
	Let $\lamn\in \ocmpardos$ regular. Then, \begin{equation} \label{equation basis degree-zero}
	\{ \psi_{\sg\tg }^{\mun } \mbox{ }|\mbox{ } \mun \in P^0(\lamn ) \mbox{ and } \sg ,\tg \in \Std_{\lamn}^0(\mun )  \}
\end{equation}
is a cellular basis for $\gblobldos_0$.
\end{lem}

We now compute the dimension of $\gblobldos_0$. By (\ref{equation degree}) we have that  $\tg\in\Std_{\lamn}^0 $ if and only if the path $p_\tg$ associated to $\tg$ satisfies the three conditions in Lemma \ref{lemma description tableaux with residue} together the additional conditions: 
\begin{enumerate}
	\item Wall to wall steps of $p_{\tg} $ do not pass through the fundamental alcove.
	\item The final straight line that completes $p_\tg$, as in Lemma \ref{lemma description tableaux with residue} (3), moves away from the symmetry axis of the Pascal triangle.
\end{enumerate} 

\begin{exa}\label{example ultima section}\rm
	Let $\kappa=(1,4)$ and $e=5$. In the figure it is shown a Pascal triangle truncated at level $30$. Let $\lamn=(1^2,1^{28})\in P_1^2(30)$, so $w_{\lamn} =5_s$.  We have depicted the set of all paths that are associated to a standard tableau in $\Std_{\lamn}^0$.  It follows that $P^0(\lamn) =\{\lamn , \mun , \boldsymbol{\nu}\}$, where $\mun =(1^7,1^{23})$ and  $\boldsymbol{\nu}=(1^{12},1^{18})$. We have 	$w_{\mun}=3_s$, $w_{\boldsymbol{\nu}}=1_s$,
	$|\Std_{\lamn }^0(\lamn) |=1$, $|\Std_{\lamn }^0(\mun ) |=3$ and $|\Std_{\lamn }^0(\boldsymbol{\nu}) |=2$.  
	
	Consider for a moment, a curious coincidence of numbers (in the next lemma we give a full explanation). There are $3=\vert P^0(\lamn) \vert $ two-column partitions: $\lambda :=(1^4)$, $\mu :=(2^1,1^2)$ and $\nu :=(2^2)$. The number of standard tableaux of their shapes are: 1  of shape $\lambda$, 3 of shape $\mu $ and 2 of shape $\nu$, which is the same as $|\Std_{\lamn }^0(\lamn) |$, $|\Std_{\lamn }^0(\mun ) |$ and $|\Std_{\lamn }^0(\boldsymbol{\nu}) |$. 
		\begin{equation}
		\scalebox{.7}{\UlSecD } 
	\end{equation}
\end{exa}

As the previous example suggests, at this point we will need to work with two-column partitions. The reader should notice the difference between such partitions and one-column bipartitions, that is, the elements in $\ocmpardos$. Following the previous notation, we reserve bold symbols for one-column bipartitions and   normal (non-bold) symbols for two-column partitions. Given a positive integer $n$ we denote by $P_2(n)$ the set of all two-column partitions of $n$ and by $\Std(\lambda)$ the set of all standard tableaux of shape $\lambda$, for $\lambda \in P_2(n)$. We also stress that an element $\lambda \in P_2(n)$ is given by $\lambda =(2^j,1^{n-2j})$, for some $0\leq j \leq \lceil  \frac{n-1}{2} \rceil$. The following lemma is an easy consequence of the preceding discussion.

\begin{lem}  \label{lemma bijection degree-zero}
	Let $\lamn\in \ocmpardos$ regular. Assume that  $w_{\lamn}=k_s$, for some $k\geq 2$. Then,  
	\begin{enumerate}
		\item  $P^0(\lamn) =\{ \mun \in \ocmpardos \mbox{ } | \mbox{ } \mun \in W \cdot \lamn \mbox{ and }w_{\mun} =(k-2j)_s\mbox{ for some }  0\leq j \leq \lceil \frac{k-2}{2} \rceil  \}$.  
		\item There is a bijection $P^0(\lamn) \rightarrow P_2(k-1)$ given by $ \mun \rightarrow \mu$, where $w_{\mun } =(k-2j)_s$ and $\mu =(2^j,1^{k-1-2j})$.
		\item Given $\mun \in P^0 (\lamn)$ with $w_{\mun} =(k-2j)_s$, there is a bijection  $\Std_{\lamn}^0(\mun ) \rightarrow \Std(\mu )$ denoted by $\tg \rightarrow \tau_{\tg}$, determined by the following rule: An integer $1\leq j\leq k-1$ is located in the second column of $\tau_\tg$ if and only if the $j$-th wall to wall step of $p_{\tg} $ points towards the symmetry axis of the Pascal triangle. 
		\item The above statements remain true if we replace $s$ by $t$.
 	\end{enumerate}
\end{lem}

\begin{exa} \rm
Let us illustrate Lemma \ref{lemma bijection degree-zero}(3) with an example. We keep  parameters and notation of Example \ref{example ultima section}.  The bijections are given as follows:

$$ \begin{array}{rcrcr}
	\Std_{\lamn}^0(\lamn) \rightarrow \Std (\lambda)  & & \Std_{\lamn}^0(\mun) \rightarrow \Std (\mu ) & &  \Std_{\lamn}^0(\boldsymbol{\nu })\rightarrow  \Std (\nu) \\
	& & & & \\
\scalebox{.9}{\PasA }  \rightarrow  
{\begin{ytableau}
1 \\
2\\
3\\
4
\end{ytableau}}	 & &  \scalebox{.9}{\PasB }  \rightarrow  \begin{ytableau}
1 & 2\\
3 \\
4 
\end{ytableau} & &  \scalebox{.9}{\PasE } \rightarrow  \begin{ytableau}
1 & 2\\
3 & 4 
\end{ytableau}  \\
  & & \scalebox{.9}{\PasC }\rightarrow \begin{ytableau}
1 & 3\\
2 \\
4 
\end{ytableau} & & \scalebox{.9}{\PasF } \rightarrow  \begin{ytableau}
1 & 3\\
2 & 4 
\end{ytableau}  \\  
 & & & & \\
& &  \scalebox{.9}{\PasD }\rightarrow \begin{ytableau}
1 & 4\\
2 \\
3 
\end{ytableau}    & &
\end{array} 
$$
\end{exa}

\begin{cor} \label{Corollary dim agree}
Let $\lamn\in \ocmpardos$ regular with  $w_{\lamn}=k_s$ (or $w_{\lamn}=k_t$), for some $k\geq 2$. Then,  $\dim \gblobldos_0 = C_{k-1}$, where $C_{k-1}$ denotes the $(k-1)$-th Catalan number.
\end{cor}
\begin{dem}
The result follows by combining Lemma \ref{Lemma basis degree algebra}, Lemma \ref{lemma bijection degree-zero} and the well-known formula
\begin{equation}
	C_{k-1}=\sum_{j=0}^{\lceil \frac{k-2}{2} \rceil} |\Std (2^j,1^{k-1-2j})|^2.
\end{equation}
\end{dem}

\begin{remark} \rm
For $\lamn$ regular  and  $w_{\lamn} =k_s$ (or $w_{\lamn} =k_t$ ), the condition $k\geq 2$ in the previous results is not restrictive since if $0\leq k <2$ it is easy to see that $\Std_{\lamn}^0 =\{\tg^{\lamn}\}$, and therefore $\gblobldos_0 \cong \mathbb{F}_p$. 
	
\end{remark}

\subsection{$\gblobldos_0$ and the Temperley-Lieb algebra}
We will now focus on determining a presentation for $\gblobldos_0$. We start by defining certain elements  $ \UB_j^{\lamn} \in \gblobldos$.  We will show later that these elements  actually generate $\gblobldos_0$. In order to introduce these elements we first need the following  

\begin{defi} \rm
	Let $\lamn = (1^{\lambda_1}, 1^{\lambda_2})\in \ocmpardos$ regular with  $w_{\lamn}=k_s$ (or $w_{\lamn}=k_t$), for some $k\geq 1$. We set $m_{\lamn}=\min \{ \lambda_1 , \lambda_2  \} $ and define 
	\begin{equation}
		f_{\lamn }= \left\{\ \begin{array}{rl}
			2m_{\lamn} - (\kappa_1-\kappa_2),&  \mbox{if } m_{\lamn} =\lambda^1; \\
			2m_{\lamn}+ (\kappa_1-\kappa_2) +e  ,&  \mbox{if } m_{\lamn} =\lambda^2.
		\end{array}  \right. 
	\end{equation}   
	The number $f_{\lamn}$ corresponds to the level when for first time $p_{\lamn} $ touches a hyperplane. 
	Given an integer $1\leq j < k-1$ we define $\underline{j}=f_{\lamn} +je $ and
$\UB_j^{\lamn }\in \gblobldos$ as
$$(\psi_{\underline{j}})(\psi_{\underline{j}-1}\psi_{\underline{j}+1} )(\psi_{\underline{j}-2}\psi_{\underline{j}}\psi_{\underline{j}+2} )\cdots (\psi_{\underline{j}-e+1}\psi_{\underline{j}-e+3}\cdots \psi_{\underline{j}+e-3}\psi_{\underline{j}+e-1} ) \cdots (\psi_{\underline{j}-2}\psi_{\underline{j}}\psi_{\underline{j}+2} )(\psi_{\underline{j}-1}\psi_{\underline{j}+1} )(\psi_{\underline{j}})e(\lamn).$$
\end{defi}

We refer to elements $\UB_j^{\lamn}$ as the diamond of weight $\lamn$ at position $\underline{j}$. The following example justifies such a name, as well as the fact that $\UB_j^{\lamn }\in \gblobldos$. 

\begin{exa} \rm
	Let $\kappa=(1,4)$ and $e=5$. Let  $\lamn =(1^2,1^{28}) \in P_1^2(30) $. In this case we  have $w_{\lamn} =ststs= 5_s$ and $f_{\lamn} =7$. Then, we have three underlined integers given by $\underline{1}=12$, $\underline{2}=17$ and $\underline{3}=22$.  The elements $\UB_j^{\lamn}$ are depicted below.
	
$$ \UB_1^{\lamn} = \UlSecGa \qquad \UB_2^{\lamn} = \UlSecGb  $$	

$$ \UB_3^{\lamn}=  \UlSecGc $$ 
\end{exa}

Our next goal is to show that the elements  $\UB_j^{\lamn}$ generate $\gblobldos_0$. We should start by proving that $\UB_j^{\lamn}$ actually  belongs to  $\gblobldos_0$, that is, $\deg  \UB_j^{\lamn} =0$. However, we prefer to begin by  finding  out the relations satisfied by such elements. The fact that $\deg  \UB_j^{\lamn} =0$ will follow from such relations. Before embarking in the proof of the relations satisfied by the $\UB_j^{\lamn}$'s we need two lemmas. 

\begin{lem}
	In $\gblobdos$  the diagrammatic relation (\ref{Trenza diagramtical}) reduces to 
	 \begin{equation}  \label{equation Trenza simplificada A}
	 	\TrenzaSimplificada \quad  =-\quad \TrenzaSimplificadaD  \qquad \qquad \TrenzaSimplificadaB \quad =0, 
	 \end{equation}
	 
	 \vspace{.3cm}
	 
	 \begin{equation} \label{equation Trenza simplificada B}
	 \TrenzaSimplificadaC \quad = - \quad \TrenzaSimplificadaE \qquad \qquad   	\TrenzaSimplificadaA \quad =0.
	 \end{equation}
	 
\end{lem}

\begin{dem}
	We will only prove (\ref{equation Trenza simplificada A}), (\ref{equation Trenza simplificada B}) is treated similarly. It follows by  (\ref{Trenza diagramtical}) that both relations in (\ref{equation Trenza simplificada A}) are equivalent. We will show the equation on the right. By Proposition \ref{Theorem Graded cellular basis blob} we know that $e(\bi)\neq 0 $ if and only if $\bi =\bi^\tg$, for some standard tableau $\tg$. The relevant tableaux in $\gblobdos$ have two components. Suppose that  $\bi =(\ldots, i,i,i+1, \ldots)\in I_e^n$. Then, $\bi \neq \bi^\tg$ for all standard tableau $\tg$, since if two $i$'s occur consecutively  in the residue sequence of a standard tableau the next value in the sequence must be $i-1$. We conclude that $e(\bi) =0$. We notice that a subsequence of type $(i,i,i+1)$ appears in the middle of the diagram on the right-hand side of (\ref{equation Trenza simplificada A}). Then, such a diagram is equal to $0$. 
	\end{dem}

\begin{lem}  \label{Lemma doble cruce con dot al medio}
	In $\gblob$ we have 
	\begin{equation}
		\DobleCruceconDot =\Cruce  \qquad \qquad \qquad    \qquad \DobleCruceconDotA = - \Cruce
	\end{equation}
\end{lem}

\begin{dem}
	The result follows by a direct application of (\ref{Diagrammatic Crossing dot}) and (\ref{Diagrammatic Quadratic Relation}). 
\end{dem}

\begin{pro}  \label{Proposition relations Us}
Let $\lamn \in \ocmpardos$ regular with  $w_{\lamn}=k_s$ (or $w_{\lamn}=k_t$), for some $k\geq 2$. The elements $\UB_1^{\lamn}, \UB_2^{\lamn}, \ldots , \UB_{k-2}^{\lamn }$ satisfy 
\begin{align}
\label{equation UA} \UB_i^{\lamn }\UB_j^{\lamn }   &= \UB_j^{\lamn }\UB_i^{\lamn },& \quad \text{ if } |i-j|>1; \\
\label{equation UC} \left(\UB_i^{\lamn }\right)^2  &=(-1)^{e-1}2\UB_i^{\lamn },  &  \\
\label{equation UB} \UB_i^{\lamn }\UB_j^{\lamn } \UB_i^{\lamn } &  = \UB_j^{\lamn }, & \quad  \mbox{ if } |i-j|=1;
\end{align}
\end{pro}

\begin{dem}
	Suppose that $\psi_s$ and $\psi_r $ occur in   $\UB_i^{\lamn }$ and $\UB_j^{\lamn }$, respectively. We note that if $|i-j|>1$ then $|s-r|>1$. By (\ref{kl8}) we conclude that $\UB_i^{\lamn }$ and $\UB_j^{\lamn }$ commute if $|i-j|>1$. This proves (\ref{equation UA}). 
	
	We will now focus on proving (\ref{equation UC}).  We call a subdiagram of the form 
	\begin{equation}
		\scalebox{.5}{\DobleCruce } 
	\end{equation}

a \emph{Double Crossing}. Such a diagram  can be \emph{replaced}	 by using  (\ref{Diagrammatic Quadratic Relation}). In accordance with such a relation we classify double crossing in four types: 
\begin{enumerate}
	\item (ZDC) Zero Double Crossings ($i=j$);
	\item (DDC) Distant Double Crossings ($|i-j|>1$);
	\item (IADC) Increasing Adjacent Double Crossings ($j=i+1$);
	\item (DADC) Decreasing Adjacent Double Crossings ($j=i-1$).
\end{enumerate}

 The diagram associated to $(\UB_i^{\lamn })^2$ is formed by two diamonds, one above the other, as is illustrated in the left-hand side of  (\ref{Diamond A}).  Such a diagram assumes that $e=8$. For brevity, we have only drawn the part of the diagram where intersections occur. We have also omitted the residues on the bottom. We only need to know that they form a sequence of length $2e$ of the form 
	$$(a,a-1, \ldots , a-(e-1), a,a-1, \ldots ,  a-(e-1) ).$$
The idea is to apply the relations to "disarm"  intersections from the center to the ends. We start by applying (\ref{Diagrammatic Quadratic Relation}) to the IADC just in the center of the diagram representing 
	$(\UB_i^{\lamn })^2$. In this way, we obtain the second equality in (\ref{Diamond A}). We remark that  by (\ref{Diagrammatic Crossing dot}) a dot can freely pass through an intersection as long as  the residues involved in such a intersection are distinct. We use this fact $e-2$ times to move the dots up in the previous diagrams. Then, we can apply (\ref{Diagrammatic Quadratic Relation}) to replace in both diagrams $(e-2)(e-1)/2-1$ DDC's by straight lines. This shows 
	the third equation in (\ref{Diamond A}).
	\begin{equation} \label{Diamond A}
	(\UB_i^{\lamn })^2 =\DiamondA = \DimanondARightPoint -\DimanondALeftPoint    =	\DiamondBA  - \DiamondB
\end{equation}
	
	Let us now concentrate on the sub-diagrams in the middle of the rightmost diagrams in (\ref{Diamond A}). Concretely, we are interested in 
\begin{equation}\label{Diamond D}
	\scalebox{.8}{\DiamondCA }\qquad \mbox{ and } \qquad \scalebox{.8}{\DiamondC }.
\end{equation}

	We notice that in both diagrams the bottom sequence is 
$$ (a,a,a-1,a-1,\ldots , a-(e-1),a-(e-1)).   $$

Let us focus on the diagram on the left of (\ref{Diamond D}). By applying (\ref{Diagrammatic Quadratic Relation}) to the leftmost DADC in such a diagram we obtain the first equality of (\ref{Diamond E}).
	Once again by (\ref{Diagrammatic Quadratic Relation}) we obtain the second equality in (\ref{Diamond E}), since the highlighted  double-crossing is of type ZDC.
\begin{equation}\label{Diamond E}
	\scalebox{.8}{\DiamondCA} = \scalebox{.8}{\DiamondDA } - \scalebox{.8}{\DiamondD } =  \scalebox{.8}{\DiamondDA }.
\end{equation}
By repeating the same argument and then by applying  Lemma \ref{Lemma doble cruce con dot al medio} we obtain 
\begin{equation}\label{Diamond G}
	\scalebox{.8}{\DiamondCA } =  \scalebox{.8}{\DiamondE } = (-1)^{e-1}  \scalebox{.8}{\DiamondFC } .
\end{equation}
Similar computations reveal that  	
\begin{equation} \label{Diamond J}
	\scalebox{.8}{\DiamondC }=(-1)^e \scalebox{.8}{ \DiamondFC }.
\end{equation} 
Finally,  by combining (\ref{Diamond A}), (\ref{Diamond G}) and (\ref{Diamond J}) we have

\begin{equation}
	(\UB_i^{\lamn })^2  =  (-1)^{e-1}\, \DiamondFinal - (-1)^{e} \, \DiamondFinal = (-1)^{e-1}2 \UB_i^{\lamn }. 
\end{equation}

We now prove relation (\ref{equation UB}). We focus on the equation $\UB_i^{\lamn}\UB_{i+1}^{\lamn}\UB_i^{\lamn}=\UB_i^{\lamn} $. The relation $\UB_i^{\lamn}\UB_{i-1}^{\lamn}\UB_i^{\lamn}=\UB_i^{\lamn} $ is treated similarly. The diagram representing $\UB_i^{\lamn}\UB_{i+1}^{\lamn}\UB_i^{\lamn}$ is formed by three diamonds, as is illustrated in the left-hand side of (\ref{blobTl}). Let us refer to them as the top diamond, the middle diamond and the bottom diamond. As before, in that diagram we have only drawn the part where intersections occur and we have assumed that $e=8$. We  have also omitted the sequence of residues at the bottom of such a diagram. The only thing we need to know about such a sequence is that it is of length $3e$ and of the form 
	$$(a,a-1, \ldots , a-(e-1), a,a-1, \ldots ,  a-(e-1) ,  a,a-1, \ldots ,  a-(e-1) ).$$
The intersection at the bottom corner of the top diamond, the intersection at the left corner of the middle diamond and the intersection at the top corner of the bottom diamond form a \emph{braid} of the form 
\begin{equation} \label{Dibujo proof three diamonds}
	\scalebox{.7}{\TrenzaSimplificadaX }.
\end{equation}

Thus, we can use (\ref{equation Trenza simplificada A}) to reduce it to three straight lines multiplied by a factor $(-1)$. Once the above braid disappeared, we have another braid of the same type located at the right of the previous one. By repeating  the  previous step $e$-times we obtain the second equality in (\ref{blobTl}). By applying (\ref{Diagrammatic Quadratic Relation}) to the $e(e-3)$ DDC's appearing in the middle of the resultant diagram we obtain the third equality in (\ref{blobTl}).

\begin{equation}  \label{blobTl}
\UB_i^{\lamn}\UB_{i+1}^{\lamn}\UB_i^{\lamn}=	\scalebox{.6}{\dibujoTL } = (-1)^e\mbox{ }  \scalebox{.6}{\dibujoTLA } =  \mbox{ }(-1)^e\scalebox{.6}{\dibujoTLB }.
\end{equation}

We now concentrate on the rightmost diagram in  (\ref{blobTl}). The $e$ double crossing in the middle of such a diagram are of type DADC. Once we apply (\ref{Diagrammatic Quadratic Relation}) to replace them by straight lines (with some dots) we have double crossings of type ZDC. Then, only one diagram does not vanish, the one which has a dot in the middle of each double crossing of type ZDC. In this way we derive the first equality in (\ref{blobTlA}). The second equality is  a consequence of Lemma \ref{Lemma doble cruce con dot al medio}.   

\begin{equation} \label{blobTlA}
	\scalebox{.7}{\dibujoTLB } =   \scalebox{.7}{\dibujoTLC }  =  (-1)^{e} \scalebox{.7}{\dibujoTLD }
\end{equation}

Finally,  by combining (\ref{blobTl}) and (\ref{blobTlA}) we obtain 

\begin{equation}
	\UB_i^{\lamn}\UB_{i+1}^{\lamn}\UB_i^{\lamn}= (-1)^{2e} \scalebox{.7}{\dibujoTLD } =  \UB_i^{\lamn}.
\end{equation}
\end{dem}

\begin{cor}
	The elements $\UB_i^{\lamn}$ belong to $\gblobldos_0$. 
\end{cor}

\begin{dem}
	The result follows by taking the degree on both sides of (\ref{equation UC}). 
\end{dem}

\medskip
Given $\mu \in P_2(n)$ we define $\tau^\mu \in \Std (\mu)$ to be the unique standard tableau in which the numbers $1,2,\ldots , n$  are filled in increasingly along the rows from top to bottom. 

\begin{lem}  \label{Lemma genera A}
Let $\lamn \in \ocmpardos$ regular with $w_{\lamn}=k_s$ (or $k_t$). Let $\mun \in P^0 (\lamn)$ with $w_{\mun} =(k-2j)_s$, for some $1\leq j \leq \lceil \frac{k-2}{2} \rceil $. Let $\mu =(2^{j}, 1^{k-1-2j})$. We denote  by $\tg \in \Std_{\lamn}^0(\mun )$ the tableau which is mapped  to $\tau^{\mu}$ via the bijection in Lemma \ref{lemma bijection degree-zero}(3). Then, 
\begin{equation}\label{equation genera}
	 \psi_{\tg\tg}^{\mun} =(-1)^{je} \UB_1^{\lamn}\UB_3^{\lamn} \cdots \UB_{2j-1}^{\lamn}.  
\end{equation}

\end{lem}

\begin{dem}
We first notice that the restriction $1\leq j$ is made in order to discard the case $\mun =\lamn $. In order to know what $\psi_{\tg\tg}^{\mun}$ is we need to know a reduced expression for $d_{\tg}$. We obtain such a  expression  by using Algorithm \ref{Algoritmo}. To do this, we have to draw the paths $p_{\mun} $ and $p_{\tg} $  associated to 
$\tg^{\mun}$ and $\tg$, respectively.	
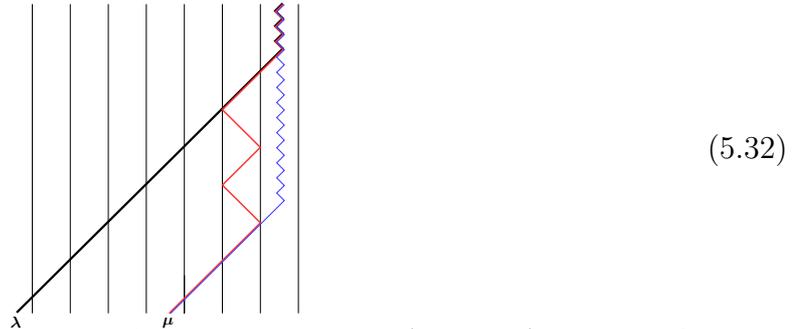
\begin{equation} \label{genera A}
	\scalebox{.5}{	\Genera }
	\end{equation}
	
	 In (\ref{genera A}) we have illustrated the situation, for $w_{\lamn }=7_s$, $w_{\mun } =(7-2\cdot j)_s=3_s$ and $e=5$. So that $j=2$. In such a figure  $p_{\mun} $ and $p_{\tg} $ correspond to the blue path and the red path, respectively. As a reference, we have also drawn the path $p_{\lamn}$ associated to $\tg^{\lamn}$ (the black path). We can now use  Algorithm \ref{Algoritmo} to obtain the diagram associated to $\psi_{\tg\tg}^{\mun}$. For instance, in the situation illustrated in (\ref{genera A}) we obtain
	
	\begin{equation} \label{genera B}
	\psi_{\tg\tg}^{\mun}= 	\scalebox{.5}{\GeneraB }. 
	\end{equation}
	Some remarks are in order. We have only drawn the part of $\psi_{\tg\tg}^{\mun}$ where intersections occur, so there are both  non-drawn straight lines to the left and to the right. We have also omitted the residues on the bottom of the diagram. This is not important, actually, we only need to know that the sequence below such a diagram is of the form $(a,a-1,a-2, \ldots)$. We recall the numbers in such a sequence are read modulo $e$. Finally, green lines do not mean anything. They are there in order to indicate a part of the diagram. 	
	
	\medskip
	By using the same arguments as the ones utilized in the proof of Proposition \ref{Proposition relations Us} we can reduce the region delimited by the green lines to 
	\begin{equation} \label{genera C}
		(-1)^{je} \quad \scalebox{.7}{\GeneraC } 
	\end{equation}
	Then, (\ref{equation genera}) follows by combining (\ref{genera B}) and (\ref{genera C}). The general case follows essentially in the same way. Of course, the  diagram associated to $\psi_{\tg\tg}^{\mun}$ will change.  For example, the 
height of the region that corresponds to the one bounded by the green lines depends on the distance between the symmetry axis  of the Pascal triangle and the relevant wall of the fundamental alcove.  Concretely, the further away the wall of the symmetry axis, the higher the region bounded by green lines will be. In any case, we still have that such a region reduces to a horizontal line of intersections, as in (\ref{genera C}), proving the result.
\end{dem}

\begin{cor} \label{coro genera A}
	Let $\lamn \in \ocmpardos$ regular with $w_{\lamn}=k_s$ (or $k_t$). Let $\mun \in P^0 (\lamn)$ with $w_{\mun} =(k-2j)_s$, for some $1\leq j \leq \lceil \frac{k-2}{2} \rceil $. Suppose that $\ug ,\vg \in \Std_{\lamn}^0(\mun ) $. Then, $\psi_{\ug\vg}^{\mun }$ can be written as a product of elements $\UB_j^{\lamn }$. Consequently, the set $\{e(\lamn), \UB_1^{\lamn} , \ldots ,\UB_{k-2}^{\lamn }\}$ generates $\gblobldos_0$. 
\end{cor}

\begin{dem}
Let $\mu = (2^j,1^{k-1-2j})$ be the two-column partition associated to $\mun$ via Lemma \ref{lemma bijection degree-zero}(2). We also denote by $\tg \in  \Std_{\lamn}^0(\mun )$ the standard tableaux which is mapped to $\tau^\lambda $ via  Lemma \ref{lemma bijection degree-zero}(3). A direct application of Algorithm \ref{Algoritmo} reveals that  
\begin{equation}
	\psi_{\ug\vg}^{\mun}  = A\psi_{\tg\tg}^{\mun} B, 
\end{equation}
where $A$ and $B$ are a product of elements $\UB_j^{\lamn}$. Then, $\psi_{\ug\vg}^{\mun }$ can be written as a product of elements $\UB_j^{\lamn }$ by Lemma \ref{Lemma genera A}. 

We can rephrase the previous paragraph by saying that each one of the elements in the cellular basis of $\gblobldos_0$ of Lemma \ref{Lemma basis degree algebra} (with the exception of $\psi_{\tg^{\lamn}\tg^{\lamn}}^{\lamn }$) belongs to the subalgebra generated by the $\UB_{i}^{\lamn}$'s. The last claim in the corollary is then a consequence of the equality $e(\lamn ) =\psi_{\tg^{\lamn}\tg^{\lamn}}^{\lamn }$. 
\end{dem}

\begin{exa} \rm 
	Let us illustrate Corollary \ref{coro genera A}. We take $\kappa =(1,4)$ and $e=5$. Let $\lamn =(1^3,1^{38})$ and $\mun =(1^{13},1^{18})$. We notice that $w_{\lamn} = 7_s$ and $w_{\mun } = 3_s$. It follows that the two-column partition corresponding to $\mun $ is $\mu =(2^2, 1^2)$. In both pictures of (\ref{genera doble}) we have drawn the paths associated to $\tg^{\lamn}$ (black),  $\tg^{\mun} $ (blue) and $\tg$ (red), where $\tg  \in \Std_{\lamn}^0(\mun )$  is the  standard tableau which is mapped to $\tau^\mu $. 
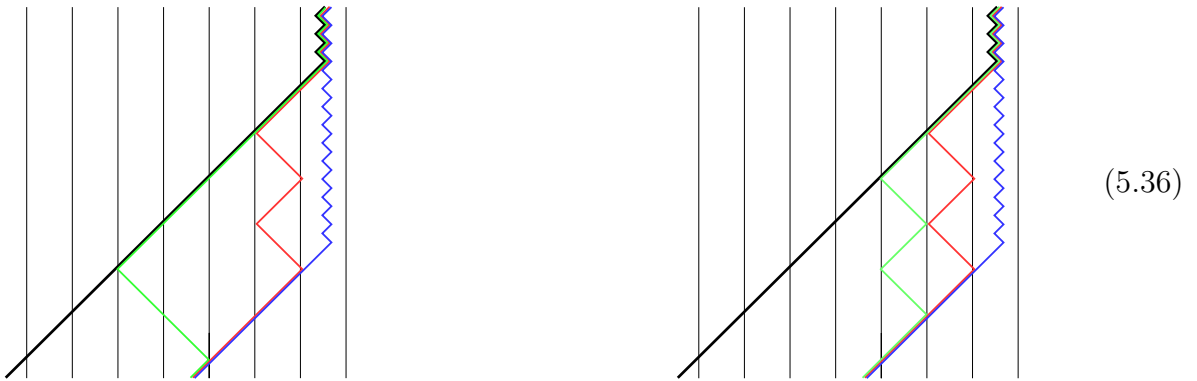
	\begin{equation} \label{genera doble}
		\scalebox{.3}{\GeneraD }   \scalebox{.3}{\GeneraE }  
	\end{equation} 
	We  denote by $\ug \in \Std_{\lamn}^0(\mun )$  (resp. $\vg$) the standard tableaux associated to the green path of  the left (resp. right) picture. The key point here is that in both cases, the red path is in between the blue and green paths. By performing Algorithm \ref{Algoritmo} and ignoring the straight lines to the left and to the right in $\psi_{\ug \vg}^{\mun}$ and also its residue sequence on the bottom, we obtain
	
	\begin{equation}
	\psi_{\ug \vg}^{\mun} = 	\scalebox{.5}{\GeneraF } = \UB_4^{\lamn} \UB_3^{\lamn} \UB_5^{\lamn} \UB_2^{\lamn}\UB_4^{\lamn} \psi_{\tg \tg}^{\mun} \UB_2^{\lamn}\UB_4^{\lamn} . 
	\end{equation}
	\end{exa}

We are now in position  to demonstrate that $\gblobldos_0$ is isomorphic to the Temperley-Lieb algebra \cite{temperley1971relations}. 
 
\begin{defi}
Let $n$ be a positive integer and $q\in \mathbb{F}_p^{\times}$. The Temperley-Lieb algebra $Tl_n^p(q)$ is the 	$\mathbb{F}_p$-algebra associative  on the  generators  $1,U_1$, $\ldots , U_{n-1} $ subject to the following relations
\begin{align}
\label{relation tl dos} U_iU_{j} U_i & = U_i, &   \mbox{ if } \vert i-j \vert =1  ; \\
\label{relation tl uno} U_i^2	& = -(q+q^{-1})U_i,      \\
\label{relation tl tres} U_iU_j &=U_jU_i,  & \mbox{ if }  \vert i-j \vert >1 .
\end{align}
\end{defi}

\begin{teo} \label{teo iso tl blob idempotent in degree zero}
Let $\lamn \in \ocmpardos$ regular with $w_{\lamn}=k_s$ (or $k_t$), for some $k\geq 2$. 
%
Then, there exists an isomorphism of $\mathbb{F}_p$-algebras  $\Phi : Tl_{k-1}^p(1) \mapsto \gblobldos_0$ determined by $1\mapsto e(\lamn)$ and $U_i \mapsto \UB_i^{\lamn}$ (resp. $U_i \mapsto -\UB_i^{\lamn}$ ) if $e$ is even (resp. if $e$ is odd). 
\end{teo}

\begin{dem}
	By combining Proposition \ref{Proposition relations Us}  and Corollary \ref{coro genera A}  we have that $\Phi$ is a surjective homomorphism. On the other hand,  it is well-known that the dimension 	of  $Tl_{k-1}^p(1)$ is the $(k-1)$-th Catalan number $C_{k-1}$. The result now follows by Corollary \ref{Corollary dim agree}.
\end{dem}

\medskip
We recall that the  function $f_p:\ZZ_{\geq 0} \times \ZZ_{\geq 0} \rightarrow \{0,1\}$ was defined in (\ref{function p contains}).

\begin{lem}\label{lemma dec temperley-Lieb}
The Temperley-Lieb algebra $Tl_n^p(1)$ is a cellular algebra with cell and simple modules parameterized  by $P_2(n)$. Let  $\lambda=(2^{j},1^{n-2j})$, $\mu = (2^k,1^{n-2k}) \in P_2(n)$. The corresponding decomposition number, $d_{\lambda , \mu }^{p,Tl}$, is given by 
\begin{equation}  \label{decom temperley-Lieb}
d_{\lambda , \mu }^{p,Tl} = f_p(n-2k , j-k). 
\end{equation} 
In particular, $d_{\lambda , \mu }^{p,Tl}=0 \mbox{ or }1$, for all $\lambda , \mu \in P_2(n)$.
\end{lem}

\begin{dem}
	The cellularity of $Tl_n^p( 1)$ is established in \cite{graham1996cellular}. The formula in (\ref{decom temperley-Lieb}) is implicit in \cite[Theorem 24.15]{james1987representation}, using the fact that $Tl_n^p( 1)$ can be realized as a quotient of the group algebra $\mathbb{F}_p \mathfrak{S}_n$. It was also proved in   \cite[Proposition 4.5]{cox2003blob}. We remark that with respect to both references we are working in the transpose setting.
	\end{dem}

\begin{teo} \label{teo equality at evaluating at zero}
	Let $\lamn \in \ocmpardos$ regular. Suppose that $\mun , \nun \in P^0(\lamn)$. Let $\mu $ and $\nu $ be the two-column partitions associated to $\mun$ and $\nun$ according to Lemma \ref{lemma bijection degree-zero}(2). Then, 
	\begin{equation}
		d_{\mun, \nun}^{p,\lamn} (0)= d_{\mu , \nu}^{p,Tl} = h_{w_{\mun},w_{\nun}}^p(0). 
	\end{equation}
\end{teo}

\begin{dem}
	The first equality follows by Theorem \ref{teo iso tl blob idempotent in degree zero}. The second equality follows by combining Lemma \ref{lemma dec temperley-Lieb} and equation (\ref{equation p-poly evaluated at zero}). 
\end{dem}

\begin{cor}
	In type $\tilde{A}_1$,  Blob vs Soergel Conjecture and Blob vs Light Leaves Conjecture hold. 
\end{cor}
	
\begin{dem}
On the one hand, the algorithm outlined in Section \ref{section an important formula}  gives us $p$-Kazhdan-Lusztig polynomials and graded ranks of intersection forms in type $\tilde{A}_1$. On the other hand, the algorithm outlined in Section \ref{section algorithm to compute gdn blob}
gives us graded decomposition numbers and graded dimensions of simple modules of $\gblobldos$. By Theorem \ref{teo dim cell module hecke} and Theorem \ref{teo equality at evaluating at zero} both algorithms produce the same polynomials.
\end{dem}

 \bibliographystyle{myalpha}
\bibliography{gen}

\def\cprime{$'$} \def\cprime{$'$} \def\cprime{$'$}
\begin{thebibliography}{AMRW17}

\bibitem[AMRW17]{achar2017koszul}
P.~Achar, S.~Makisumi, S.~Riche, and G.~Williamson.
\newblock Koszul duality for {K}ac-{M}oody groups and characters of tilting
  modules.
\newblock {\em arXiv preprint arXiv:1706.00183}, 2017.

\bibitem[BC17]{bowman2017modular}
C.~Bowman and A.~Cox.
\newblock Modular decomposition numbers of cyclotomic {H}ecke and diagrammatic
  {C}herednik algebras: A path theoretic approach.
\newblock {\em arXiv preprint arXiv:1706.07128}, 2017.

\bibitem[BCS17]{bowman2017family}
C.~Bowman, A.~Cox, and L.~Speyer.
\newblock A family of graded decomposition numbers for diagrammatic {C}herednik
  algebras.
\newblock {\em International Mathematics Research Notices}, 2017(9):2686--2734,
  2017.

\bibitem[BK09]{BrKl}
J.~Brundan and A.~Kleshchev.
\newblock Blocks of cyclotomic {H}ecke algebras and {K}hovanov-{L}auda
  algebras.
\newblock {\em Invent. Math.}, 178(3):451--484, 2009.

\bibitem[Bow17]{bowman2017many}
C.~Bowman.
\newblock The many graded cellular bases of {H}ecke algebras.
\newblock {\em arXiv preprint arXiv:1702.06579v6}, 2017.

\bibitem[CGM03]{cox2003blob}
A.~Cox, J.~Graham, and P.~Martin.
\newblock The blob algebra in positive characteristic.
\newblock {\em Journal of Algebra}, 266(2):584--635, 2003.

\bibitem[DJM98]{dipper1998cyclotomic}
R.~Dipper, G.~James, and A.~Mathas.
\newblock Cyclotomic q--schur algebras.
\newblock {\em Mathematische Zeitschrift}, 229(3):385--416, 1998.

\bibitem[Don98]{donkin1998q}
S.~Donkin.
\newblock {\em The q-Schur algebra}, volume 253.
\newblock Cambridge University Press, 1998.

\bibitem[EK10]{EKh}
B.~{Elias} and M.~{Khovanov}.
\newblock {Diagrammatics for Soergel categories.}
\newblock {\em {Int. J. Math. Math. Sci.}}, 2010:58, 2010.

\bibitem[EL17a]{ElLi}
B.~Elias and N.~Libedinsky.
\newblock Indecomposable {S}oergel bimodules for universal {C}oxeter groups.
\newblock {\em Trans. Amer. Math. Soc.}, 369(6):3883--3910, 2017.
\newblock With an appendix by Ben Webster.

\bibitem[EL17b]{elias2017modular}
B.~Elias and I.~Losev.
\newblock Modular representation theory in type {A} via {S}oergel bimodules.
\newblock {\em arXiv preprint arXiv:1701.00560}, 2017.

\bibitem[Eli16]{EDC}
B.~Elias.
\newblock The two-color {S}oergel calculus.
\newblock {\em Compos. Math.}, 152(2):327--398, 2016.

\bibitem[EW14]{EW2}
B.~Elias and G.~Williamson.
\newblock The {H}odge theory of {S}oergel bimodules.
\newblock {\em Ann. of Math. (2)}, 180(3):1089--1136, 2014.

\bibitem[EW16]{EW}
B.~Elias and G.~Williamson.
\newblock Soergel calculus.
\newblock {\em Represent. Theory}, 20:295--374, 2016.

\bibitem[GJSV13]{gainutdinov2013physical}
A.~M. Gainutdinov, J.~L. Jacobsen, H.~Saleur, and R.~Vasseur.
\newblock A physical approach to the classification of indecomposable
  {V}irasoro representations from the blob algebra.
\newblock {\em Nuclear Physics B}, 873(3):614--681, 2013.

\bibitem[GL96]{graham1996cellular}
J.~J. Graham and G.~I. Lehrer.
\newblock Cellular algebras.
\newblock {\em Inventiones mathematicae}, 123(1):1--34, 1996.

\bibitem[Gro]{Gr}
I.~Grojnowski.
\newblock Affine $sl_p$ controls the representation theory of the symmetric
  group and related {H}ecke algebras.
\newblock Preprint.
\newblock arXiv:9907129.

\bibitem[HM10]{hu2010graded}
J.~Hu and A.~Mathas.
\newblock Graded cellular bases for the cyclotomic
  {K}hovanov--{L}auda--{R}ouquier algebras of type {A}.
\newblock {\em Advances in Mathematics}, 225(2):598--642, 2010.

\bibitem[Hum92]{humphreys1992reflection}
J.~E. Humphreys.
\newblock {\em Reflection groups and {C}oxeter groups}, volume~29.
\newblock Cambridge {U}niversity {P}ress, 1992.

\bibitem[Jam87]{james1987representation}
G.~D. James.
\newblock {\em The representation theory of the symmetric groups}.
\newblock Berlin, 1987.

\bibitem[JW17]{JW}
L.~T. Jensen and G.~Williamson.
\newblock The {$p$}-canonical basis for {H}ecke algebras.
\newblock In {\em Categorification and higher representation theory}, volume
  683 of {\em Contemp. Math.}, pages 333--361. Amer. Math. Soc., Providence,
  RI, 2017.

\bibitem[KL79]{KLPolynomials}
D.~Kazhdan and G.~Lusztig.
\newblock Representations of {C}oxeter groups and {H}ecke algebras.
\newblock {\em Invent. Math.}, 53(2):165--184, 1979.

\bibitem[KL09]{KhLa}
M.~Khovanov and A.~D. Lauda.
\newblock A diagrammatic approach to categorification of quantum groups. {I}.
\newblock {\em Represent. Theory}, 13:309--347, 2009.

\bibitem[Lib08a]{Lieq}
N.~Libedinsky.
\newblock \'{E}quivalences entre conjectures de {S}oergel.
\newblock {\em J. Algebra}, 320(7):2695--2705, 2008.

\bibitem[Lib08b]{LLL}
N.~Libedinsky.
\newblock Sur la cat\'egorie des bimodules de {S}oergel.
\newblock {\em J. Algebra}, 320(7):2675--2694, 2008.

\bibitem[Lib10]{LibRA}
N.~Libedinsky.
\newblock Presentation of right-angled {S}oergel categories by generators and
  relations.
\newblock {\em J. Pure Appl. Algebra}, 214(12):2265--2278, 2010.

\bibitem[Lib15]{LLC}
N.~Libedinsky.
\newblock Light leaves and {L}usztig's conjecture.
\newblock {\em Adv. Math.}, 280:772--807, 2015.

\bibitem[Lib17]{Li17}
N.~Libedinsky.
\newblock Gentle introduction to {S}oergel bimodules {I}: the basics.
\newblock Preprint, arXiv:1702.00039, 2017.

\bibitem[LP07]{lenart2007affine}
C.~Lenart and A.~Postnikov.
\newblock Affine {W}eyl groups in {K}-theory and representation theory.
\newblock {\em International Mathematics Research Notices}, 2007, 2007.

\bibitem[LRH18]{lobos2018graded}
D.~Lobos and S.~Ryom-Hansen.
\newblock Graded cellular basis and jucys-murphy elements for generalized blob
  algebras.
\newblock {\em arXiv preprint arXiv:1812.11143. To appear in Journal of pure
  and applied algebra.}, 2018.

\bibitem[LW17a]{libedinsky2017anti}
N.~Libedinsky and G.~Williamson.
\newblock The anti-spherical category.
\newblock {\em arXiv preprint arXiv:1702.00459}, 2017.

\bibitem[LW17b]{LW17}
N.~Libedinsky and G.~Williamson.
\newblock A non-perverse {S}oergel bimodule in type {$A$}.
\newblock {\em C. R. Math. Acad. Sci. Paris}, 355(8):853--858, 2017.

\bibitem[LW18]{lusztig2018billiards}
G.~Lusztig and G.~Williamson.
\newblock Billiards and tilting characters for $sl_3$.
\newblock {\em SIGMA}, 14(015), 2018.

\bibitem[Mar08]{Martin4}
P.~Martin.
\newblock On diagram categories, representation theory and statistical
  mechanics.
\newblock In {\em Noncommutative rings, group rings, diagram algebras and their
  applications}, volume 456 of {\em Contemp. Math.}, pages 99--136. Amer. Math.
  Soc., Providence, RI, 2008.

\bibitem[MS94]{martin1994blob}
P.~Martin and H.~Saleur.
\newblock The blob algebra and the periodic {T}emperley-{L}ieb algebra.
\newblock {\em Letters in mathematical physics}, 30(3):189--206, 1994.

\bibitem[MW00]{martin2000structure}
P.~P. Martin and D.~Woodcock.
\newblock On the structure of the blob algebra.
\newblock {\em Journal of Algebra}, 225(2):957--988, 2000.

\bibitem[MW03]{martin2003generalized}
P.~P. Martin and D.~Woodcock.
\newblock Generalized blob algebras and alcove geometry.
\newblock {\em LMS Journal of Computation and Mathematics}, 6:249--296, 2003.

\bibitem[Pla13]{plaza2013graded}
D.~Plaza.
\newblock Graded decomposition numbers for the blob algebra.
\newblock {\em Journal of Algebra}, 394:182--206, 2013.

\bibitem[PRH14]{plaza2014graded}
D.~Plaza and S.~Ryom-Hansen.
\newblock Graded cellular bases for {T}emperley--{L}ieb algebras of type {A}
  and {B}.
\newblock {\em Journal of Algebraic Combinatorics}, 40(1):137--177, 2014.

\bibitem[RH10]{ryom2010ariki}
S.~Ryom-Hansen.
\newblock The {A}riki--{T}erasoma--{Y}amada tensor space and the blob algebra.
\newblock {\em Journal of Algebra}, 324(10):2658--2675, 2010.

\bibitem[Rou08]{rouquier20082}
R.~Rouquier.
\newblock 2-{K}ac-{M}oody algebras.
\newblock {\em arXiv preprint arXiv:0812.5023}, 2008.

\bibitem[RW15]{RW}
S.~Riche and G.~Williamson.
\newblock Tilting modules and the $p$-canonical basis.
\newblock Preprint, arXiv:1512.08296, 2015.

\bibitem[Soe92]{SHC}
W.~Soergel.
\newblock The combinatorics of {H}arish-{C}handra bimodules.
\newblock {\em J. Reine Angew. Math.}, 429:49--74, 1992.

\bibitem[Soe97]{SoeKL}
W.~Soergel.
\newblock Kazhdan-{L}usztig polynomials and a combinatoric[s] for tilting
  modules.
\newblock {\em Represent. Theory}, 1:83--114 (electronic), 1997.

\bibitem[Soe07]{soergel2007kazhdan}
W.~Soergel.
\newblock {K}azhdan-{L}usztig-polynome und unzerlegbare bimoduln {\"u}ber
  polynomringen.
\newblock {\em Journal of the {I}nstitute of {M}athematics of {J}ussieu},
  6(3):501--525, 2007.

\bibitem[TL71]{temperley1971relations}
H.~N. Temperley and E.~H. Lieb.
\newblock Relations between the `percolation' and `colouring' problem and other
  graph-theoretical problems associated with regular planar lattices: some
  exact results for the `percolation' problem.
\newblock In {\em Proceedings of the Royal Society of London A: Mathematical,
  Physical and Engineering Sciences}, volume 322, pages 251--280. The Royal
  Society, 1971.

\bibitem[Wil12]{Wi12}
G.~Williamson.
\newblock Some examples of parity sheaves.
\newblock Oberwolfach report for the meeting ``Enveloping Algebras and
  Geometric Representation Theory'', available under:
  http://people.mpimbonn.mpg.de/geordie/GWilliamson-EnvAlg.pdf, 2012.

\bibitem[Wil17a]{williamson2017algebraic}
G.~Williamson.
\newblock Algebraic representations and constructible sheaves.
\newblock {\em Japanese Journal of Mathematics}, 12(2):211--259, 2017.

\bibitem[Wil17b]{WExplosion}
G.~Williamson.
\newblock Schubert calculus and torsion explosion.
\newblock {\em J. Amer. Math. Soc.}, 30(4):1023--1046, 2017.

\end{thebibliography}

\end{document}